\documentstyle[amstex,amssymb,11pt]{article}

  \input xypic
  \input cyracc.def

  \newfam\cyifam
  \font\tencyi=wncyi10
  \font\sevencyi=wncyi7
  \font\fivecyi=wncyi5
  \def\cyi{\fam\cyifam\tencyi\cyracc}
  \textfont\cyifam=\tencyi \scriptfont\cyifam=\sevencyi
    \scriptscriptfont\cyifam=\fivecyi

  \newfam\cyrfam
  \font\tencyr=wncyr10
  \font\sevencyr=wncyr7
  \font\fivecyr=wncyr5
  \def\cyr{\fam\cyrfam\tencyr\cyracc}
  \textfont\cyrfam=\tencyr \scriptfont\cyrfam=\sevencyr
    \scriptscriptfont\cyrfam=\fivecyr
 \newcommand{\lon}{\longrightarrow}
 
 \newcommand{\rar}{\rightarrow}
 \newcommand{\hook}{\hookrightarrow}
 \newcommand{\Proof}{{\bf Proof}.\, }

 \newcommand{\Def}{\mathsf{Def}}

 \newcommand{\Z}{{\Bbb Z}}
 \newcommand{\bS}{{\Bbb S}}
 \newcommand{\p}{{\partial}}

 \newcommand{\R}{{\Bbb R}}
 \newcommand{\ot}{\otimes}
 \newcommand{\tl}{\tilde}
 \newcommand{\tlk}{\tilde{\kappa}}
 
 \newcommand{\Beq}{\begin{equation}}
 \newcommand{\Eeq}{\end{equation}}
 \newcommand{\Beqr}{\begin{eqnarray}}
 \newcommand{\Eeqr}{\end{eqnarray}}
 \newcommand{\Beqrn}{\begin{eqnarray*}}
 \newcommand{\Eeqrn}{\end{eqnarray*}}
 \newcommand{\Ba}{\begin{array}}
 \newcommand{\Ea}{\end{array}}
 \newcommand{\Bi}{\begin{itemize}}
 \newcommand{\Ei}{\end{itemize}}
 \newcommand{\Bc}{\begin{center}}
 \newcommand{\Ec}{\end{center}}
 \newcommand{\fg}{{\frak g}}

 \newcommand{\f}{{\cal O}}

 \newcommand{\cA}{{\cal A}}
 \newcommand{\cB}{{\cal B}}
 \newcommand{\cC}{{\cal C}}
 \newcommand{\caD}{{\cal D}}
 \newcommand{\cE}{{\cal E}}
 
 \newcommand{\cG}{{\cal G}}
 \newcommand{\caH}{{\cal H}}
 \newcommand{\cI}{{\cal I}}
 
 \newcommand{\caL}{{\cal L}}
 \newcommand{\cM}{{\cal M}}
 \newcommand{\cP}{{\cal P}}

 \newcommand{\al}{\alpha}
 \newcommand{\be}{\beta}
 \newcommand{\ga}{\gamma}
 \newcommand{\ka}{\kappa}

 \newcommand{\var}{\varepsilon}

 \newcommand{\tv}{\tilde{v}}

 \newcommand{\bp}{\bar{\partial}}

 \newcommand{\Ker}{{\mathsf Ker}\, }
 \newcommand{\Img}{{\mathsf Im}\, }
 
 \newcommand{\Hom}{{\mathrm Hom}}

 \newcommand{\cT}{{\cal T}}

 \newcommand{\sip}{\smallskip}
 \newcommand{\bip}{\bigskip}

\begin{document}

 \sloppy

 \title{Operads, deformation theory\\
 and  $F$-manifolds}
 \author{ S.A.\ Merkulov}
 \date{}
 \maketitle

 {\hfill {\cyi Posvyawaet\-sya \i.~I.~Maninu }}

 {\hfill {\cyi s blagodarnost{\cprime}{yu} i  vos{h}iweniem. }}






 \begin{center}
 {\bf \S 0. Introduction}
 \end{center}
 {\bf 0.1. Little disks operad and Hertling-Manin's $F$-manifolds.}
 Frobenius manifolds created by Dubrovin in 1991 from rich theoretical
 physics material have been found since in many different fragments of
 mathematics
 --- quantum cohomology and mirror symmetry, complex geometry,
  symplectic geometry, singularity theory,
 integrable systems --- raising hopes for  unifying them into one picture. It
 also became clear that the notion of Frobenius manifold is not broad enough to
 cover {\em all}\, objects of the associated working categories; say, on the
 $B$-side of the mirror symmetry it applies only to extended moduli spaces of
 {\em Calabi-Yau}\, manifolds, the latter forming a rather small subcategory of the
 category of complex manifolds. In 1998 Hertling and Manin \cite{HM} introduced
 a weaker notion of $F$-{\em manifold} which is, by definition, a pair
 $(M,\mu_2)$ consisting of a smooth supermanifold $M$ and a
  smooth
 $\f_M$-linear associative graded commutative multiplication on the
 tangent sheaf, $\mu_2: \ot^2_{\f_M} \cT_M \rar \cT_M$, satisfying
 the  integrability condition,
 $$
 [\mu_2,\mu_2]=0,
 $$
 where $[\mu_2,\mu_2]: \ot^4_{\f_M} \cT_M \rar \cT_M$ is given
 explicitly by
 \Beqrn
  [\mu_2,\mu_2](X,Y,Z,W)&:=&
  [\mu_2(X, Y), \mu_2(Z,W)] - \mu_2([\mu_2(X, Y), Z], W)\\
  &&-(-1)^{(|X|+|Y|)|Z|} \mu_2(Z,[\mu_2(X,Y), W])
   \\
 && - \mu_2(X, [Y,\mu_2(Z, W)]) - (-1)^{|Y|(|Z|+|W|)} \mu_2[X,
 \mu_2(Z,W)], Y) \\
 && + (-1)^{|Y||Z|} \mu_2(X, \mu_2(Z, [Y,W]))  + \mu_2(X, \mu_2([Y,Z], W)) \\
 && + (-1)^{|Y||Z|}\mu_2([X,Z],\mu_2( Y, W)) +
 (-1)^{|W|(|Y|+|Z|)}\mu_2([X,W],\mu_2(Y, Z)).
  \Eeqrn
  A  non-trivial part of the above definition is an implicit assertion  that
  $[\mu_2,\mu_2]$ is a tensor, i.e.\ $\f_M$-polylinear in all
  four inputs. It is here where the assumption that $\mu_2$ is
  both graded commutative an associative plays a key role.

  \sip

  Any Frobenius manifold is an $F$-manifold. Any $F$-manifold with
 semi-simple product $\mu_2$ can be made into a Frobenius manifold
 \cite{HM}. Hertling in his book \cite{H} explained in detail how
 $F$-manifolds turn up in the singularity theory.

 \sip

 In this paper we show that (cohomology/strong homotopy, see below)
 $F$-manifolds arise naturally in every mathematical structure which
 admits an action of the chain operad of the little disks operad (or
 its more compact version, $\cG_{\infty}$-operad \cite{GJ}). In particular, we
 prove

 \bip

 {\bf Theorem A}. {\em Let $\mathsf A$ be either a complex or symplectic
 structure on a compact manifold. Then the smooth part, $M_{\mathrm reg}$,
  of the extended moduli space $M$
 of deformations of \, $\mathsf A$ is canonically an $F$-manifold.}

 \bip

 {\bf 0.2. Cohomology $F$-manifolds.} It is often  not a pleasure  to
 work with objects like $M_{\mathrm reg}\subset M$ in Theorem~A;
 moreover, their existence is not guaranteed for many reasonable
 deformation problems.

 \sip

  The germ, $(M,*\,)$, of, in general, singular
  moduli space $M$ (if it exists at all) at the distinguished point $*$
  always admits a smooth dg resolution, $(\cM, *, \eth)$ \cite{Me2}.
 The latter  consists of
  a germ of a {\em smooth}\, graded  pointed manifold, $(\cM, *\, )$, and a
  germ of a {\em smooth}\,
 degree 1 vector field, $\eth$, satisfying two conditions,
 $$
 [\eth,\eth]=0, \ \ \ {\mathrm and}\ \ \eth I\subset I^2,
 $$
 where $I$ is the ideal of the distinguished point. The relation
 between $(M,*\,)$ (which may not exist as an analytic space) and
 $(\cM,*, \eth)$ (which always exists for any deformation problem!) is
 given by the well known formula of ``nonlinear cohomology'',
 $$
 M\simeq \frac{{\mathsf Zeros}(\eth)}{\Img \eth},
 $$
 representing $M$ as the quotient of the zero set of
 the vector field $\eth$ by
 the integrable distribution
 $$
 \Img \eth:= \{X\in \cT_{\cM}: X=[\eth, Y]\ {\mathrm for\ some}\ Y\in \cT_{\cM}\},
 $$
 which, as it is easy to check, is tangent to ${\mathsf Zeros}(\eth)$.

 \sip

 With any dg manifold $(\cM,\eth)$ one can associate two
 cohomology sheaves: the  cohomology structure sheaf,
 $$
 H(\f_{\cM}) :=
 \frac{\Ker \eth: \f_{\cM} \rar \f_{\cM}}
 {\Img \eth: \f_{\cM} \rar \f_{\cM}}\, ,
 $$
 and the  cohomology tangent sheaf,
 $$
 \caH\cT_{\cM}:= \frac{\Ker Lie_{\eth}: \cT_{\cM} \rar \cT_{\cM}}
 {\Img Lie_{\eth}: \cT_{\cM} \rar \cT_{\cM}}\, ,
 $$
 which is a sheaf of $H(\f_{\cM})$-modules (in fact, a sheaf of Lie
 $H(\f_{\cM})$-algebras).

   \sip

 We define a {\em cohomology $F$-manifold}\, to be a dg manifold
 $(\cM, \eth)$ together with a graded commutative associative
 $H(\f_{\cM})$-polylinear product $\mu_2 : \caH\cT_{\cM}\times
 \caH\cT_{\cM} \rar \caH\cT_{\cM}$, such that the integrability
 condition, $[\mu_2,\mu_2]=0$, holds. This notion also makes sense in
 the category of formal dg manifolds.

 \bip

 {\bf Theorem B}. {\em If the operad $\cG_{\infty}$ acts on a dg vector space $(V,d)$,
 then the formal graded manifold associated with the cohomology vector space $H(V,d)$ is canonically
 a cohomology $F$-manifold.}

 \bip

 Let us emphasize again that the notion of (cohomology)
 $F$-manifold is diffeomorphism invariant. Though the input in
 Theorem B  belongs to the category of vector spaces which
 one can geometrically interpret as pointed affine (=flat)
 manifolds, the output lies in the category of general smooth graded
 manifolds with morphisms being arbitrary (not necessary, linear)
 smooth maps\footnote{There is no way to remember the original flat
 structure of the input unless the combination $(G_{\infty} \
 {\mathrm action}, V,d)$ is formal as a $L_{\infty}$-algebra, and
 one makes a particular choice of a homotopy class of formality
 maps.}. Thus the output  of Theorem B belongs  to the realm of
 differential geometry.

 \bip

 Recent proofs of Deligne's conjecture \cite{Ko3,KS0,MS,Ta,V} together
 with Theorem~B
 imply\footnote{Here and everywhere in this paper  $k$ stands for a field of
 characteristic $0$. Every vector space is implicitly assumed to be over $k$.}

 \sip

 {\bf Corollary C}. (i) {\em Let $A$ be an associative $k$-algebra.
 The formal manifold
 associated with the Hochschild cohomology $H^{\bullet}(A,A)$ is naturally a
 cohomology $F$-manifold.}

 \sip

 (ii) {\em Let $X$ be a compact topological space. The formal manifold
 associated with its singular cohomology $H^{\bullet}(X,k)$ is
 naturally a cohomology $F$-manifold.}

 \bip

 {\bf 0.3. $F_{\infty}$-manifolds.} Instead of passing to cohomology
 sheaves as above, one can adopt the notion of $F$-manifold to the
 category of dg manifolds by constructing its strong homotopy version.
 We do it in this paper with the help of the $\cG_{\infty}$-operad
 (cf.\ Theorem~B).

 \sip

 Let $(\cM,\eth,*)$ be a formal dg manifold and let
 $$
 \mu_{\bullet}=\{\mu_n\}_{n\geq 1}: \ot^{\bullet}_{\f_{\cM}}\cT_{\cM}
 \rar \cT_{\cM}
 $$
  be a structure of $\cC_{\infty}$-algebra on the tangent sheaf. We call
  it {\em geometric}\, if $\mu_1=Lie_{\eth}$ and  $\mu_{\bullet\geq 2}$
  are morphisms of $\f_\cM$-modules, i.e.\ are tensors.
   If all $\mu_n$ except $\mu_2$ vanish, this
  structure reduces to the structure of graded commutative associative
  product as in Sect.\ 0.1.
Note that  $(\ot^{\bullet}_{\f_{\cM}}\cT^*_{\cM}\ot_{\f_{\cM}} \cT_{\cM},
 Lie_{\eth})$ is a complex of (sheaves of) $\f_\cM$-modules.
 Its cohomology is denoted by
 $H(\ot^{\bullet}_{\f_{\cM}}\cT^*_{\cM}\ot_{\f_{\cM}} \cT_{\cM})$.

  \sip

Choosing a torsion-free affine connection $\nabla$ on $\cM$, one
can construct
 an extension of the Hertling-Manin's ``bracket''  $[\mu_2,\mu_2]$
   to geometric
  $\cC_{\infty}$-structures,
  $$
 [\mu_{\bullet},\mu_{\bullet}]^\nabla(X_1,\ldots,X_{\bullet},Y_1,\ldots,
 Y_{\bullet}):=[\mu_{\bullet}(X_1,\ldots,X_{\bullet}),
 \mu_{\bullet}(Y_1,\ldots,Y_{\bullet})]
 + {\mathrm correction\ terms},
  $$
 producing thereby a collection of {\em tensors}\footnote{The
 assumption that $\mu_{\bullet}$ is a $\cC_{\infty}$-structure is
 important. The construction does not work for
 geometric $\cA_{\infty}$-structures.},
 $[\mu_{\bullet},\mu_{\bullet}]^\nabla:\ot^{\bullet+\bullet}_{\f_{\cM}}\cT_{\cM}
 \rar \cT_{\cM}$, satisfying the following two conditions,
\begin{itemize}
    \item $Lie_{\eth}[\mu_{\bullet},\mu_{\bullet}]^\nabla =0$,
    \item the cohomology class,
$$
[[\mu_{\bullet},\mu_{\bullet}]]\in H(\ot^{\bullet}_{\f_{\cM}}\cT^*_{\cM}\ot_{\f_{\cM}}
 \cT_{\cM}),
$$
    produced by  $[\mu_{\bullet},\mu_{\bullet}]^\nabla$ does not depend on the choice
    of the connection $\nabla$ and hence gives a well-defined invariant of the geometric
    $\cC_\infty$-structure. Moreover, this invariant depends only on the homotopy class of that
    structure.
\end{itemize}

 The correction terms to $[[\mu_{\bullet},\mu_{\bullet}]]$ can, in principle, be
 read off from the
 structural equations of the $\cG_{\infty}$-operad, as explained in
 Sect.\ 4. However, all the basic properties of the bracket
 $[\mu_{\bullet},\mu_{\bullet}]^\nabla$,
 such as its existence, $\f_{\cM}$-linearity, $Lie_{\eth}$-closedness etc.,
  can be proved without doing this sort of explicit calculations.

 \sip


 \bip

 {\bf  Definition D.} {\em An $F_{\infty}$-manifold is a dg manifold,
 $(\cM,\eth,*)$, together with a  homotopy class of geometric
   $\cC_{\infty}$-structures,
  $\{\mu_{\bullet}: \ot^{\bullet}_{\f_{\cM}}\cT_{\cM}
 \rar \cT_{\cM}\}$, satisfying the integrability condition,
  $$
  [[\mu_{\bullet},\mu_{\bullet}]]=0.
 $$
 }

 Clearly, any $F_{\infty}$-manifold gives naturally rise to a
 cohomology $F$-manifold.  In fact, the cohomology $F$-manifolds
 discussed above in Sect.\ 0.2 are precisely of this type:

 \bip

 {\bf Theorem E}. {\em All statements of Theorem~B and Corollary~C
 remain true if one replaces
 \begin{center}
  cohomology  $F$-manifold \ \ \ $\lon$ \ \ \ $F_{\infty}$-manifold.
 \end{center}
 }
  \bip

 {\bf 0.4. Content.}
 In \S 1 and \S 2 we remind basic notions and notations of the (homotopy) theory of operads
 and  discuss in detail some particular examples. The main result in  \S 3
 is an explicit graphical description, Proposition~3.6.1,
  of the cobar construction for the operad
 of non-commutative Gerstenhaber algebras and a surprisingly nice geometric
 interpretation, Theorem 3.9.2, of the derived category of algebras over that
 operad. In \S 4 we outline an operadic guide to the extended deformation theory
 (as a more informative alternative to the classical idea of deformation functor)
 and, in that context, prove all the claims made in the Introduction.

 \bip

 \bip
 \begin{center}
 {\bf \S 1.  Operads and their  algebras}
 \end{center}

 {\bf 1.1. Operads.} By an operad in this paper we always
 understand what is usually called a {\em nonunital or
 pseudo-operad} \cite{Mar}, that is, a pair of collections,
 $$
 \f=\left( \left\{\f(n)\right\}_{n\geq 1}\
 ,\  \{\circ^{n,n'}_{i}\}_{ n,n'\geq 1 \atop 1\leq i\leq n}\right),
 $$
 where each $\f(n)$ is a $\Z$-graded vector space
 equipped with a linear action of the
 permutation group $\bS_n$ (the collection $\left\{\f(n)\right\}_{n\geq 1}$
 will sometimes be called an $\Bbb S$-{\em module}),
 and each $\circ^{n,n'}_{i}$ is a linear equivariant
 map,
 $$
 \circ_i^{n,n'}: \f(n)\ot \f(n') \lon \f(n+n'-1),
 $$
 such that, for any $f\in \f(n)$, $f'\in \f(n')$ and $f''\in \f(n'')$, one has
 $$
 \left(f \circ_i^{n,n'} f'\right) \circ_{j+n'-1}^{n+n'-1,n''} f'' =
 (-1)^{|f'||f''|}
 \left(f \circ_j^{n,n''} f''\right) \circ_i^{n+n''-1,n'} f', \ \ \forall\  1\leq i < j\leq n,
 $$
 and
 $$
 f \circ_i^{n,n'+n''-1} \left( f'\circ_j^{n',n''} f''\right) =
 \left(f \circ_i^{n,n'} f'\right) \circ_{i+j-1}^{n+n'-1,n''} f'', \ \ \forall\  1\leq i\leq n,
 1\leq j\leq n'.
 $$

 Equivariance of $\circ_i^{n,n'}$ above means that for any $\sigma\in \bS_n$ and
 $\sigma'\in \bS_{n'}$
 one has
 $$
 (\sigma f)\circ_i^{n,n'}(\sigma' f')
 = (\sigma_i \sigma')
 ( f\circ_i^{n,n'} f')
 $$
 where $(\sigma_i \sigma')\in \bS_{n+n'-1}$ is given by inserting the permutation $\sigma'$
 into the $i$th place of $\sigma$.

 \sip

 An {\em ideal}\, in an operad $\f$ is a collection $I$
 of $\bS_n$-invariant
 subspaces $\{\cI(n)\subset \f(n)\}_{n\geq 1}$
 such that $f\circ_i^{n,n'} f'\in \cI(n+n'-1)$
 whenever $f\in \cI(n)$ or $f'\in \cI(n')$; in particular, $\cI$ is a suboperad of $\f$.
 It is clear that the quotient ${\Bbb S}$-module
 $\{\f(n)/ \cI(n)\}_{n\geq 1}$ has a naturally induced structure of an operad called a
 {\em quotient operad}.

 \sip

 An operad $\f$ with $\f(1)=0$ is called {\em simply connected}.

 \bip
 \sip

 {\bf 1.2. Free operads and trees.}
 A morphism of operads, $f:\f \rar \f'$, is, by definition, a morphism of the associated
 ${\Bbb S}$-modules, $\{f(n): \f(n)\lon \f'(n)\}_{n\geq 2}$, which commutes in the obvious way
 with all the operations $\circ_{i}^{n,n'}$. Operads form a category.
 \sip

 The forgetful functor
 $$
 \Ba{ccc}
 {\mathrm Category\ of\ operads} & \lon & {\mathrm Category\ of}\ {\Bbb S}\hspace{-1mm}-
 \hspace{-1mm}{\mathrm modules} \\
 \left( \left\{\f(n)\right\}_{n\geq 1}\
 ,\  \{\circ^{n,n'}_{i}\}_{ n,n'\geq 1 \atop 1\leq i\leq n}\right) & \lon &
 \left( \left\{\f(n)\right\}_{n\geq 1}\right),
 \Ea
 $$
 has a left adjoint functor,
  $Free$, which associates to an arbitrary collection,
 $\cE=\{\cE(n)\}_{n\geq 1}$,
 of graded vector $\bS_n$-spaces  the {\em free operad}, $Free(\cE)$. It is best described in terms
 of trees as follows (see \cite{GK,GJ, KS0} for more details).

 \sip

 An {\em $[n]$-tree}\, $T$ is, by definition, the data $(V_T, N_T, \phi_T)$
 consisting of
 \Bi
 \item a stratified finite set $V_T= V^i_T \sqcup V^t_T$ whose elements are called
 {\em vertices};
  elements of the subset $V^i_T$ (resp.\ $V^t_T$) are called {\em internal}\, (resp.\ {\em tail})
 { vertices};
 \item a bijection $\phi: V^t_T \rar \{1,2,\ldots, n\}=:[n]$;
 \item a map $N_T: V_T \rar V_T$ satisfying the conditions: (i) $N_T$ has only one fixed point
 $root_T$ which lies in $V^i_T$ and is called the {\em root vertex}, (ii) $N_T^k(v)=root_T$,
 $\forall\ v\in V_T$ and $k\gg 1$, (iii) for all $v\in V_T^i$ the cardinality, $\# v$,
 of the set $N_T^{-1}(v)$ is greater than or equal to $1$, while for all $v\in V_T^t$ one has
  $\# v=0$.
 \Ei

 The number $\# v$ is often called the {\em valency}\, of the vertex $v$; the pairs
 $(v, N_T(v))$ are called edges.

 \sip

 Given an $\bS$-module $\cE=\{\cE(n)\}_{n\geq 1}$, we
 can associate to an $[n]$-tree $T$ the vector space
 $$
 \cE(T):= \bigotimes_{v\in V_T^i} \cE(\# v).
 $$
 Its elements are interpreted as $[n]$-trees
 whose internal vertices are decorated with elements of $\cE$.
 The permutation group $\bS_n$ then acts on this space via
 relabelling the tail vertices (i.e\ changing
 $\phi_T$ to $\sigma\circ \phi_T$, $\sigma\in \bS_n$).

 \sip

 Now, as an $\bS$-module the free operad $Free(\cE)$ is defined as
 $$
 Free(\cE)(n) = \bigoplus_{[[n]-{\mathrm trees}\ T]} \cE(T),
 $$
 where the summation goes over all isomorphism classes of $[n]$-trees. The composition,
 say $f\circ_i^{n,n'} f'$, is given by gluing the root vertex of the decorated $[n']$-tree
 $f'\in Free(\cE)(n')$  with the $i$-labelled
 tail vertex of the decorated $[n]$-tree $f$. The new numeration, $\phi: V_T^t\rar [n+n'-1]$,
  of tails is clear.

 \sip

 Any free operad is naturally graded, $Free(\cE)=\bigoplus_{p=1}^{\infty} Free^p(\cE)$,
 where $Free^p(\cE)$ is the $\bS$-submodule of $Free(\cE)$
 spanned by all possible isomorphism classes of $\cE$-decorated trees with precisely $p$ internal
 vertices.

 \bip
 \sip

 {\bf 1.3. Example}.
  Let $V$ be a $\Z$-graded vector space. The associated ${\Bbb S}$-module,
 $$
 \cE_V=\{ \cE_V(n):= Hom(V^{\ot n}, V)\},
 $$ has a natural structure of operad with compositions,
 $f\circ_i^{n,n'} f'$,  given by inserting the output of $f'\in Hom(V^{\ot n'},V)$
 into the $i$-th input of $f\in Hom(V^{\ot n},V)$.

 \sip

 An {\em algebra over an operad}\,  $\f$ is, by definition, a $\Z$-graded vector space $V$
 together with
 a morphism of operads $\f\rar \cE_V$.

 \bip

 {\bf 1.4. Example}. Let $\cA$ be an $\bS$-module given by
 $$
 \cA(n):= \left\{ \Ba{ll}
 k[\bS_2][0] & {\mathrm if}\ n=2\\
 0 & {\mathrm otherwise},
 \Ea
 \right.
 $$
 where here and  below the symbol $k[\bS_n][p]$  stands for the graded vector space whose
 only  non-vanishing homogeneous component lies in degree $-p$ and equals the regular
 representation $k[\bS_n]$ of
 the permutation group $\bS_n$\footnote{More generally, for a $\Z$-graded vector space
 $V=\oplus_{i\in \Z} V^i$, the symbol $V[p]$ stands for the $\Z$-graded vector space with
 $V[p]^i:= v^{i+p}$.}.
  If we identify the natural basis, $id$ and $(12)$, of  $k[\bS_2]$
 with {\em planar}\,  $[2]$-corollas,
 $$
 
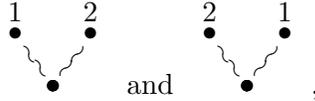
\begin{xy}
 <5mm,0cm>*{\bullet};<0cm,7mm>*{\bullet}**@{~},
 <5mm,0cm>*{\bullet};<10mm,7mm>*{\bullet}**@{~},
 <5mm,0cm>*{\bullet};<0mm,10mm>*{1}**@{},
 <5mm,0cm>*{\bullet};<10mm,10mm>*{2}**@{},
 \end{xy}
 \ \ \ {\mathrm and}\ \ \
 \begin{xy}
 <5mm,0cm>*{\bullet};<0cm,7mm>*{\bullet}**@{~},
 <5mm,0cm>*{\bullet};<10mm,7mm>*{\bullet}**@{~},
 <5mm,0cm>*{\bullet};<0mm,10mm>*{2}**@{},
 <5mm,0cm>*{\bullet};<10mm,10mm>*{1}**@{},
 \end{xy} \ \ ,
 $$
 then the associated free operad $\{Free(\cA)(n), \circ_i^{n,n'}\}$
 can be represented as a linear span
 of all possible (isomorphism classes of) binary {\em planar}\, $[n]$-trees, e.g.
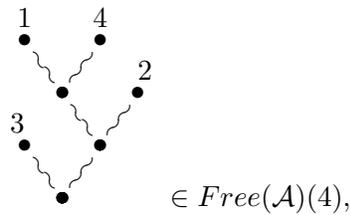
 $$
 \begin{xy}
 <5mm,0mm>*{\bullet};<0cm,7mm>*{\bullet}**@{~},
 <5mm,0cm>*{\bullet};<10mm,7mm>*{\bullet}**@{~},
 <10mm,7mm>*{\bullet};<5mm,14mm>*{\bullet}**@{~},
 <10mm,7mm>*{\bullet};<15mm,14mm>*{\bullet}**@{~},
 <5mm,14mm>*{\bullet};<0cm,21mm>*{\bullet}**@{~},
 <5mm,14mm>*{\bullet};<10mm,21mm>*{\bullet}**@{~},
 <5mm,0cm>*{\bullet};<-1mm,10mm>*{3}**@{},
 <5mm,0cm>*{\bullet};<0mm,24mm>*{1}**@{},
 <5mm,0cm>*{\bullet};<10mm,24mm>*{4}**@{},
 <5mm,0cm>*{\bullet};<16mm,17mm>*{2}**@{},
 \end{xy}
 \ \in Free(\cA)(4),
 $$
 with the compositions $\circ_i^{n,n'}$
 given simply by gluing  the root vertex of a planar $[n']$-tree to the $i$th tail vertex
 of an $[n]$-tree (the new numeration of tails is clear).
 Indeed an isomorphism class of an $\{id,(12)\}$-decorated abstract($\equiv$space) binary
 tree of
 Subsect.\ 1.2  has a natural representative which lies in a fixed plane in the  space
 and which is consistent with the interpretation of the labelling set
 $\{id,(12)\}$ as the set of planar
 [2]-corollas; more importantly, the resulting correspondence
 $$
 \left\{\Ba{c}
 {\mathrm isomorphism\  classes\  of}\\
  \{id,(12)\}\hspace{-1mm}-\hspace{-1mm}{\mathrm decorated\ abstract\ binary\ trees}
 \Ea \right\}
 \lon
 \left\{\Ba{c}
 {\mathrm isomorphism\  classes\  of}\\
  {\mathrm planar\ binary\ numbered\  trees}
 \Ea \right\}
 $$
 is one-to-one.
 \sip

 Algebras over $Free(\cA)$ are not that interesting objects --- they are
  just graded vector spaces $V$ together with  a fixed
  element of $\Hom(V^{\ot 2}, V)$ which can be arbitrary.

 \sip

 Let $\cI_{\cA}$ be the ideal in $Free(\cA)$ generated by 3! vectors of the form
 $$
 
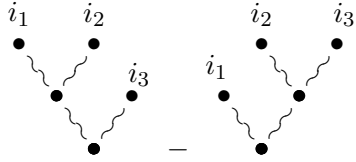
\begin{xy}
 <5mm,0mm>*{\bullet};<0cm,7mm>*{\bullet}**@{~},
 <5mm,0cm>*{\bullet};<10mm,7mm>*{\bullet}**@{~},
 <0mm,7mm>*{\bullet};<-5mm,14mm>*{\bullet}**@{~},
 <0mm,7mm>*{\bullet};<5mm,14mm>*{\bullet}**@{~},
 <5mm,0cm>*{\bullet};<-5mm,18mm>*{i_1}**@{},
 <5mm,0cm>*{\bullet};<5mm,18mm>*{i_2}**@{},
 <5mm,0cm>*{\bullet};<11mm,10mm>*{i_3}**@{},
 \end{xy}
 \
 -
 \
 \begin{xy}
 <5mm,0mm>*{\bullet};<0cm,7mm>*{\bullet}**@{~},
 <5mm,0cm>*{\bullet};<10mm,7mm>*{\bullet}**@{~},
 <10mm,7mm>*{\bullet};<5mm,14mm>*{\bullet}**@{~},
 <10mm,7mm>*{\bullet};<15mm,14mm>*{\bullet}**@{~},
 <5mm,0cm>*{\bullet};<-1mm,11mm>*{i_1}**@{},
 <5mm,0cm>*{\bullet};<5mm,18mm>*{i_2}**@{},
 <5mm,0cm>*{\bullet};<16mm,18mm>*{i_3}**@{},
 \end{xy}
 $$
 The quotient operad $Free(\cA)/\cI_{\cA}$ is denoted by $\cA ss$ for the obvious reason ---
 its algebras are nothing but the usual graded associative algebras. As an $\bS$-module,
 $\cA ss(n)\simeq k[\bS_n]$.

 \bip

 {\bf 1.5. Example}. Consider an  $\bS$-module,
 $$
 \cC omm(n):= {\mathbf 1}_{n}[0]
 $$
 where ${\mathbf 1}_{n}$ stands for the trivial representation of
 the permutation group $\bS_n$. This $\bS$-module
 can be made into an operad $\cC omm$ by defining
 the compositions $\circ_i^{n,n'}$ to be the identity maps. It is not hard to check
 that $\cC omm$-algebras are graded commutative associative algebras in the usual sense.

 \sip

 For later reference it will be convenient to represent the operad $\cC omm$ as a quotient of
 a  free operad. For this purpose we first consider an $\bS$-module $\cC$,
 $$
 \cC(n):= \left\{ \Ba{ll}
 {\mathbf 1}_{2}[0] & {\mathrm if}\ n=2\\
 0 & {\mathrm otherwise}.
 \Ea
 \right.
 $$
 If we identify a basis vector of ${\mathbf 1}_{n}[0]$ with the unique (up to an isomorphism)
 {\em space}\,  corolla (i.e.\ the one embedded in $\R^3$)
 $$
 
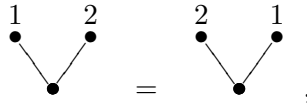
\begin{xy}
 <5mm,0cm>*{\bullet};<0cm,7mm>*{\bullet}**@{-},
 <5mm,0cm>*{\bullet};<10mm,7mm>*{\bullet}**@{-},
 <5mm,0cm>*{\bullet};<0mm,10mm>*{1}**@{},
 <5mm,0cm>*{\bullet};<10mm,10mm>*{2}**@{},
 \end{xy}
 \ \ \ = \ \ \
 \begin{xy}
 <5mm,0cm>*{\bullet};<0cm,7mm>*{\bullet}**@{-},
 <5mm,0cm>*{\bullet};<10mm,7mm>*{\bullet}**@{-},
 <5mm,0cm>*{\bullet};<0mm,10mm>*{2}**@{},
 <5mm,0cm>*{\bullet};<10mm,10mm>*{1}**@{},
 \end{xy} \ \ ,
 $$
 then the associated free operad $\{Free(\cC)(n)\}$ can be represented as a linear span
 of all possible isomorphism classes of binary {\em space}\,  $[n]$-trees; for example,
 $Free(\cC)(3)$ is a 3-dimensional vector space spanned by the following space $[3]$-trees
 $$
 
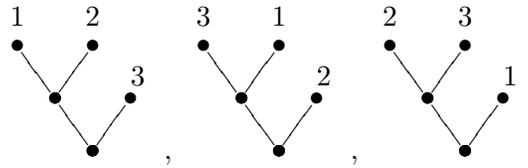
\begin{xy}
 <5mm,0mm>*{\bullet};<0cm,7mm>*{\bullet}**@{-},
 <5mm,0cm>*{\bullet};<10mm,7mm>*{\bullet}**@{-},
 <0mm,7mm>*{\bullet};<-5mm,14mm>*{\bullet}**@{-},
 <0mm,7mm>*{\bullet};<5mm,14mm>*{\bullet}**@{-},
 <5mm,0cm>*{\bullet};<-5mm,18mm>*{1}**@{},
 <5mm,0cm>*{\bullet};<5mm,18mm>*{2}**@{},
 <5mm,0cm>*{\bullet};<11mm,10mm>*{3}**@{},
 \end{xy}\ \
 , \ \
 \begin{xy}
 <5mm,0mm>*{\bullet};<0cm,7mm>*{\bullet}**@{-},
 <5mm,0cm>*{\bullet};<10mm,7mm>*{\bullet}**@{-},
 <0mm,7mm>*{\bullet};<-5mm,14mm>*{\bullet}**@{-},
 <0mm,7mm>*{\bullet};<5mm,14mm>*{\bullet}**@{-},
 <5mm,0cm>*{\bullet};<-5mm,18mm>*{3}**@{},
 <5mm,0cm>*{\bullet};<5mm,18mm>*{1}**@{},
 <5mm,0cm>*{\bullet};<11mm,10mm>*{2}**@{},
 \end{xy}\ \ , \ \
 \begin{xy}
 <5mm,0mm>*{\bullet};<0cm,7mm>*{\bullet}**@{-},
 <5mm,0cm>*{\bullet};<10mm,7mm>*{\bullet}**@{-},
 <0mm,7mm>*{\bullet};<-5mm,14mm>*{\bullet}**@{-},
 <0mm,7mm>*{\bullet};<5mm,14mm>*{\bullet}**@{-},
 <5mm,0cm>*{\bullet};<-5mm,18mm>*{2}**@{},
 <5mm,0cm>*{\bullet};<5mm,18mm>*{3}**@{},
 <5mm,0cm>*{\bullet};<11mm,10mm>*{1}**@{},
 \end{xy}\ .
 $$
 The composition in $Free(\cC)$ is given   by gluing  the root vertex of one space
 tree to a tail vertex of another one. The new numeration of tail vertices is clear.

 \sip

 Let $\cI_{\cC}$ be the ideal in $Free(\cC)$ generated by 2 vectors of the form
 $$
 
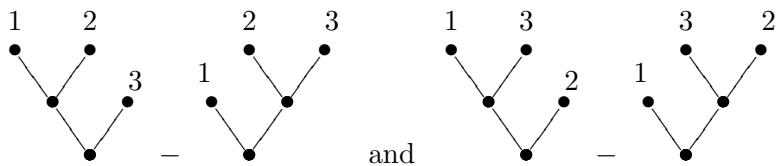
\begin{xy}
 <5mm,0mm>*{\bullet};<0cm,7mm>*{\bullet}**@{-},
 <5mm,0cm>*{\bullet};<10mm,7mm>*{\bullet}**@{-},
 <0mm,7mm>*{\bullet};<-5mm,14mm>*{\bullet}**@{-},
 <0mm,7mm>*{\bullet};<5mm,14mm>*{\bullet}**@{-},
 <5mm,0cm>*{\bullet};<-5mm,18mm>*{1}**@{},
 <5mm,0cm>*{\bullet};<5mm,18mm>*{2}**@{},
 <5mm,0cm>*{\bullet};<11mm,10mm>*{3}**@{},
 \end{xy}
 \
 -
 \
 \begin{xy}
 <5mm,0mm>*{\bullet};<0cm,7mm>*{\bullet}**@{-},
 <5mm,0cm>*{\bullet};<10mm,7mm>*{\bullet}**@{-},
 <10mm,7mm>*{\bullet};<5mm,14mm>*{\bullet}**@{-},
 <10mm,7mm>*{\bullet};<15mm,14mm>*{\bullet}**@{-},
 <5mm,0cm>*{\bullet};<-1mm,11mm>*{1}**@{},
 <5mm,0cm>*{\bullet};<5mm,18mm>*{2}**@{},
 <5mm,0cm>*{\bullet};<16mm,18mm>*{3}**@{},
 \end{xy}
 \ \ \ {\mathrm and} \ \ \
 \begin{xy}
 <5mm,0mm>*{\bullet};<0cm,7mm>*{\bullet}**@{-},
 <5mm,0cm>*{\bullet};<10mm,7mm>*{\bullet}**@{-},
 <0mm,7mm>*{\bullet};<-5mm,14mm>*{\bullet}**@{-},
 <0mm,7mm>*{\bullet};<5mm,14mm>*{\bullet}**@{-},
 <5mm,0cm>*{\bullet};<-5mm,18mm>*{1}**@{},
 <5mm,0cm>*{\bullet};<5mm,18mm>*{3}**@{},
 <5mm,0cm>*{\bullet};<11mm,10mm>*{2}**@{},
 \end{xy}
 \
 -
 \
 \begin{xy}
 <5mm,0mm>*{\bullet};<0cm,7mm>*{\bullet}**@{-},
 <5mm,0cm>*{\bullet};<10mm,7mm>*{\bullet}**@{-},
 <10mm,7mm>*{\bullet};<5mm,14mm>*{\bullet}**@{-},
 <10mm,7mm>*{\bullet};<15mm,14mm>*{\bullet}**@{-},
 <5mm,0cm>*{\bullet};<-1mm,11mm>*{1}**@{},
 <5mm,0cm>*{\bullet};<5mm,18mm>*{3}**@{},
 <5mm,0cm>*{\bullet};<16mm,18mm>*{2}**@{},
 \end{xy}
 $$
 The quotient operad $Free(\cC)/\cI_{\cC}$ is clearly isomorphic to $\cC omm$.

 \bip

 {\bf 1.6. Example}. Let $\caL$ be an $\bS$-module given by
 $$
 \caL(n):= \left\{ \Ba{ll}
 {\mathbf 1}_2[-1]
 & {\mathrm if}\ n=2,\\
 0 & {\mathrm otherwise}.
 \Ea
 \right.
 $$
 If we identify, as in Example 1.4,  a basis vector of the one dimensional
 vector space ${\mathbf 1}_2[1]$
 with the unique (up to isomorphism) {\em space}\, $[2]$-corolla

 $$
 
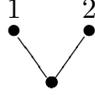
\begin{xy}
 <5mm,0cm>*{\bullet};<0cm,7mm>*{\bullet}**@{-},
 <5mm,0cm>*{\bullet};<10mm,7mm>*{\bullet}**@{-},
 <5mm,0cm>*{\bullet};<0mm,10mm>*{1}**@{},
 <5mm,0cm>*{\bullet};<10mm,10mm>*{2}**@{},
 \end{xy}
 $$
 then the associated free operad $\{Free(\caL)(n)\}$ can be represented as a linear span
 of all possible (isomorphism classes of) binary {\em space}\,  $[n]$-trees with
 the composition given by gluing  the root vertex of one space
 tree to a tail vertex
 of another one.

 \sip

 Let $\cI_{\caL}$ be the ideal in $Free(\caL)$ generated by the following $\bS_3$-invariant vector,
 $$
 
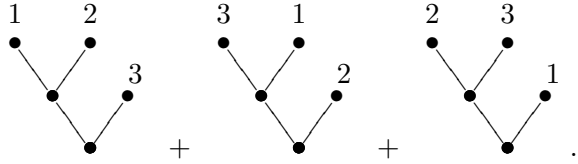
\begin{xy}
 <5mm,0mm>*{\bullet};<0cm,7mm>*{\bullet}**@{-},
 <5mm,0cm>*{\bullet};<10mm,7mm>*{\bullet}**@{-},
 <0mm,7mm>*{\bullet};<-5mm,14mm>*{\bullet}**@{-},
 <0mm,7mm>*{\bullet};<5mm,14mm>*{\bullet}**@{-},
 <5mm,0cm>*{\bullet};<-5mm,18mm>*{1}**@{},
 <5mm,0cm>*{\bullet};<5mm,18mm>*{2}**@{},
 <5mm,0cm>*{\bullet};<11mm,10mm>*{3}**@{},
 \end{xy}\ \ + \ \
 \begin{xy}
 <5mm,0mm>*{\bullet};<0cm,7mm>*{\bullet}**@{-},
 <5mm,0cm>*{\bullet};<10mm,7mm>*{\bullet}**@{-},
 <0mm,7mm>*{\bullet};<-5mm,14mm>*{\bullet}**@{-},
 <0mm,7mm>*{\bullet};<5mm,14mm>*{\bullet}**@{-},
 <5mm,0cm>*{\bullet};<-5mm,18mm>*{3}**@{},
 <5mm,0cm>*{\bullet};<5mm,18mm>*{1}**@{},
 <5mm,0cm>*{\bullet};<11mm,10mm>*{2}**@{},
 \end{xy}\ \ + \ \
 \begin{xy}
 <5mm,0mm>*{\bullet};<0cm,7mm>*{\bullet}**@{-},
 <5mm,0cm>*{\bullet};<10mm,7mm>*{\bullet}**@{-},
 <0mm,7mm>*{\bullet};<-5mm,14mm>*{\bullet}**@{-},
 <0mm,7mm>*{\bullet};<5mm,14mm>*{\bullet}**@{-},
 <5mm,0cm>*{\bullet};<-5mm,18mm>*{2}**@{},
 <5mm,0cm>*{\bullet};<5mm,18mm>*{3}**@{},
 <5mm,0cm>*{\bullet};<11mm,10mm>*{1}**@{},
 \end{xy}\ .
 $$

 Algebras over the associated quotient operad, $\caL ie:= Free(\caL)/ \cI_{\caL}$, are graded
 vector spaces $V$ equipped with a degree $-1$ element $\nu\in Hom(\odot^2 V, V)$
 satisfying the Jacobi condition,
 $$
 \nu\left(\nu(v_1,v_2), v_3\right) +
 (-1)^{|v_3|(|v_1|+|v_2|)} \nu\left(\nu(v_3,v_1), v_2\right) +
 (-1)^{|v_1|(|v_2|+|v_3|)} \nu\left(\nu(v_2,v_3), v_1\right)=0 .
 $$
 Setting
 $$
 [v_1\bullet v_2]:= (-1)^{|v_1|} \nu(v_1,v_2)
 $$
 we recover the notion of (odd) Lie algebra  \cite{Ma}.
 It is, of course, the same thing as the usual graded Lie algebra structure on the shifted
 graded vector space $V[1]$ but for our purposes it is more suitable {\em not}\, to make this shift;
 thus  in the present paper by a graded Lie algebra we always
 understand an algebra over the operad $\caL ie$, i.e. a
 pair $(V, [\ \bullet\ ])$ with $[\ \bullet \ ]: \odot^2 V \rar V$ having
 degree $-1$ and satisfying (odd) Jacobi identity.

 \bip

 {\bf 1.7. Example}. Let $\cA\caL$ be an $\bS$-module given by $\cA\caL(n):= \cA(n)\oplus \caL(n)$.
 Its only non-vanishing component $\cA\caL(2)$ is a 3-dimensional vector space spanned by
 two planar corollas in degree 0  and one space corolla in degree $-1$,
 $$
 
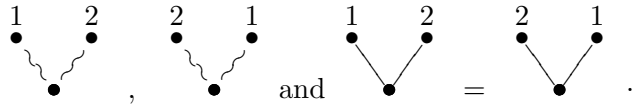
\begin{xy}
 <5mm,0cm>*{\bullet};<0cm,7mm>*{\bullet}**@{~},
 <5mm,0cm>*{\bullet};<10mm,7mm>*{\bullet}**@{~},
 <5mm,0cm>*{\bullet};<0mm,10mm>*{1}**@{},
 <5mm,0cm>*{\bullet};<10mm,10mm>*{2}**@{},
 \end{xy}
 \ \ \ , \ \ \
 \begin{xy}
 <5mm,0cm>*{\bullet};<0cm,7mm>*{\bullet}**@{~},
 <5mm,0cm>*{\bullet};<10mm,7mm>*{\bullet}**@{~},
 <5mm,0cm>*{\bullet};<0mm,10mm>*{2}**@{},
 <5mm,0cm>*{\bullet};<10mm,10mm>*{1}**@{},
 \end{xy} \ \  {\mathrm and} \ \
 \begin{xy}
 <5mm,0cm>*{\bullet};<0cm,7mm>*{\bullet}**@{-},
 <5mm,0cm>*{\bullet};<10mm,7mm>*{\bullet}**@{-},
 <5mm,0cm>*{\bullet};<0mm,10mm>*{1}**@{},
 <5mm,0cm>*{\bullet};<10mm,10mm>*{2}**@{},
 \end{xy}
 \ \ \ {=}\ \ \
 \begin{xy}
 <5mm,0cm>*{\bullet};<0cm,7mm>*{\bullet}**@{-},
 <5mm,0cm>*{\bullet};<10mm,7mm>*{\bullet}**@{-},
 <5mm,0cm>*{\bullet};<0mm,10mm>*{2}**@{},
 <5mm,0cm>*{\bullet};<10mm,10mm>*{1}**@{}, \ \ \hspace{10mm} .
 \end{xy}
 $$
 The associated free operad $Free(\cA\caL)$
 can be represented as a linear span
 of all possible isomorphism classes of  binary  $[n]$-trees in the 3-space $\R^3$
 with the condition that  all ``planar'' corollas are perpendicular to a fixed line in $\R^3$.
 The composition in $Free(\cA\caL)$ is given again  by gluing  the root vertex of one such partially
 planar/ partially space
 tree to a tail vertex of another one.

 \sip

 Let $\cI_{\cA\caL}$ be the ideal in $Free(\cA\caL)$ generated by the following
 3! vectors,
 $$
 
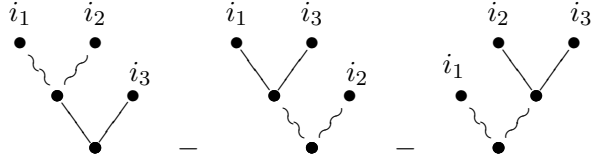
\begin{xy}
 <5mm,0mm>*{\bullet};<0cm,7mm>*{\bullet}**@{-},
 <5mm,0cm>*{\bullet};<10mm,7mm>*{\bullet}**@{-},
 <0mm,7mm>*{\bullet};<-5mm,14mm>*{\bullet}**@{~},
 <0mm,7mm>*{\bullet};<5mm,14mm>*{\bullet}**@{~},
 <5mm,0cm>*{\bullet};<-5mm,18mm>*{i_1}**@{},
 <5mm,0cm>*{\bullet};<5mm,18mm>*{i_2}**@{},
 <5mm,0cm>*{\bullet};<11mm,10mm>*{i_3}**@{},
 \end{xy}\ \ - \ \
 \begin{xy}
 <5mm,0mm>*{\bullet};<0cm,7mm>*{\bullet}**@{~},
 <5mm,0cm>*{\bullet};<10mm,7mm>*{\bullet}**@{~},
 <0mm,7mm>*{\bullet};<-5mm,14mm>*{\bullet}**@{-},
 <0mm,7mm>*{\bullet};<5mm,14mm>*{\bullet}**@{-},
 <5mm,0cm>*{\bullet};<-5mm,18mm>*{i_1}**@{},
 <5mm,0cm>*{\bullet};<5mm,18mm>*{i_3}**@{},
 <5mm,0cm>*{\bullet};<11mm,10mm>*{i_2}**@{},
 \end{xy}\ \ - \ \
 \begin{xy}
 <5mm,0mm>*{\bullet};<0cm,7mm>*{\bullet}**@{~},
 <5mm,0cm>*{\bullet};<10mm,7mm>*{\bullet}**@{~},
 <10mm,7mm>*{\bullet};<5mm,14mm>*{\bullet}**@{-},
 <10mm,7mm>*{\bullet};<15mm,14mm>*{\bullet}**@{-},
 <5mm,0cm>*{\bullet};<-1mm,11mm>*{i_1}**@{},
 <5mm,0cm>*{\bullet};<5mm,18mm>*{i_2}**@{},
 <5mm,0cm>*{\bullet};<16mm,18mm>*{i_3}**@{},
 \end{xy}
 $$

 Algebras over the  quotient operad,
 $$
 \cG erst := Free(\caL)/ <\cI_{\cA}, I_{\caL}, I_{\cA\caL}>,
 $$
 are called (non-commutative) {\em Gerstenhaber}\, algebras. These are
 triples, $(V, \circ, [\, \bullet \, ])$, consisting
 of a graded vector space $V$, a degree $0$ associative product,
 $\circ: V\ot V \rar V$ and a degree $-1$ Lie bracket, $[\ \bullet\ ]: \odot^2 V \rar V$
 which satisfy the following compatibility condition,
 $$
 [a\bullet (b\circ c)]=[a\bullet b]\circ c +(-1)^{(\tl{a}+1)\tl{b}}b\circ
 [a\bullet c],
 $$
 for all homogeneous $a,b,c\in V$.

 \bip

 {\bf 1.8. Example.} A Gerstenhaber algebra $(V, \circ, [\, \bullet\, ])$ is called
 {\em graded commutative}\, if such is the product $\circ$. Let us denote by $\cG$ the operad
 which governs graded commutative Gerstenhaber algebras.

 \bip

 {\bf  1.9. Remark.} There is a canonical map of operads,
 $$
 \cG \lon \cG erst,
 $$
 corresponding to the obvious functor
 $$
 \left\{\Ba{c}
 {\mathrm A\ category\ of}\\
 \cG\hspace{-1mm}
 -\hspace{-1mm}{\mathrm algebras}
 \Ea\right\} \lon
 \left\{\Ba{c}
 {\mathrm A\ category\ of}\\
 \cG erst\hspace{-1mm}-\hspace{-1mm}{\mathrm algebras} \Ea\right\}
 $$
 which simply forgets graded commutativity of the associated product.

 \sip

Both operads $\cG$ and $\cG erst$ are composed from a pair of simpler
operads, $(\cC omm, \caL ie)$ and, respectively, $(\cA ss, \caL ie)$.
The difference, however, is that the composition of $(\cC omm, \caL ie)$
into $\cG$ satisfies the  distributive law \cite{Mar1}, while
the composition of $(\cA ss, \caL ie)$
into $\cG erst$ {\em does not}. Indeed,  ``opening'' the expression
 $$
 [(a_1\circ a_2) \bullet (a_3\circ a_4)]\ \ \simeq \ \
 \begin{xy}
 <7mm,-7mm>*{\bullet};<0cm,1mm>*{\bullet}**@{-},
 <7mm,-7mm>*{\bullet};<14mm,1mm>*{\bullet}**@{-},
 <0mm,1mm>*{\bullet};<-5mm,8mm>*{\bullet}**@{~},
 <0mm,1mm>*{\bullet};<5mm,8mm>*{\bullet}**@{~},
 <14mm,1mm>*{\bullet};<19mm,8mm>*{\bullet}**@{~},
 <14mm,1mm>*{\bullet};<9mm,8mm>*{\bullet}**@{~},
 <7mm,-7mm>*{\bullet};<-5mm,11mm>*{1}**@{},
 <7mm,-7mm>*{\bullet};<5mm,11mm>*{2}**@{},
 <7mm,-7mm>*{\bullet};<9mm,11mm>*{3}**@{},
 <7mm,-7mm>*{\bullet};<19mm,11mm>*{4}**@{},
 \end{xy}
 $$
 in two possible ways,
 \Beqrn
 [(a_1\circ a_2) \bullet (a_3\circ a_4)] &=&  a_1\circ[a_2\bullet (a_3\circ a_4)] +
 (-1)^{|a_2|(|a_3|+|a_4|+1)}
 [a_1\bullet (a_3\circ a_4)]\circ a_2 \\
 &=&
  a_1\circ[a_2\bullet a_3]\circ a_4 + (-1)^{|a_3|(|a_2|+1)} a_1\circ a_3\circ [a_2\bullet  a_4)] \\
 && + (-1)^{|a_2|(|a_3|+|a_4|+1)}\left( [a_1\bullet a_3]\circ a_4\circ a_2 +
 (-1)^{ |a_3|(|a_1|+1)}
 a_3 \circ [a_1\bullet  a_4]\circ a_2\right), \\
 \Eeqrn
 \Beqrn
 [(a_1\circ a_2) \bullet (a_3\circ a_4)] &=&  [(a_1\circ a_2)\bullet a_3]\circ a_4 +
 (-1)^{|a_3|(|a_1|+|a_2|+1)}
 a_3\circ [(a_1\circ a_2)\bullet  a_4] \\
 &=&
  a_1\circ [a_2\bullet a_3]\circ a_4 + (-1)^{|a_2|(|a_3|+1)}[a_1\bullet a_3]\circ a_2 \circ a_4\\
 &&+(-1)^{|a_3|(|a_1|+|a_2|+1)} \left(a_3\circ a_1 \circ [
 a_2\bullet  a_4] +(-1)^{|a_2|(|a_4|+1)}a_3\circ [a_1\bullet
 a_4]\circ a_2\right) 
\Eeqrn 
and then decomposing the associated
 relation in $\cG erst(4)$ into irreducibles, one gets
 $$
 
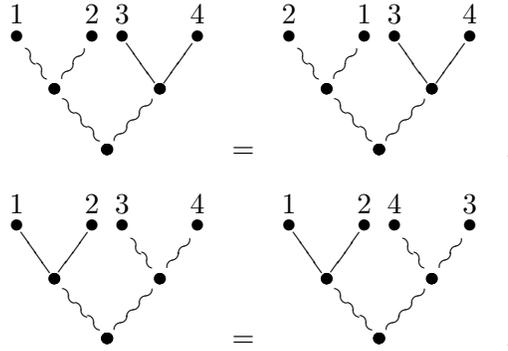
\begin{xy}
 <7mm,0mm>*{\bullet};<0cm,8mm>*{\bullet}**@{~},
 <7mm,0cm>*{\bullet};<14mm,8mm>*{\bullet}**@{~},
 <0mm,8mm>*{\bullet};<-5mm,15mm>*{\bullet}**@{~},
 <0mm,8mm>*{\bullet};<5mm,15mm>*{\bullet}**@{~},
 <14mm,8mm>*{\bullet};<19mm,15mm>*{\bullet}**@{-},
 <14mm,8mm>*{\bullet};<9mm,15mm>*{\bullet}**@{-},
 <7mm,0cm>*{\bullet};<-5mm,18mm>*{1}**@{},
 <7mm,0cm>*{\bullet};<5mm,18mm>*{2}**@{},
 <7mm,0cm>*{\bullet};<9mm,18mm>*{3}**@{},
 <7mm,0cm>*{\bullet};<19mm,18mm>*{4}**@{},
 \end{xy}
 \ \ = \ \
 \begin{xy}
 <7mm,0mm>*{\bullet};<0cm,8mm>*{\bullet}**@{~},
 <7mm,0cm>*{\bullet};<14mm,8mm>*{\bullet}**@{~},
 <0mm,8mm>*{\bullet};<-5mm,15mm>*{\bullet}**@{~},
 <0mm,8mm>*{\bullet};<5mm,15mm>*{\bullet}**@{~},
 <14mm,8mm>*{\bullet};<19mm,15mm>*{\bullet}**@{-},
 <14mm,8mm>*{\bullet};<9mm,15mm>*{\bullet}**@{-},
 <7mm,0cm>*{\bullet};<-5mm,18mm>*{2}**@{},
 <7mm,0cm>*{\bullet};<5mm,18mm>*{1}**@{},
 <7mm,0cm>*{\bullet};<9mm,18mm>*{3}**@{},
 <7mm,0cm>*{\bullet};<19mm,18mm>*{4}**@{},
 \end{xy}\ \ \ ,
 $$
 $$
 \begin{xy}
 <7mm,0mm>*{\bullet};<0cm,8mm>*{\bullet}**@{~},
 <7mm,0cm>*{\bullet};<14mm,8mm>*{\bullet}**@{~},
 <0mm,8mm>*{\bullet};<-5mm,15mm>*{\bullet}**@{-},
 <0mm,8mm>*{\bullet};<5mm,15mm>*{\bullet}**@{-},
 <14mm,8mm>*{\bullet};<19mm,15mm>*{\bullet}**@{~},
 <14mm,8mm>*{\bullet};<9mm,15mm>*{\bullet}**@{~},
 <7mm,0cm>*{\bullet};<-5mm,18mm>*{1}**@{},
 <7mm,0cm>*{\bullet};<5mm,18mm>*{2}**@{},
 <7mm,0cm>*{\bullet};<9mm,18mm>*{3}**@{},
 <7mm,0cm>*{\bullet};<19mm,18mm>*{4}**@{},
 \end{xy}
 \ \ = \ \
 \begin{xy}
 <7mm,0mm>*{\bullet};<0cm,8mm>*{\bullet}**@{~},
 <7mm,0cm>*{\bullet};<14mm,8mm>*{\bullet}**@{~},
 <0mm,8mm>*{\bullet};<-5mm,15mm>*{\bullet}**@{-},
 <0mm,8mm>*{\bullet};<5mm,15mm>*{\bullet}**@{-},
 <14mm,8mm>*{\bullet};<19mm,15mm>*{\bullet}**@{~},
 <14mm,8mm>*{\bullet};<9mm,15mm>*{\bullet}**@{~},
 <7mm,0cm>*{\bullet};<-5mm,18mm>*{1}**@{},
 <7mm,0cm>*{\bullet};<5mm,18mm>*{2}**@{},
 <7mm,0cm>*{\bullet};<9mm,18mm>*{4}**@{},
 <7mm,0cm>*{\bullet};<19mm,18mm>*{3}**@{},
 \end{xy}\ \ \ .
 $$
 The resulting relations in $\cG erst(4)$ are non-trivial unless
 the product $\circ$ is graded commutative.


 \bip

 {\bf 1.10. Example}. Let $\cA\cC$ be an $\bS$-module given by
 $\cA\cC(n): = \cA(n)[-1] \oplus \cC omm(n)$.
 Its only non-vanishing component $\cA\cC(2)$ is a 3-dimensional vector space spanned by
 two planar corollas of degree $1$ and one space corolla of degree 0,
 $$
 
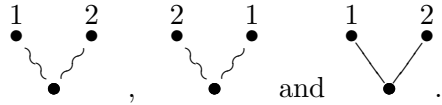
\begin{xy}
 <5mm,0cm>*{\bullet};<0cm,7mm>*{\bullet}**@{~},
 <5mm,0cm>*{\bullet};<10mm,7mm>*{\bullet}**@{~},
 <5mm,0cm>*{\bullet};<0mm,10mm>*{1}**@{},
 <5mm,0cm>*{\bullet};<10mm,10mm>*{2}**@{},
 \end{xy}
 \ \ \ , \ \ \
 \begin{xy}
 <5mm,0cm>*{\bullet};<0cm,7mm>*{\bullet}**@{~},
 <5mm,0cm>*{\bullet};<10mm,7mm>*{\bullet}**@{~},
 <5mm,0cm>*{\bullet};<0mm,10mm>*{2}**@{},
 <5mm,0cm>*{\bullet};<10mm,10mm>*{1}**@{},
 \end{xy} \ \  {\mathrm and} \ \
 \begin{xy}
 <5mm,0cm>*{\bullet};<0cm,7mm>*{\bullet}**@{-},
 <5mm,0cm>*{\bullet};<10mm,7mm>*{\bullet}**@{-},
 <5mm,0cm>*{\bullet};<0mm,10mm>*{1}**@{},
 <5mm,0cm>*{\bullet};<10mm,10mm>*{2}**@{},
 \end{xy}.
 $$
 As in example 1.7,
 the associated free operad $Free(\cA\caL)$
 can be represented as a linear span
 of all possible isomorphism classes of  binary  $[n]$-trees in the 3-space ${\mathbb R}^3$
 with the condition that  all planar corollas are perpendicular to a fixed line in ${\mathbb R}^3$.
 The composition in $Free(\cA\cC)$ is given  by gluing  the root vertex of one such partially
 planar/ partially space
 tree to a tail vertex of another one.

 \sip

 Let $\cI_{\cA\cC}$ be the ideal in $Free(\cA\cC)$ generated by the following
 3+3=3! vectors,
 $$
 
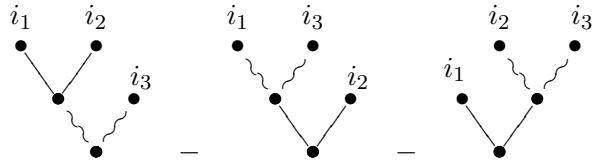
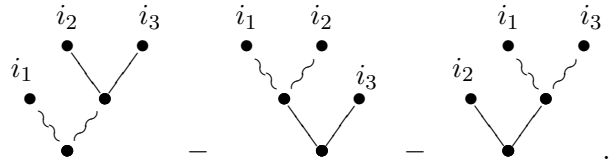
\begin{xy}
 <5mm,0mm>*{\bullet};<0cm,7mm>*{\bullet}**@{~},
 <5mm,0cm>*{\bullet};<10mm,7mm>*{\bullet}**@{~},
 <0mm,7mm>*{\bullet};<-5mm,14mm>*{\bullet}**@{-},
 <0mm,7mm>*{\bullet};<5mm,14mm>*{\bullet}**@{-},
 <5mm,0cm>*{\bullet};<-5mm,18mm>*{i_1}**@{},
 <5mm,0cm>*{\bullet};<5mm,18mm>*{i_2}**@{},
 <5mm,0cm>*{\bullet};<11mm,10mm>*{i_3}**@{},
 \end{xy}\ \ - \ \
 \begin{xy}
 <5mm,0mm>*{\bullet};<0cm,7mm>*{\bullet}**@{-},
 <5mm,0cm>*{\bullet};<10mm,7mm>*{\bullet}**@{-},
 <0mm,7mm>*{\bullet};<-5mm,14mm>*{\bullet}**@{~},
 <0mm,7mm>*{\bullet};<5mm,14mm>*{\bullet}**@{~},
 <5mm,0cm>*{\bullet};<-5mm,18mm>*{i_1}**@{},
 <5mm,0cm>*{\bullet};<5mm,18mm>*{i_3}**@{},
 <5mm,0cm>*{\bullet};<11mm,10mm>*{i_2}**@{},
 \end{xy}\ \ - \ \
 \begin{xy}
 <5mm,0mm>*{\bullet};<0cm,7mm>*{\bullet}**@{-},
 <5mm,0cm>*{\bullet};<10mm,7mm>*{\bullet}**@{-},
 <10mm,7mm>*{\bullet};<5mm,14mm>*{\bullet}**@{~},
 <10mm,7mm>*{\bullet};<15mm,14mm>*{\bullet}**@{~},
 <5mm,0cm>*{\bullet};<-1mm,11mm>*{i_1}**@{},
 <5mm,0cm>*{\bullet};<5mm,18mm>*{i_2}**@{},
 <5mm,0cm>*{\bullet};<16mm,18mm>*{i_3}**@{},
 \end{xy}
 $$
 and
 $$
 \begin{xy}
 <5mm,0mm>*{\bullet};<0cm,7mm>*{\bullet}**@{~},
 <5mm,0cm>*{\bullet};<10mm,7mm>*{\bullet}**@{~},
 <10mm,7mm>*{\bullet};<5mm,14mm>*{\bullet}**@{-},
 <10mm,7mm>*{\bullet};<15mm,14mm>*{\bullet}**@{-},
 <5mm,0cm>*{\bullet};<-1mm,11mm>*{i_1}**@{},
 <5mm,0cm>*{\bullet};<5mm,18mm>*{i_2}**@{},
 <5mm,0cm>*{\bullet};<16mm,18mm>*{i_3}**@{},
 \end{xy}
 \ \ - \ \
 \begin{xy}
 <5mm,0mm>*{\bullet};<0cm,7mm>*{\bullet}**@{-},
 <5mm,0cm>*{\bullet};<10mm,7mm>*{\bullet}**@{-},
 <0mm,7mm>*{\bullet};<-5mm,14mm>*{\bullet}**@{~},
 <0mm,7mm>*{\bullet};<5mm,14mm>*{\bullet}**@{~},
 <5mm,0cm>*{\bullet};<-5mm,18mm>*{i_1}**@{},
 <5mm,0cm>*{\bullet};<5mm,18mm>*{i_2}**@{},
 <5mm,0cm>*{\bullet};<11mm,10mm>*{i_3}**@{},
 \end{xy}\ \ - \ \
 \begin{xy}
 <5mm,0mm>*{\bullet};<0cm,7mm>*{\bullet}**@{-},
 <5mm,0cm>*{\bullet};<10mm,7mm>*{\bullet}**@{-},
 <10mm,7mm>*{\bullet};<5mm,14mm>*{\bullet}**@{~},
 <10mm,7mm>*{\bullet};<15mm,14mm>*{\bullet}**@{~},
 <5mm,0cm>*{\bullet};<-1mm,11mm>*{i_2}**@{},
 <5mm,0cm>*{\bullet};<5mm,18mm>*{i_1}**@{},
 <5mm,0cm>*{\bullet};<16mm,18mm>*{i_3}**@{},
 \end{xy}.
 $$

 Algebras over the  quotient operad,
 $$
 \caH := Free(\cA\cC)/ <\cI_{\cA[-1]}, I_{\cC}, I_{\cA\cC}>,
 $$
 are
 triples, $(V, \circ, \bullet )$, consisting
 of a graded vector space $V$, a degree $0$ associative graded commutative product,
 $\circ: V\odot V \rar V$ and a degree $1$ associative product $ \bullet:  V\ot V \rar V$
 which satisfy the following compatibility conditions,
 $$
 a\bullet (b\circ c)= (a\bullet b)\circ c +(-1)^{(\tl{a}+1)\tl{b}}b\circ
 (a\bullet c),
 $$
 $$
 (a\circ b)\bullet c= a\circ ( b\bullet c) +(-1)^{(\tl{c}+1)\tl{b}}
 (a\bullet c)\circ b,
 $$
 for all homogeneous $a,b,c\in V$.

 \bip

 \bip

 \begin{center}
 {\bf \S 2.  Strongly homotopy algebras}
 \end{center}

 \bip

 {\bf 2.1. Dg operads.} A differential graded (shortly, dg) operad is an operad,
 $$
 \f=\left( \left\{\f(n)\right\}_{n\geq 1}\
 ,\  \{\circ^{n,n'}_{i}\}_{ n,n'\geq 1 \atop 1\leq i\leq n}\right),
 $$
  in the sense of \S 1 together
 with a degree 1 equivariant linear map $d: \f(n) \rar \f(n)$, $\forall n$, satisfying the conditions,
 \Beqrn
 d^2 &=& 0, \\
 d\left(f\circ_i^{n,n'} f'\right) &=&
 (df)\circ_i^{n,n'} f' + (-1)^{|f|} f\circ_i^{n,n'} df', \ \  \forall f\in \f(n), f'\in \f(n').
 \Eeqrn
 The associated cohomology $\bS$-module $H(\f):= \{H^{\bullet}(\f(n))\}_{n\geq 1}$
 has an induced operad structure.

 \sip
 A morphism, $f: (\f,d) \rar (\f', d')$, of dg operads is, by definition, a morphism
 of operads $f: \f\rar \f'$
 which commutes in the obvious sense with the differentials.
 A morphism, $f: (\f,d) \rar (\f', d')$ is called a {\em quasi-isomorphism}\, if the induced
 morphism of the cohomology operads, $[f]: H(\f) \rar H(\f')$, is an isomorphism.

 \sip

 If $(V,d)$ is a dg vector space, then $(\cE_V,d_{\mathrm ind})$ is naturally
 a dg operad  where $d_{\mathrm ind}: Hom(\ot^{\bullet} V, V) \rar Hom(\ot^{\bullet} V, V)$
 is the differential which is
 naturally induced by $d$ and which we denote from now on by the same symbol $d$.

 \sip

 An {\em algebra}\, over a dg operad $(\f,d)$ is a dg vector space $(V,d)$ together with
 a morphism of dg operads $(\f,d)\rar (\cE_V, d)$.

 \sip

 It was shown in \cite{Mar} that for any simply connected dg operad $(\f, d)$ there
  exists a unique
 (up to an isomorphism) triple $\f_{\infty}:=(Free(\cE), d, f)$ where
 \Bi
 \item[(i)] $Free(\cE)$ is the free operad generated by an $\bS$-module $\cE=\{\cE(n)\}_{n\geq 2}$;
 \item[(ii)] $d$ is the differential in $Free(\cE)$ which is decomposable in the sense
 that
 $df \in Free^{\geq 2}(\cE)(n)$ for any $f\in \cE(n)$, $n\geq 2$.
 \item[(iii)] $f: (Free(\cE), d) \rar (\f,d)$ is a quasi-isomorphism of dg operads.
 \Ei
 This operad $\f_{\infty}$ is called the {\em minimal resolution}\footnote{
 An operad  $\f_{\infty}:=(Free(\cE), d)$ satisfying relations (i) and (ii) is often called
 {\em minimal}.} \, of the operad
 $\f$. Such minimal resolutions play a very important role in the homotopy theory of operadic
 algebras which we discuss below after considering  a few examples.

 \bip

 {\bf 2.2. Remark.}
 It is clear that to define a particular dg operad, $(\f, d)$,  is the same thing as to define
 its algebras, i.e.\ the image of the map  $(\f, d)\rar \cE_V$
 for some ``variable'' graded vector space $V$. Moreover, for this purely descriptive  purpose it is
 enough to assume that $\dim V< \infty$. We {\em always}\, make such an assumption when applying
 this method to concrete examples; in particular, there is never a  problem with
 replacing $V$ by its  dual, $V^*=\Hom(V,k)$.

 \bip

 {\bf 2.3. Example: operad $\cA_{\infty}$.}
  Let $\hat{\cA}$ be an $\bS$-module given by,
 $$
 \cA(n):= \left\{ \Ba{ll}
 0 & {\mathrm if}\ n=1 \\
 k[\bS_n][n-2] & {\mathrm if}\ n\geq2.
 \Ea
 \right.
 $$
 If we identify the natural basis of  $k[\bS_n][n-2]$
 with {\em planar}\,  $[n]$-corollas,
 $$
 
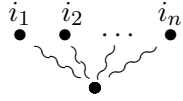
\begin{xy}
 <10mm,0cm>*{\bullet};<0cm,7mm>*{\bullet}**@{~},
 <10mm,0cm>*{\bullet};<6mm,7mm>*{\bullet}**@{~},
 <10mm,0cm>*{\bullet};<13mm,7mm>*{\cdots}**@{~},
 <10mm,0cm>*{\bullet};<20mm,7mm>*{\bullet}**@{~},
 <10mm,0cm>*{\bullet};<0mm,10mm>*{i_1}**@{},
 <10mm,0cm>*{\bullet};<7mm,10mm>*{i_2}**@{},
 <10mm,0cm>*{\bullet};<13mm,10mm>*{}**@{},
 <10mm,0cm>*{\bullet};<20mm,10mm>*{i_n}**@{},
 \end{xy}
 $$
 then the associated free operad $\cA_{\infty}:=(\{Free(\hat{\cA})(n), \circ_i^{n,n'}\})$
 can be represented as a linear span
 of all possible (isomorphism classes of) {\em planar}\,  $[n]$-trees
 with the compositions $\circ_i^{n,n'}$
 given simply by gluing  the root vertex of a planar $[n']$-tree to the $i$th tail vertex
 of an $[n]$-tree.

 \sip

 One can make $\cA_{\infty}$ into a dg operad with the differential $d$ given on generators by
 $$
 \begin{xy}
 <0mm,0cm>*{};<0cm,1mm>*{d}**@{},
 \end{xy}\ \
 \begin{xy}
 <10mm,0cm>*{\bullet};<0cm,7mm>*{\bullet}**@{~},
 <10mm,0cm>*{\bullet};<6mm,7mm>*{\bullet}**@{~},
 <10mm,0cm>*{\bullet};<13mm,7mm>*{\cdots}**@{~},
 <10mm,0cm>*{\bullet};<20mm,7mm>*{\bullet}**@{~},
 <10mm,0cm>*{\bullet};<0mm,10mm>*{1}**@{},
 <10mm,0cm>*{\bullet};<7mm,10mm>*{2}**@{},
 <10mm,0cm>*{\bullet};<13mm,10mm>*{}**@{},
 <10mm,0cm>*{\bullet};<20mm,10mm>*{n}**@{},
 \end{xy}
 \ = \
 \sum_{l+p=n+1\atop l,p\geq 2} \sum_{s=0}^{p-1}(-1)^{l+s(l+1)} \
 \begin{xy}
 <23mm,0cm>*{\bullet};<5mm,1cm>*{\bullet}**@{~},
 <23mm,0cm>*{\bullet};<10mm,1cm>*{\ldots}**@{~},
 <23mm,0cm>*{\bullet};<15mm,1cm>*{\bullet}**@{~},
 <23mm,0cm>*{\bullet};<21mm,1cm>*{\bullet}**@{~},
 <23mm,0cm>*{\bullet};<29mm,1cm>*{\bullet}**@{~},
 <23mm,0cm>*{\bullet};<34mm,1cm>*{\ldots}**@{~},
 <23mm,0cm>*{\bullet};<39mm,1cm>*{\bullet}**@{~},
 <23mm,0cm>*{\bullet};<5mm,13mm>*{1}**@{},
 <23mm,0cm>*{\bullet};<15mm,13mm>*{s}**@{},
 <23mm,0cm>*{\bullet};<29mm,13mm>*{s\hspace{-1mm}
 +\hspace{-1mm}l\hspace{-1mm}+\hspace{-1mm}1}**@{},
 <23mm,0cm>*{\bullet};<39mm,13mm>*{n}**@{},
 <21mm,1cm>*{{\bullet}};<15mm,21mm>*{\bullet}**@{~},
 <21mm,1cm>*{\bullet};<21mm,21mm>*{\cdots}**@{~},
 <21mm,1cm>*{\bullet};<27mm,21mm>*{\bullet}**@{~},
 <21mm,1cm>*{\bullet};<14mm,24mm>*{s\hspace{-1mm}+\hspace{-1mm}1}**@{},
 <21mm,1cm>*{\bullet};<27mm,24mm>*{s\hspace{-1mm}+\hspace{-1mm}l}**@{},
 \end{xy}
 $$
 The associated cohomology operad,
 $H(\cA_{\infty}, d)$,  is in fact isomorphic to $\cA ss$ \cite{GK}. Hence, the
 natural morphism of dg operads,
 $$
 f: (\cA_{\infty}, d) \lon (\cA ss, 0),
 $$
 defined to be identity on $[2]$-corollas and zero on $[n\geq 3]$-corollas, is a
 quasi-isomorphism. Thus $(\cA_{\infty}, d)$ is the minimal resolution of $\cA ss$.

 \sip

 An algebra over the dg operad $\cA_{\infty}$ is called a {\em strongly homotopy associative
 algebra}\,
 or, shortly, an $A_{\infty}$-{\em algebra}. This is a dg vector space $(V,d)$ equipped with
 degree $2-n$ multilinear operations $\mu_n: \ot^n V\rar V$, $n \geq 2$, such that
 for any $N\geq 1$ and any $v_1, \ldots, v_N\in V$,
 $$
 \sum_{l+p=N+1\atop l,p\geq 1}
 (-1)^r \mu_p(v_1,\ldots, v_s, \mu_l(v_{s+1}, \ldots, v_{s+l}), v_{s+l+1},\ldots, v_N)=0,
 $$
 where $\mu_1=d$ and $r=l+s(l+1) + p(|v_1| + \ldots + |v_s|$. If all $\mu_n$
 except $n=2$
 vanish, the above equation translates into the associative condition for the product
 $v_1\circ v_2:=\mu_2(v_1, v_2)$. Strongly homotopy associative algebras have been invented by
 Stasheff \cite{St} in his study of spaces homotopy equivalent to loop spaces.

 \bip

 {\bf 2.4. Example: operad $\cC_{\infty}$.}
  The dg operad $\cA_{\infty}$ has a commutative analog, $\cC_{\infty}$, which provides us
 with  the minimal resolution of the operad $\cC omm$ from Example~1.5. Using  remark~2.2
 one can  describe $\cC_{\infty}$ as follows:
  a $\cC_{\infty}$-{\em algebra}\, is,
 by definition, an $\cA_{\infty}$-algebra $(V, \{\mu_n\}_{n\geq 1})$ such that
 every multilinear operation $\mu_n:  \ot^n V \rar V$ is a Harrison cochain, that is,
 vanishes on every {\em shuffle product} which is given on generators by the formula
 $$
 (v_1\ot\ldots \ot v_i){\cyr x} (v_{i+1}\ot\ldots \ot v_n)=
 \sum_{\sigma\in Sh(i,n)} (-1)^{\sigma} v_{\sigma(1)}\ot \ldots
 \ot v_{\sigma(n)}.
 $$
 Here $(-1)^{\sigma}$ is the standard Koszul sign and $Sh(i,n)$ stands for a subset of $\bS_n$
 consisting of all permutations satisfying $\sigma^{-1}(1)<\ldots<\sigma^{-1}(i)$,
  $\sigma^{-1}(i+1)<\ldots<\sigma^{-1}(n)$. In particular, $\mu_2$ must be graded
 commutative, $\mu_2(v_1,v_2)= (-1)^{|v_1||v_2|}\mu_2(v_2,v_1)$.

 \bip

 {\bf 2.5. Example: operad $\caL_{\infty}$.} Let $\hat{\caL}$ be an $\bS$-module given by,
 $$
 \hat{\caL}(n):= \left\{ \Ba{ll}
 0 & {\mathrm if}\ n=1 \\
 {\mathbf 1}_n[2n-3] & {\mathrm if}\ n\geq2.
 \Ea
 \right.
 $$
 If we identify a basis vector of the one dimensional vector space   ${\mathbf 1}_n[2n-3]$
 with the (unique, up to an isomorphism) {\em space}\,  $[n]$-corolla,
 $$
 
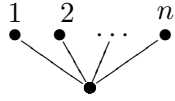
\begin{xy}
 <10mm,0cm>*{\bullet};<0cm,7mm>*{\bullet}**@{-},
 <10mm,0cm>*{\bullet};<6mm,7mm>*{\bullet}**@{-},
 <10mm,0cm>*{\bullet};<13mm,7mm>*{\cdots}**@{-},
 <10mm,0cm>*{\bullet};<20mm,7mm>*{\bullet}**@{-},
 <10mm,0cm>*{\bullet};<0mm,10mm>*{1}**@{},
 <10mm,0cm>*{\bullet};<7mm,10mm>*{2}**@{},
 <10mm,0cm>*{\bullet};<13mm,10mm>*{}**@{},
 <10mm,0cm>*{\bullet};<20mm,10mm>*{n}**@{},
 \end{xy}
 $$
 then the associated free operad $\caL_{\infty}:=(\{Free(\hat{\caL})(n), \circ_i^{n,n'}\})$
 can be represented as a linear span
 of all possible (isomorphism classes of) {\em space}\, $[n]$-trees
 with the compositions $\circ_i^{n,n'}$
 given  by gluing  the root vertex of a planar $[n']$-tree to the $i$th tail vertex
 of an $[n]$-tree. The point is that $\caL_{\infty}$ can be naturally made into a dg
 operad with the differential $d$ given on the generators by
 $$
 \begin{xy}
 <0mm,0cm>*{};<0cm,1mm>*{d}**@{},
 \end{xy}\ \
 \begin{xy}
 <10mm,0cm>*{\bullet};<0cm,7mm>*{\bullet}**@{-},
 <10mm,0cm>*{\bullet};<6mm,7mm>*{\bullet}**@{-},
 <10mm,0cm>*{\bullet};<13mm,7mm>*{\cdots}**@{-},
 <10mm,0cm>*{\bullet};<20mm,7mm>*{\bullet}**@{-},
 <10mm,0cm>*{\bullet};<0mm,10mm>*{1}**@{},
 <10mm,0cm>*{\bullet};<7mm,10mm>*{2}**@{},
 <10mm,0cm>*{\bullet};<13mm,10mm>*{}**@{},
 <10mm,0cm>*{\bullet};<20mm,10mm>*{n}**@{},
 \end{xy}
 \ = \
 \sum_{I_1\sqcup I_2=(1,\ldots,n)\atop \# I_1\geq 2, \# I_2\geq 2}\hspace{-3mm}
\begin{xy}
 <28mm,0cm>*{\bullet};<21mm,1cm>*{\bullet}**@{-},
 <28mm,0cm>*{\bullet};<26mm,1cm>*{\bullet}**@{-},
 <28mm,0cm>*{\bullet};<31mm,1cm>*{\cdots}**@{-},
 <28mm,0cm>*{\bullet};<36mm,1cm>*{\bullet}**@{-},
 <28mm,0cm>*{\bullet};<31mm,15mm>*{I_2}**@{},
 <28mm,0cm>*{\bullet};<31mm,12mm>*{\overbrace{\ \ \ \ \ \ \ \ \ \  }}**@{},
 <21mm,1cm>*{{\bullet}};<13mm,22mm>*{\bullet}**@{-},
 <21mm,1cm>*{\bullet};<18mm,22mm>*{\bullet}**@{-},
 <21mm,1cm>*{\bullet};<23mm,22mm>*{\cdots}**@{-},
 <21mm,1cm>*{\bullet};<28mm,22mm>*{\bullet}**@{-},
 <21mm,1cm>*{\bullet};<21mm,27mm>*{I_1}**@{},
 <21mm,1cm>*{\bullet};<21mm,24mm>*{\overbrace{\ \ \ \ \ \ \ \ \ \  \ \ \ \ \ }}**@{},
 \end{xy} .
 $$
 The associated cohomology operad,
 $H(\caL_{\infty}, d)$,  is isomorphic to $\caL ie$ \cite{GK}. Hence, the
 natural morphism of dg operads,
 $$
 f: (\caL_{\infty}, d) \lon (\caL ie, 0),
 $$
 defined to be identity on $[2]$-corollas and zero on $[n\geq 3]$-corollas, is a
 quasi-isomorphism. Thus $(\caL_{\infty}, d)$ is the minimal resolution of the operad of
 Lie algebras.

 \sip

 An algebra over the dg operad $\caL_{\infty}$ is called a {\em strongly homotopy Lie
 algebra}\,
 or, shortly, a $\caL_{\infty}$-{\em algebra}. This is a dg vector space $(V,d)$ equipped with
 degree $3-2n$ multilinear operations $\nu_n: \odot^n V\rar V$, $n \geq 2$, such that
 for any $N\geq 1$ and any $v_1, \ldots, v_N\in V$,
 $$
 \sum_{I_1\sqcup I_2=(1,\ldots,n)\atop  \# I_1\geq 1, \# I_2\geq 0}
 (-1)^{\sigma}  \mu_{\# I_2 +1}\left(
 \mu_{\# I_1}(v_{I_1}), v_{I_2}\right)=0,
 $$
 where  $\mu_1:=d$, $(-1)^{\sigma}$ is the standard Koszul sign associated with a
 permutation of the elements $v_1, \ldots , v_N$, and
 $$
 v_{I}:= v_{i_1}\ot\ldots\ot v_{i_l}, \ \ \ {\mathrm for}\ \
 I=(i_1,\ldots,i_l)\subset (1,\ldots, n).
 $$

 \bip

 \sip

 {\bf 2.6. Example: operad $\cG_{\infty}$.} Using remark~2.2 one can describe
 the minimal resolution, $\cG_{\infty}$, of the operad, $\cG$, of graded
 commutative Gerstenhaber algebras as follows \cite{GJ}\footnote{Actually, in
\cite{GJ} this operad was called {\em Fulton-MacPherson} operad. We, however, follow in this paper
the terminology and notations used in \cite{TT}}.

 \sip

 Let $V$ be a finite dimensional graded vector space and
 $$
 {\mathsf Lie}({V}^*[-2])= \bigoplus_{k=1}^{\infty}{\mathsf Lie}^k({V}^*[-2])
 $$
 the free graded Lie algebra generated by the shifted dual vector space, i.e.
 $$
 {\mathsf Lie}^1({V}^*[-2]):= {V}^*[-2],  \ \ \
 {\mathsf Lie}^k({V}^*[-2]):=
 \left[ {V}^*[-2]\bullet {\mathsf Lie}^{k-1}({V}^*[-2])
 \right].
 $$
 The Lie bracket on ${\mathsf Lie}({V}^*[-2])$ extends in the usual way
 to the completed (with respect to the natural filtration,
 ${\mathsf Lie}^{\geq k}({V}^*[-2])$)
 graded commutative associative algebra
 $$
 \hat{\odot}^{\bullet} {\mathsf Lie}({V}^*[-2]) = \prod_{k=0}^{\infty} \odot^k
 {\mathsf Lie}({V}^*[-2]),
 $$
 making the latter into a graded commutative Gerstenhaber algebra.

 \bip

 {\bf 2.6.1. Proposition} \cite{GJ} (see also \cite{Ta,TT}). {\em  A
 $G_{\infty}$-algebra structure on a finite-dimensional vector space
 $V$ is, by definition, a differential,
 of the  free $\cG$-algebra
 $\hat{\odot}^{\bullet} {\mathsf Lie}({V}^*[-2]) $.
 }

 \bip

 There is a one-to-one correspondence between
 arbitrary derivations, $D$,  of
 $\hat{\odot}^{\bullet} {\mathsf Lie}({V}^*[-2]) $
 and
  arbitrary collections of linear maps,
 $$
 m^*_{k_1,\ldots, k_n}: {V}^*[-2] \lon  {\mathsf Lie}^{k_1}({V}^*[-2])
 \odot \ldots\odot {\mathsf Lie}^{k_n}({V}^*[-2]).
 $$
 Upon dualization the latter go into linear homogeneous maps,
 $$
 m_{k_1,\ldots, k_n}: \frac{V^{\ot k_1}}{\mathrm shuffle\ products}
  \ot \ldots \ot \frac{V^{\ot k_n}}{\mathrm shuffle\ products}
 \lon V,
 $$
 of degree $3-n-k_1-\ldots-k_n$.
 The condition $D^2=0$ translates  into a well-defined
 set of quadratic equations for $m_{k_1,\ldots, k_n}$ which say,
 in particular, that $m_1$ is a differential on $V$ and that the product,
 $v_1\cdot v_2:=(-1)^{\tv_1}m_2(v_1,v_2)$, together with
 the Lie bracket,
  $[v_1\bullet v_2]:=-(-1)^{\tv_1}
 m_{1,1}(v_1,v_2)$, satisfy the Poisson identity up to a homotopy
 given by $m_{2,1}$. Hence the associated
 cohomology space $H(V,\mu_1)$ is a graded commutative Gerstenhaber algebra
 with respect to the  binary operations induced by  $m_2$ and $m_{1,1}$.

 \bip

 {\bf 2.7. A tower of approximations to the $\cG_{\infty}$ operad}.
 Let $V$ be a (finite-dimensional, see Remark 2.2) vector space
 and
 $\hat{\odot}^{\bullet} \left({\mathsf Lie}{V}^*[-2]\right)$
 the associated free $\cG$-algebra.
 It is easy to see that the
  multiplicative ideal, $I:= <{\mathsf Lie}^{\geq 2}{V}^*[-2]>$, generated by the commutant
  of ${\mathsf Lie}{V}^*[-2]$,
  as well as its multiplicative power $I^n$, $n\geq 2$, are also {\em Lie}\, ideals
  of $\hat{\odot}^{\bullet} \left({\mathsf Lie}V^*[-2]\right)$.
  Hence the quotient, $\hat{\odot}^{\bullet} \left({\mathsf Lie}{V}^*[-2]\right)/I^n$,
  $n\geq 1$, has a canonical $\cG$-algebra structure (note, however, that in the case
  $n=1$ the induced Lie brackets vanish).

 \bip

 {\bf 2.7.1. Definition of operad $\cG_{\infty}^{(n)}$:}
 A $\cG_{\infty}^{(n)}$-{\em algebra}\, structure on a finite-dimensional vector space
 $V$ is a differential
 of the  quotient $\cG$-algebra
 $\hat{\odot}^{\bullet} \left({\mathsf Lie}{V}^*[-2]\right)/I^n$, $n\geq 1$.

 \bip

 The following statements are obvious.

 \bip

 {\bf 2.7.2. Lemma.} $\cG_{\infty}^{(1)}= \caL_{\infty}$.

 \bip

 {\bf 2.7.3. Lemma.} {\em For any $n\geq 2$, there is a commutative diagram},
 $$
 \xymatrix{
 \cG_{\infty}^{(n)} \ar[r] & \cG_{\infty}  \\
  \caL_{\infty}\ar[u] \ar[r]^{=}  &  \caL_{\infty}\ar[u]
 }
 $$
 {\em with all arrows being canonical cofibrations\,\footnote{The notion of
 cofibration is explained in Sect.\ 2.8.} of operads.}

 \bip

 We shall give below in this paper a fairly explicit description of the next two floors,
 $\cG_{\infty}^{(2)}$ and $\cG_{\infty}^{(3)}$, of the tower,
 $$
 \caL_{\infty}=\cG_{\infty}^{(1)} \lon  \cG_{\infty}^{(2)} \lon \cG_{\infty}^{(3)}
  \lon \ldots \lon \cG_{\infty}^{(n)}\lon  \ldots,
 $$
  of cofibrant approximations to the operad $\cG_{\infty}={\mathrm colim} \cG_{\infty}^{(n)}$.
 Interestingly enough, the operad $\cG_{\infty}^{(2)}$ is closely related
 (through the cobar construction)
 to the operad $\cG erst$ and governs Frobenius$_{\infty}$ manifolds introduced in \cite{Me2},
 while the operad
 $\cG_{\infty}^{(3)}$ governs strong homotopy generalizations
  of Hertling-Manin's $F$-manifolds.

 \bip

 \sip

 {\bf 2.8. Homotopy theory.} The categories of operads and of their algebras belong to a
 class of so called closed model categories  \cite{Qu}  which
 have a particularly  nice homotopy theory. Here is a brief outline of all the relevant
 notions and facts we use in the paper (see, e.g., \cite{DS,GeMa} for more details and proofs).

 \bip

 {\bf 2.8.1. Definition}. A {\em closed model category}
  is, by definition, a category ${\cC at}$ with three distinguished classes of morphisms ---
 (i) {\em weak equivalences}, ${\frak E}$, (ii) {\em fibrations}, ${\frak F}$,
 and (iii) {\em cofibrations},
  ${\frak F}^{\circ}$,
 --- which are closed under composition and contain all identity maps. The following axioms
 must be satisfied:

 \noindent{\bf CMC1:} Finite limits and colimits exist in $\cC at$.

 \noindent{\bf CMC2:} If $f$ and $g$ are morphisms in $\cC at$ such that their composition
 $fg$ is defined, then if any two of the three maps $f$,$g$, $fg$ are weak equivalences
  then so is the third morphism.

 \noindent{\bf CMC3:} Given any commutative diagram of the form,
 $$
 \xymatrix{
 A \ar[r]^i \ar[d]_f & B \ar[d]^g \ar[r]^p & A \ar[d]_f \\
 A'\ar[r]^{i'} & B'\ar[r]^{p'} & A'
 }
 $$
 with $pi$ and $p'i'$ being the identity maps. If $g$ is in ${\frak E}$, ${\frak F}$ or
 ${\frak F}^{\circ}$,
 then so is $f$.

 \noindent{\bf CMC4:} Given any commutative solid arrow diagram of the form,
 $$
 \xymatrix{
 A \ar[r] \ar[d]_f & B \ar[d]^g \\
 A'\ar[r]\ar@{..>}[ur]^h & B'
 }
 $$
 A dotted arrow $h$ commuting with all other maps exists in either of the following two situations:
 (i) $f\in  {\frak F}^{\circ}\cap {\frak E}$ and $g\in {\frak F}$, or
 (ii) $f\in {\frak F}^{\circ}$ and $g\in  {\frak F}\cap {\frak E}$.

 \noindent{\bf CMC5} Any morphism can be factored in two ways: (i) $F=pi$ with
 $i\in {\frak F}^{\circ}$ and $p\in {\frak F}\cap {\frak E}$, and (ii) $f=pi$
 $i\in  {\frak F}^{\circ}\cap {\frak E}$ and $p\in {\frak F}$.

 \bip

 {\bf 2.8.2. Definitions}. (i) If a pair of morphisms, $f:A\rar A'$ and $g: B\rar B'$, satisfies
 the condition {\bf CMC4}, then we say that $f$ has the {\em left lifting property}\, (LLP)
 with respect to $g$ or that $g$ has the {\em right lifting property}\, (RLP)
 with respect to $f$.

 (ii) The morphisms in ${\frak F}\cap {\frak E}$ are called {\em acyclic fibrations}.
 The morphisms in ${\frak F}^{\circ}\cap {\frak E}$ are called {\em acyclic cofibrations}.

 (iii) The axiom {\bf CMC1} implies, in particular, that every closed
 model category $\cC at$ has both an initial object $\emptyset$ and a terminal object $*$.
 An object $A$ of $\cC at$ is called {\em fibrant} (resp., {\em cofibrant}) if
 $A\rar *$ is a fibration (resp., $\emptyset\rar A$ is a cofibration).

 \bip

 {\bf 2.8.3. Facts} \cite{GJ, Hi}. (i) The category of dg operads is a closed
 model category with
 \Bi
 \item weak equivalences ${\frak E}$ =\{ the morphisms, $f: (\f, d)\rar (\f',d')$,
  of dg operads which induce isomorphism, $[f]: H(\f)\rar H(\f')$, in cohomology\};
 \item fibrations ${\frak F}$=\{surjective morphisms of dg operads\};
 \item cofibrations ${\frak F}^{\circ}$=\{the morphisms which have LLP with respect to all acyclic
 fibrations \}.
 \Ei

 (ii) Given a dg operad $\f$, the associated category of $\f$-algebras is a closed model category
 with the classes  of morphism ${\frak E}$, ${\frak F}$ and ${\frak F}^{\circ}$ defined in
 close analogy to (i).

 (iii) Every object in the closed model categories (i) and (ii) is obviously fibrant.

 \bip

 {\bf 2.8.4. Homotopy and derived categories of a closed model category.}
 From now on $\cC at $  stands for  a closed model category. Moreover, we assume for simplicity
 that every object in $\cC at$ is fibrant (as in the two examples above).

 \sip

 Two morphisms, $f,g: A\rar B$, are called (right) {\em homotopic}\, if there exists a path object,
 for $B$ (that is an object  $B^I$ together with a diagram
 $$
 B\stackrel{i}{\lon} B^I\stackrel{p}{\lon} B\times B, \ \ \ i\in {\frak E},
 $$
 which factors the diagonal map $B\stackrel{(id, id)}{\lon} B\times B$) such that
 the product map $(f,g): A\rar B\times B$ lifts to a map $H: A\lon B^I$. Such a map
 $H$ is called a (right) homotopy from $f$ to $g$. This, in fact, defines an
 equivalence relation
 $\sim$
 in $\Hom_{\cC at}(A,B)$ for any objects $A,B$. However, the associated homotopy
 classes of maps,
 $$
 \pi(A,B):= \frac{\Hom_{\cC at}(A,B)}{\sim},
 $$
 do not necessarily compose, $\pi(A,B)\times \pi(B,C)\rar \pi(A,C)$, unless the
 objects involved are cofibrant.

 \sip

 By {\bf CMC5}(i), the map $\emptyset\rar A$ can be factored,
 $\emptyset\rar QA \stackrel{p_A}\rar A$,
 with $QA$ being cofibrant and $p_A$ a weak equivalence. Such an object $QA$ is called a {\em
 cofibrant
 resolution}\,  of $A$. Usually, cofibrant resolutions are constructed by the method
 of ``adding a new variable and killing a cycle''.
 A nice illustration of the method at work is, for example, Markl's  \cite{Mar} original
 construction of
 the minimal resolution $\f_{\infty}$ of a simply connected dg operad $\f$; Markl's
 minimal resolutions
 give an important class of cofibrant objects in the category of dg operads.

 \sip

 Another remarkable fact is that not only every object, but also every morphism,
 $f:A\rar B$, has a ``cofibrant resolution'' $Qf$ making the following diagram commutative,
 $$
 \xymatrix{
 QA \ar[r]^{Qf} \ar[d]_{p_A} & QB \ar[d]^{p_B} \\
 A\ar[r]^f  & B
 }
 $$
 Moreover, the homotopy class of such maps $[Qf]$ is defined uniquely by the homotopy class
 $[f]$.

 \sip

 The {\em homotopy category}, ${\mathbf Ho}(\cC at)$, is the category with the same objects as
 $\cC at$ and with morphisms given by
 $$
 \Hom_{{\mathbf Ho}(\cC at)}(A,B) := \pi(QA,QB).
 $$
 Clearly, there is a canonical functor $\al: {\cC at} \rar {\mathbf Ho}(\cC at)$,
 which is the identity on objects and sends morphisms $f$ to $[Qf]$.

 \sip

 The derived category, ${\mathbf D}(\cC at)$, is the category obtained from $\cC at$ by
 localization with respect to weak equivalences; put another way, this is a category
  together with a functor $F: \cC at \rar {\mathbf D}(\cC at)$ satisfying two conditions,
 \Bi
 \item $F(f)$ is an isomorphism for each weak equivalence $f$;
 \item every functor $G: \cC at  \rar \cC at'$ sending weak equivalences into isomorphisms
 factors uniquely through $({\mathbf D}(\cC at),F)$,
 $$
 G: \cC at \stackrel{F}{\lon} {\mathbf D}(\cC at) \stackrel{G'}{\lon} \cC at',
 $$
 for some functor $G'$.
 \Ei

 Note that the definition of ${\mathbf D}(\cC at)$ involves only one class, ${\frak E}$,
 of the three classes which define the closed model structure in $\cC at$. Nevertheless,
 one of the central results in Quillen's \cite{Qu} theory of closed model categories
 asserts the equivalence,
 $$
 {\mathbf Ho}(\cC at) \simeq {\mathbf D}(\cC at),
 $$
 of the two categories associated to $\cC at$.

 \bip

 {\bf 2.8.5. Transfer Theorem}.  {\em Let $\cP$ be a cofibrant dg operad
 and  $f:V\rar V' $ a quasi-isomorphism of dg vector spaces.
  For any $\cP$-algebra structure on $V$ (resp.\ $V'$),
 there exists a $\cP$-algebra structure on $V'$ (resp.\ $V$)
  so that $V$ and $V'$ are equivalent as $\cP$-algebras
 (i.e.\ there exists a $\cP$-algebra $V''$ and a pair of quasi-isomorphisms of $\cP$-algebras,
 $V\leftarrow V''\rar V'$). }

 \bip

 This is a well known fact. We show the proof only for
 completeness (cf.\ \cite{BM}).

 \sip

 {\bf Proof}.
 The dg operads $\cE_V$ and $\cE_{V'}$, if viewed only as
 dg $\bS$-modules,
  have two natural maps,
 $$
 \cE_V  \stackrel{f^*}{\lon} \cE_{V,V'},  \ \ \
 \cE_{V'}  \stackrel{f_*}{\lon} \cE_{V,V'},
 $$
  to the dg $\bS$-module  $\cE_{V,V'}:=\{Hom(V^{\ot n}, V'), d\}_{n\geq 1}$.
 Define a dg $\bS$-module, $\cE_f=\{\cE_f(n)\}_{n\geq 1}$, by setting
 $$
 \cE_f(n):= \left\{ (v,v')\in \cE_V(n)\times \cE_{V'}(n) \mid fv=v'f^{\ot n}\right\},
 $$
 or, alternatively, by the pullback diagram of dg $\bS$-modules,
 $$
 \xymatrix{
 \cE_f \ar[r]^{i_1} \ar[d]_{i_2} & \cE_{V} \ar[d]^{f_*} \\
 \cE_{V'}\ar[r]^{f^*}  & \cE_{V,V'}.
 }
 $$
 It is easy to check that
  $\cE_f$ inherits from the $\cE_V$ and $\cE_{V'}$ not only the $\bS$-module
 structure, but also the compositions $\circ_i^{n,n'}$. Thus $\cE_f$  is a dg operad
 with maps $i_1$ and $i_2$ above being morphisms of dg operads. Moreover, $i_1$ and $i_2$
 are weak
 equivalences (recall that we are working over a field of characteristic zero). Hence,
 as objects in ${\bf D}(\mathrm Oper)$, the dg operads $\cE_{V}$, $\cE_f$ and $\cE_{V'}$
 are isomorphic. Since the derived category of dg operads is equivalent to the homotopy
 category, we get isomorphisms of sets,
 $$
 \pi(\cP, \cE_V) \simeq \pi (\cP,\cE_f) \simeq \pi(\cP, \cE_{V'}).
 $$
 If, say, $V$ is a $\cP$-algebra, the homotopy equivalence class of the structure map
  $\phi: \cP \rar \cE_V$
  gives rise to an element $[\phi_f]$  in $\pi(\cP, \cE_{f})$. As $\cP$ is cofibrant,
 the latter has a representative, $\phi_f\in \Hom_{\mathrm Oper}(\cP, Q\cE_f)$, where
 $p_{f}: Q\cE_f \rar \cE_f$ is some cofibrant resolution of $\cE_f$.  By construction,
 the composition $i_1\circ p_f\circ \phi_f: \cP\rar \cE_V$ is homotopy equivalent
 to the original structure map  $\phi: \cP \rar \cE_V$. Finally, another composition,
 $i_1\circ\p_f\circ\phi_f: \cP\rar \cE_{V'}$ makes $V'$ into a $\cP$ algebra which,
 in the derived category of $\cP$-algebras, is obviously isomorphic to $\phi$.
  Analogously one proves the dual statement.


 \hfill$\Box$

 \sip


 {\bf 2.8.6.  Sh algebras.} An algebra over a cofibrant operad is called a
  {\em strongly homotopy}\,
 (or, shortly, sh) algebra. By the Theorem above, sh algebraic structures
  can be transferred by quasi-isomorphisms of complexes.

 \bip

 {\bf 2.9. Markl's theory of sh maps.} The beauty of sh algebras,
 the transfer property 2.8.5, is spoiled by the fact that to
 compare such structures on quasi-isomorphic dg spaces $V$ and $V'$
 (we refer to 2.8.5 again)
  one has to resort to a {\em chain}\, of strict $\cP$-algebra morphisms,
  $V\leftarrow V'' \rar\ V'$, involving a third party  which is often hard to construct
 explicitly. One may try to overcome this deficiency by appropriately extending the notion
 of {\em map}\, between sh algebras.

 \sip

 Markl made in \cite{Mar2,Mar3} an interesting suggestion which, in
 the setting of the proof of the Transfer Theorem, can be
 illustrated as follows. First one observes that
  the operad $\cE_f$ is in fact a two coloured operad with one colour  associated
  to
 $V$ and another one to $V'$. Next one constructs a  two coloured
 cofibrant resolution, $Q\cE_f$, and then defines the set of
 {\em sh maps}\, between the $\cP$-algebras $V$ and $V'$ as the set
 of all algebras over the dg operad $Q\cE_f$. In this way Markl was
 able to prove stronger versions of the Transfer Theorem
 \cite{Mar2}. The  problem, however, with this approach is that it
 is not yet clear whether or not such sh maps can be composed
 making the pair (sh algebras, sh maps) into a genuine category. At
 present, this is known to be true only for a class of sh algebras
 associated with  Koszul operads. In particular, it is true for
 $\cA_{\infty}$-, $\cC_{\infty}$- and $\caL_{\infty}$-algebras
  reproducing thereby  the well established theory of sh maps of these
 three classes of sh algebras.  For later reference we review below a few basic
 facts (see, e.g., \cite{Ko1,P}).

 \bip

 {\bf 2.10. Sh maps of $\cA_{\infty}$- and $\cC_{\infty}$-algebras.}
 An $A_{\infty}$-structure on a vector space $V$ can be suitably represented
 as a codifferential, $\mu: (T^{\bullet}V[1], \Delta) \rar
 (T^{\bullet}V[1],\Delta)$, of the free tensor coalgebra
 cogenerated by $V[1]$. A {\em sh map}, $f:(V,\mu_{\bullet})\rar (\hat{V},\hat{\mu}_{\bullet})$,
 of $\cA_{\infty}$-algebras is, by definition, a morphism of the associated
 differential coalgebras,
 $f: (T^{\bullet}V[1], \Delta, \mu) \rar (T^{\bullet}\hat{V}[1],\Delta, \hat{\mu})$.
 Such a map is equivalent to
 a set of linear maps
 $\{f_n:V^{\otimes n}\lon \tilde{V}, \ n\geq 1 \}$
  of degree  $1-n$
 which satisfy the equations,
 $$
 \sum_{1\le k_1<k_2<\ldots<k_i=n}(-1)^{i+r}
 \hat{\mu}_i(f_{k_1}(v_1,\ldots,v_{k_1}),f_{k_2-k_1}(v_{k_1+1},\ldots,v_{k_2}),
 \ldots,f_{n-k_{i-1}}(v_{k_{i-1}+1},\ldots,v_n))
 $$
 $$
 =\sum_{k+l=n+1}\sum_{j=0}^{k-1}(-1)^{l(\tv_1+\ldots+\tv_j+n)+j(l-1)}
 f_k(v_1,\ldots,v_j,\mu_l(v_{j+1},\ldots,v_{j+l}),v_{j+l+1},\ldots,v_n),
 $$
 for arbitrary $v_i\in V$.

 \sip

 The pair $(Ob=\cA_{\infty}{\mathrm\mbox{-}algebras}, Mor={\mathrm sh\ maps})$
 forms a category called
 the {\em category of $\cA_{\infty}$-algebras}.

 \sip

 A sh map $f=\{f_n\}:(V,\mu_{\bullet})\rar (\hat{V},\hat{\mu}_{\bullet})$ is called
 a {\em quasi-isomorphism}\, if the associated map of dg vector spaces,
 $f_1:(V,\mu_{1})\rar (\hat{V},\hat{\mu}_{1})$, induces an isomorphism in cohomology.

 \sip

 Two sh maps, $f,g: (T^{\bullet}V[1], \Delta, \mu) \rar
 T^{\bullet}\hat{V}[1],\Delta, \hat{\mu})$, are said to
 be {\em homotopic}\,  if there is a homogeneous map,
 $h: T^{\bullet}V[1] \rar T^{\bullet}\hat{V}[1]$, of degree $-1$ such that
 $$
 \Delta h= (f\ot h + h\ot g)\Delta, \ \ \ \ \ \ \
 f-g=\hat{\mu}\circ h + \mu\circ h.
 $$
 Remarkably enough, homotopy induces an equivalence relation in the set
 of sh maps $(V,\mu_{\bullet})\rar (\hat{V},\hat{\mu}_{\bullet})$ \cite{P}.
 Moreover, a sh map $f=\{f_n\}:(V,\mu_{\bullet})\rar (\hat{V},\hat{\mu}_{\bullet})$
 is a quasi-isomorphism if and only if it is a homotopy equivalence.
 Thus the derived category of $\cA_{\infty}$-algebras is simply the quotient of the
 category $\cA_{\infty}$-algebras by the above homotopy relation!

 \sip

 For $\cC_{\infty}$-algebras one has a similar list of definitions and results.

 \bip
 {\bf 2.11. Sh maps of $\caL_{\infty}$-algebras.}
 A $\caL_{\infty}$-algebra structure, $\nu=\{\nu_n: \odot^n V\rar V,\ |\nu_n|=3-2n\}$,
  on a vector space $V$ can be compactly described as a codifferential,
 $\nu: (\odot^{\bullet}V[2], \Delta) \rar
 (\odot^{\bullet}V[2],\Delta)$, of the free cocommutative tensor coalgebra
 cogenerated by $V[2]$. A {\em sh map}, $f:(V,\nu_{\bullet})\rar (\hat{V},\hat{\nu}_{\bullet})$,
 of $\caL_{\infty}$-algebras is, by definition, a morphism of the associated
 differential cocommutative  coalgebras,
 $f: (\odot^{\bullet}V[2], \Delta, \nu) \rar (\odot^{\bullet}\hat{V}[2],\Delta, \hat{\nu})$.
 Such a map is equivalent to
 a set of linear maps
 $\{f_n: \odot^n V\lon \tilde{V}, \ n\geq 1 \}$
  of degree  $2-2n$
 which satisfy the equations similar to the ones in Subsect.~2.10. The notions of
 quasi-isomorphism and homotopy are similar as well.

 \sip
 The pair $(Ob=\caL_{\infty}{\mathrm\mbox{-}algebras}, Mor={\mathrm sh\ maps})$
 forms a category called
 the {\em category of $\caL_{\infty}$-algebras}.
 \sip

 Dualizing the above formulae for a finite-dimensional vector space  $V$ one arrives at a
 beautiful
 geometric formulation of $\caL_{\infty}$-algebras and their sh maps \cite{Ko1}:
 \begin{itemize}
     \item a $\caL_{\infty}$-algebra structure on $V$ can be identified with
     a smooth degree 1 vector field $\vec{\nu}$ on the pointed flat graded
     manifold $(V[2],0)$ which satisfies the conditions
     $[\vec{\nu},\vec{\nu}]=0$ and  $\vec{\nu}|_0=0$. Explicitly, the identification,
     $$
     \vec{\nu}\ \longleftrightarrow\ \{\nu_n: \odot^n V\rar V \},
     $$
      is given by the formula,
      $$
      \vec{\nu}=\sum_{n=1}^{\infty}\ \sum_{\al,\be_1,\ldots\be_n}
      \frac{(-1)^r}{n!}t^{\be_1}
      \cdots t^{\be_n} \mu^{\al}_{\be_1,\ldots,\be_n} \frac{\p}{\p t^{\al}},
      $$
      where $\{t^{\al}, \al=1, \ldots, \dim V\}$ is the basis of $V^*[-2]$ associated
      to a basis,
      $\{e_{\al}\}$, of $V$ (so that $|t^{\al}|=-|e_{\al}|+2$),
      $$
      r=(2n-3)(|e_{\al_1}| +\ldots+|e_{\al_n}|) +\sum_{k=2}^n|e_{\al_k}|
      (|e_{\al_1}| +\ldots+|e_{\al_{k-1}}|),
      $$
      and $\mu^{\al}_{\be_1,\ldots,\be_n}\in k$ are given by
      $$
      \mu_{n}(e_{\be_1},\ldots,e_{\be_n})=\sum_{\al}\mu^{\al}_{\be_1,\ldots,\be_n} e_{\al}.
      $$
     \item A sh map of $f:(V,\nu_{\bullet})\rar (\hat{V},\hat{\nu}_{\bullet})$,
      of $\caL_{\infty}$-algebras is a smooth map of pointed graded manifolds,
      $f:(V^*[-2],0)\rar (\hat{V}^*[-2],0)$ such that $f_*(\vec{\nu})$ is well
      defined and coincides with $\vec{\hat{\nu}}$. Put another way,
      a sh map
      of $\caL_{\infty}$-algebras is just a morphism of the associated pointed
      dg manifolds.
 \end{itemize}

 \sip

 A $\caL_{\infty}$-algebra $(V,\{\nu_{\bullet}\}_{n\geq 1})$ with $\nu_1=0$ is called
 {\em minimal} (equivalently, the homological vector field $\vec{\nu}$ has zero at
 the distinguished point of order $\geq 2$).

 \bip

 {\bf 2.11.1. Facts \cite{Ko1}.} (i) {\em Every $\caL_{\infty}$-algebra is quasi-isomorphic
 to a minimal one}.

 (ii) {\em There is a one-to-one correspondence between
 quasi-isomorphisms of $\caL_{\infty}$-algebras and diffeomorphisms of the
 associated dg manifolds}.

 \bip

 {\bf 2.11.2. Fact \cite{Me2}.} {\em The canonical functor
 $$
 \left\{\begin{array}{c}
   {\mathrm The\ category} \\
   {\mathrm of}\ \caL_{\infty}{\mathrm \mbox{-}algebras} \\
 \end{array}
 \right\}
 \lon
 \left\{\begin{array}{c}
   {\mathrm The\ derived\ category} \\
   {\mathrm of}\ \caL_{\infty}{\mathrm \mbox{-}algebras} \\
 \end{array}
 \right\},
 $$
 when restricted to minimal $\caL_{\infty}$-algebras,
 becomes simply a forgetful functor,
 $$
 (M,\,*\, ,\, {\mathrm flat\ structure}\, ,\,\vec{\nu}) \lon
 (M,\, * \,  ,\, \vec{\nu}),
 $$
 which forgets the flat (=affine) structure on $(M,*)=(V^*[2],0)$.}

 \sip

 Thus a homotopy class of minimal $\caL_{\infty}$-algebras is nothing but
 a pointed formal dg manifold, $(M,\, * \,  ,\, \vec{\nu})$,
 with no preferred choice of local coordinates. Moreover, the derived (=homotopy)
 category of $L_{\infty}$-algebras is equivalent to the purely geometric
 category of formal dg manifolds.

 \bip

 \bip

 \pagebreak

 \begin{center}
 {\bf \S 3. Cobar construction for $\cG erst$}
 \end{center}

 \bip

 \sip


 {\bf 3.1. Cobar construction}. For an $\bS$-module $\f=\{\f(n)\}_{n\geq 1}$
 we set $\f\{m\}$ to be an $\bS$-module given by the tensor product,
 $$
 \f\{m\}(n):= \f(n)\ot_k \Lambda_n^{\ot m}[m(n-1)],
 $$
 where $\Lambda_n$ is the sign representation of the permutation group $\bS_n$. If $\f$
 is a dg operad, then $\f\{m\}$ is naturally a dg operad as well:
  a structure of $\f\{m\}$-algebra on a dg vector space $V$ is the same as a structure
 of $\f$-algebra on the shifted dg vector space $V[m]$.

 \sip

 Let $\f=\{\f(n),\circ_i^{n,n'}, d\}$ be a simply connected dg operad and let $\f^*[-1]$ stand for the
 $\bS$-module $\{\f(n)^*[-1]\}$. It was shown in \cite{GK} that the free operad associated
 to the $\bS$-module
 $Free(\f^*[-1]\{-1\})$ can be naturally made into a {\em differential}\,  operad,
 $\caD(\f)=(Free(\f^*[-1]\{-1\}), \delta)$,
 with the differential $\delta$ defined by that in $\f$ and the compositions $\circ_i^{n,n'}$.
 This construction gives rise to a functor, $\caD: Oper_1 \rar Oper_1$, on the category of
 simply connected dg operads
 with the property that $\caD(\caD(\f))$ is canonically quasi-isomorphic to the original
  operad $\f$.
 This functor is called a {\em cobar construction}.

 \bip

 {\bf 3.2. Koszul duals}.
 An operad $\f$ is called {\em quadratic}\, if it can be represented as a quotient,
 $$
 \f=\frac{Free(\cE)}{<R>},
 $$
 of the free operad generated by an $\bS$-module $\cE$ with $\cE(n)=0$ for $n\neq 2$
 by an ideal generated by an $\bS_3$-invariant subspace $R$ in
 $ Free(\cE)(3)$. For example, operads $\cC omm$, $\cA ss$ and $\caL ie$ are quadratic.

 \sip

 The {\em Koszul dual}\, of a quadratic operad $\f=Free(\cE)/<R>$ is, by definition,
 the quadratic operad $\f^!=Free(\hat{\cE})/<R^{\bot}>$ where $\hat{\cE}$ is the $\bS$-module
 whose only non-vanishing component is $\hat{\cE}(2)=\cE(2)^*\ot \Lambda_2$
 and $R^{\bot}$ is the annihilator of $R$, i.e.\
 the kernel of the natural map $Free(\hat{\cE})(3) \rar R^*$.

 \sip

 Applying cobar construction to the Koszul dual of a quadratic operad $\f$ one gets
 a cofibrant  dg operad $\caD(\f^!)$ together with  a canonical map of dg operads \cite{GK},
 $$
 (\caD(\f^!), \delta) \lon (\f,0).
 $$
 Whatever the operad $\f$ is, the associated $\caD(\f^!)$-algebras are strong
 homotopy ones.

 \bip

 {\bf 3.3. Koszul operads.} A quadratic operad $\f$ is called {\em Koszul} \, if the
 canonical map
 $(\caD(\f^!), \delta) \lon (\f,0)$ is a quasi-isomorphism. In such a case the cobar
 construction
 applied to the Koszul dual operad $\f^!$ provides us with the minimal resolution,
 $\f_{\infty}=\caD(\f^!)$,  of the operad $\f$ (see Sect.\ 2.1).

 \bip

 {\bf 3.4. Examples} \cite{GJ,GK}. The operads $\cA ss$, $\cC omm$, $\caL ie$ and $\cG$
 are  Koszul with
 $$
  \cA ss^! =\cA ss, \ \ \ \cC omm^!=\caL ie, \ \ \ \caL ie^!=\cC omm, \ \ \
 \cG^!= \cG\{1\}.
 $$

 \sip

 Thus the operads $\cA_{\infty}:=\caD(\cA ss)$, $\cC_{\infty}:=
 \caD(\caL ie)$, $\caL_{\infty}:=\caD(\cC omm)$ and
 $\cG_{\infty}:=\caD(G\{1\})$ are minimal resolutions of the
 operads, $\cA ss$, $\cC omm$, $\caL ie$ and $\cG$ respectively.
 This  explains all the claims made in Examples~2.3---2.6.

 \bip

 {\bf 3.5. Proposition.} {\em
 ${\cG erst}^!= \caH \{1\}$ and $\caH^!= \cG erst\{1\}$, where
 $\caH$ is the operad defined in Example~1.10.}

 \sip

 Proof is  straightforward.


 \bip

 {\bf 3.6. Cobar construction for $\cG erst^!$}. Surprisingly
 enough, the minimal operad\footnote{Conjecture: $Gerst_{\infty}$
 is a minimal resolution of $Gerst$, i.e.\ the operad $Gerst$ is
 Koszul.

 In this paper we need only the fact that the operad
 $Gerst_{\infty}$ is minimal so that $Gerst_{\infty}$-algebras are
 strong homotopy ones. } $Gerst_{\infty}:=\caD(\cG erst^!)$ turns
 out to be a much more elementary object than its ``graded
 commutative'' analogue
  $\cG_{\infty}=\caD(\cG^!)$.
 In this section we present an explicit and
  simple description of the cobar construction $Gerst_{\infty}$ in terms
 of partially planar/partially space trees (reflecting its nature
 as a composition of the operads $\cA_{\infty}$ and
 $\caL_{\infty}$, see below), and in the next section we show that
 strong homotopy $Gerst_{\infty}$-algebras admit a nice geometric
 interpretation.

 \sip

  The main reason behind that acclaimed simplicity of $\caD(\cG erst^!)=\caD(\caH\{1\})$
 is the non-distributive nature of the operad $\caH$ (cf.\ Remark~1.9): ``opening'' the left-
 and right hand sides of the following two equalities in the operad $\caH$,
 $$
 \begin{xy}
 <7mm,-7mm>*{\bullet};<0cm,1mm>*{\bullet}**@{~},
 <7mm,-7mm>*{\bullet};<14mm,1mm>*{\bullet}**@{~},
 <0mm,1mm>*{\bullet};<-5mm,8mm>*{\bullet}**@{~},
 <0mm,1mm>*{\bullet};<5mm,8mm>*{\bullet}**@{~},
 <14mm,1mm>*{\bullet};<19mm,8mm>*{\bullet}**@{-},
 <14mm,1mm>*{\bullet};<9mm,8mm>*{\bullet}**@{-},
 <7mm,-7mm>*{\bullet};<-5mm,11mm>*{1}**@{},
 <7mm,-7mm>*{\bullet};<5mm,11mm>*{2}**@{},
 <7mm,-7mm>*{\bullet};<9mm,11mm>*{3}**@{},
 <7mm,-7mm>*{\bullet};<19mm,11mm>*{4}**@{},
 \end{xy}
 \ \simeq \
 (a_1\bullet a_2) \bullet (a_3\circ a_4) =
 a_1\bullet\left( a_2 \bullet (a_3\circ a_4)\right)  \ \simeq \
 \begin{xy}
 <5mm,-7mm>*{\bullet};<0mm,0mm>*{\bullet}**@{~},
 <5mm,-7mm>*{\bullet};<10mm,0mm>*{\bullet}**@{~},
 <10mm,0mm>*{\bullet};<5mm,7mm>*{\bullet}**@{~},
 <10mm,0mm>*{\bullet};<15mm,7mm>*{\bullet}**@{~},
 <15mm,7mm>*{\bullet};<10mm,14mm>*{\bullet}**@{-},
 <15mm,7mm>*{\bullet};<20mm,14mm>*{\bullet}**@{-},
 <5mm,-7mm>*{\bullet};<-1mm,3mm>*{1}**@{},
 <5mm,-7mm>*{\bullet};<10mm,17mm>*{3}**@{},
 <5mm,-7mm>*{\bullet};<20mm,17mm>*{4}**@{},
 <5mm,-7mm>*{\bullet};<5mm,10mm>*{2}**@{},
 \end{xy}
 $$
 in two possible ways, one gets after the decomposition into
 irreducibles the following  relation in $\caH(4)$,
 $$
 \begin{xy}
 <7mm,-7mm>*{\bullet};<0cm,1mm>*{\bullet}**@{-},
 <7mm,-7mm>*{\bullet};<14mm,1mm>*{\bullet}**@{-},
 <0mm,1mm>*{\bullet};<-5mm,8mm>*{\bullet}**@{~},
 <0mm,1mm>*{\bullet};<5mm,8mm>*{\bullet}**@{~},
 <14mm,1mm>*{\bullet};<19mm,8mm>*{\bullet}**@{~},
 <14mm,1mm>*{\bullet};<9mm,8mm>*{\bullet}**@{~},
 <7mm,-7mm>*{\bullet};<-5mm,11mm>*{1}**@{},
 <7mm,-7mm>*{\bullet};<5mm,11mm>*{2}**@{},
 <7mm,-7mm>*{\bullet};<9mm,11mm>*{3}**@{},
 <7mm,-7mm>*{\bullet};<19mm,11mm>*{4}**@{},
 \end{xy}
 =0.
 $$
 This relation makes the computation of the $\bS$-module structure
 of the operad  $Gerst_{\infty}\simeq \caD(\caH\{1\})\simeq
 Free(\caH^*[-1]\{-2\})$ an easy exercise:
 $$
 \cG erst_{\infty}(n)={\mathbf 1}_{n}[2-3n]\ \ \oplus \ \
 \bigoplus_{p=2}^n {\mathrm Ind}^{\bS_n}_{\bS_p\times \bS_{n-p}}
 \left(k[\bS_p]\ot {\mathbf 1}_{n-p}\right)[n+2p-2].
 $$
 If we identify a basis vector of the one dimensional summand   ${\mathbf 1}_n[2n-3]$
 with the (unique, up to an isomorphism) {\em space}\,  $[n]$-corolla,
 $$
 
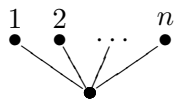
\begin{xy}
 <10mm,0cm>*{\bullet};<0cm,7mm>*{\bullet}**@{-},
 <10mm,0cm>*{\bullet};<6mm,7mm>*{\bullet}**@{-},
 <10mm,0cm>*{\bullet};<13mm,7mm>*{\cdots}**@{-},
 <10mm,0cm>*{\bullet};<20mm,7mm>*{\bullet}**@{-},
 <10mm,0cm>*{\bullet};<0mm,10mm>*{1}**@{},
 <10mm,0cm>*{\bullet};<6mm,10mm>*{2}**@{},
 <10mm,0cm>*{\bullet};<13mm,10mm>*{}**@{},
 <10mm,0cm>*{\bullet};<20mm,10mm>*{n}**@{},
 \end{xy}
 $$
 and the natural basis of all other summands with partially planar
 partially space $[n]$-corollas (cf.\ Sect.\ 1.7),

 $$
 
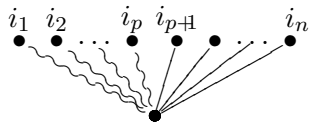
\begin{xy}
 <18mm,0cm>*{{\bullet}};<0cm,1cm>*{\bullet}**@{~},
 <18mm,0cm>*{\bullet};<5mm,1cm>*{\bullet}**@{~},
 <18mm,0cm>*{\bullet};<10mm,1cm>*{\ldots}**@{~},
 <18mm,0cm>*{\bullet};<15mm,1cm>*{\bullet}**@{~},
 <18mm,0cm>*{\bullet};<21mm,1cm>*{\bullet}**@{-},
 <18mm,0cm>*{\bullet};<26mm,1cm>*{\bullet}**@{-},
 <18mm,0cm>*{\bullet};<31mm,1cm>*{\ldots}**@{-},
 <18mm,0cm>*{\bullet};<36mm,1cm>*{\bullet}**@{-},
 <18mm,0cm>*{\bullet};<0mm,13mm>*{i_1}**@{},
 <18mm,0cm>*{\bullet};<5mm,13mm>*{i_2}**@{},
 <18mm,0cm>*{\bullet};<15mm,13mm>*{i_p}**@{},
 <18mm,0cm>*{\bullet};<21mm,13mm>*{i_{p+\hspace{-1mm}1}}**@{},
 <18mm,0cm>*{\bullet};<37mm,13mm>*{i_n}**@{},
 \end{xy} \ \ ,
 $$
 then $Gerst_{\infty}$, as an $\bS$-module.
 equals the linear span
 of all possible (isomorphism classes of) partially plane partially space
 trees formed by these corollas, while
  the compositions $\circ_i^{n,n'}$ are
 given  simply by gluing  the root vertex of an $[n']$-tree to the $i$th tail vertex
 of an $[n]$-tree. To complete the description of the operad $Gerst_{\infty}$
 we need only to compute Ginzburg-Kapranov's  cobar differential $d$.

 \sip

 {\bf 3.6.1. Proposition}. {\em The cobar differential in $\cG erst_{\infty}$
 is given on generators by}
 \Beqrn
 \begin{xy}
 <0mm,0cm>*{};<0cm,1mm>*{d}**@{},
 \end{xy}\ \
 \begin{xy}
 <10mm,0cm>*{\bullet};<0cm,7mm>*{\bullet}**@{-},
 <10mm,0cm>*{\bullet};<6mm,7mm>*{\bullet}**@{-},
 <10mm,0cm>*{\bullet};<13mm,7mm>*{\cdots}**@{-},
 <10mm,0cm>*{\bullet};<20mm,7mm>*{\bullet}**@{-},
 <10mm,0cm>*{\bullet};<0mm,10mm>*{1}**@{},
 <10mm,0cm>*{\bullet};<6mm,10mm>*{2}**@{},
 <10mm,0cm>*{\bullet};<13mm,10mm>*{}**@{},
 <10mm,0cm>*{\bullet};<20mm,10mm>*{n}**@{},
 \end{xy}
 & = &
 \sum_{I_1\sqcup I_2=(1,\ldots,n)\atop \# I_1\geq 2, \# I_2\geq 2}\hspace{-3mm}
 \begin{xy}
 <28mm,0cm>*{\bullet};<21mm,1cm>*{\bullet}**@{-},
 <28mm,0cm>*{\bullet};<26mm,1cm>*{\bullet}**@{-},
 <28mm,0cm>*{\bullet};<31mm,1cm>*{\cdots}**@{-},
 <28mm,0cm>*{\bullet};<36mm,1cm>*{\bullet}**@{-},
 <28mm,0cm>*{\bullet};<31mm,15mm>*{I_2}**@{},
 <28mm,0cm>*{\bullet};<31mm,12mm>*{\overbrace{\ \ \ \ \ \ \ \ \ \  }}**@{},
 <21mm,1cm>*{{\bullet}};<13mm,22mm>*{\bullet}**@{-},
 <21mm,1cm>*{\bullet};<18mm,22mm>*{\bullet}**@{-},
 <21mm,1cm>*{\bullet};<23mm,22mm>*{\cdots}**@{-},
 <21mm,1cm>*{\bullet};<28mm,22mm>*{\bullet}**@{-},
 <21mm,1cm>*{\bullet};<21mm,27mm>*{I_1}**@{},
 <21mm,1cm>*{\bullet};<21mm,24mm>*{\overbrace{\ \ \ \ \ \ \ \ \ \  \ \ \ \ \ }}**@{},
 \end{xy} ,
 \Eeqrn
\Beqrn
 \begin{xy}
 <18mm,0cm>*{{\bullet}};<0cm,1cm>*{\bullet}**@{~},
 <18mm,0cm>*{\bullet};<5mm,1cm>*{\bullet}**@{~},
 <18mm,0cm>*{\bullet};<10mm,1cm>*{\cdots}**@{~},
 <18mm,0cm>*{\bullet};<15mm,1cm>*{\bullet}**@{~},
 <18mm,0cm>*{\bullet};<21mm,1cm>*{\bullet}**@{-},
 <18mm,0cm>*{\bullet};<26mm,1cm>*{\bullet}**@{-},
 <18mm,0cm>*{\bullet};<31mm,1cm>*{\cdots}**@{-},
 <18mm,0cm>*{\bullet};<36mm,1cm>*{\bullet}**@{-},
 <18mm,0cm>*{\bullet};<0mm,13mm>*{1}**@{},
 <18mm,0cm>*{\bullet};<5mm,13mm>*{2}**@{},
 <18mm,0cm>*{\bullet};<15mm,13mm>*{p}**@{},
 <18mm,0cm>*{\bullet};<21mm,13.5mm>*{p\hspace{-0.5mm}+\hspace{-1mm}1}**@{},
 <18mm,0cm>*{\bullet};<37mm,13mm>*{n}**@{},
 <18mm,0cm>*{{\bullet}};<-4mm,5mm>*{d}**@{},
 \end{xy}
 &=& \vspace{30mm}
 \sum_{I_1\sqcup I_2=(p+1,\ldots,n) \atop \#I_1\geq 2, \#I_2\geq 2}\hspace{-5mm}
 \begin{xy}
 <18mm,0cm>*{{\bullet}};<0cm,1cm>*{\bullet}**@{~},
 <18mm,0cm>*{\bullet};<5mm,1cm>*{\bullet}**@{~},
 <18mm,0cm>*{\bullet};<10mm,1cm>*{\cdots}**@{~},
 <18mm,0cm>*{\bullet};<15mm,1cm>*{\bullet}**@{~},
 <18mm,0cm>*{\bullet};<21mm,1cm>*{\bullet}**@{-},
 <18mm,0cm>*{\bullet};<26mm,1cm>*{\bullet}**@{-},
 <18mm,0cm>*{\bullet};<31mm,1cm>*{\cdots}**@{-},
 <18mm,0cm>*{\bullet};<36mm,1cm>*{\bullet}**@{-},
 <18mm,0cm>*{\bullet};<0mm,13mm>*{1}**@{},
 <18mm,0cm>*{\bullet};<5mm,13mm>*{2}**@{},
 <18mm,0cm>*{\bullet};<15mm,13mm>*{p}**@{},
 <18mm,0cm>*{\bullet};<31mm,15mm>*{I_2}**@{},
 <18mm,0cm>*{\bullet};<31mm,12mm>*{\overbrace{\ \ \ \ \ \ \ \ \ \  }}**@{},
 <21mm,1cm>*{{\bullet}};<13mm,22mm>*{\bullet}**@{-},
 <21mm,1cm>*{\bullet};<18mm,22mm>*{\bullet}**@{-},
 <21mm,1cm>*{\bullet};<23mm,22mm>*{\cdots}**@{-},
 <21mm,1cm>*{\bullet};<28mm,22mm>*{\bullet}**@{-},
 <21mm,1cm>*{\bullet};<21mm,27mm>*{I_1}**@{},
 <21mm,1cm>*{\bullet};<21mm,24mm>*{\overbrace{\ \ \ \ \ \ \ \ \ \  \ \ \ \ \ }}**@{},
 \end{xy} \\
 \vspace{20mm}&&
 - \  \ \
 \sum_{I_1\sqcup I_2=(p+1,\ldots,n) \atop \#I_1\geq 0, \#I_2\geq 1}\hspace{-5mm}
 \begin{xy}
 <13mm,1cm>*{{\bullet}};<1mm,22mm>*{\bullet}**@{~},
 <13mm,1cm>*{\bullet};<6mm,22mm>*{\cdots}**@{~},
 <13mm,1cm>*{\bullet};<11mm,22mm>*{\bullet}**@{~},
 <13mm,1cm>*{\bullet};<16mm,22mm>*{\bullet}**@{-},
 <13mm,1cm>*{\bullet};<21mm,22mm>*{\cdots}**@{-},
 <13mm,1cm>*{\bullet};<26mm,22mm>*{\bullet}**@{-},
 <13mm,1cm>*{\bullet};<21mm,27mm>*{I_1}**@{},
 <13mm,1cm>*{\bullet};<21mm,24mm>*{\overbrace{\ \ \ \ \ \ \ \ \ \  }}**@{},
 <13mm,1cm>*{\bullet};<1mm,25mm>*{1}**@{},
 <13mm,1cm>*{\bullet};<11mm,25mm>*{p}**@{},
 <21mm,0cm>*{{\bullet}};<13mm,10mm>*{\bullet}**@{-},
 <21mm,0cm>*{\bullet};<20mm,10mm>*{\bullet}**@{-},
 <21mm,0cm>*{\bullet};<25mm,10mm>*{\cdots}**@{-},
 <21mm,0cm>*{\bullet};<30mm,10mm>*{\bullet}**@{-},
 <21mm,0cm>*{\bullet};<25mm,15mm>*{I_2}**@{},
 <21mm,0cm>*{\bullet};<25mm,12mm>*{\overbrace{\ \ \ \ \  \ \ \ \ \ }}**@{},
 \end{xy}
 \\
 &&
 +
 \sum_{I_1\sqcup I_2=(p+1,\ldots,n) \atop \#I_1\geq 1, \#I_2\geq 0}
  \sum_{s=1}^{p}
 \begin{xy}
 <28mm,0cm>*{{\bullet}};<-1mm,1cm>*{\bullet}**@{~},
 <28mm,0cm>*{\bullet};<4mm,1cm>*{\ldots}**@{~},
 <28mm,0cm>*{\bullet};<9mm,1cm>*{\bullet}**@{~},
 <28mm,0cm>*{\bullet};<17mm,1cm>*{\bullet}**@{~},
 <28mm,0cm>*{\bullet};<27mm,1cm>*{\bullet}**@{~},
 <28mm,0cm>*{\bullet};<32mm,1cm>*{\ldots}**@{~},
 <28mm,0cm>*{\bullet};<38mm,1cm>*{\bullet}**@{~},
 <28mm,0cm>*{\bullet};<44mm,1cm>*{\bullet}**@{-},
 <28mm,0cm>*{\bullet};<49mm,1cm>*{\ldots}**@{-},
 <28mm,0cm>*{\bullet};<53mm,1cm>*{\bullet}**@{-},
 <28mm,0cm>*{\bullet};<-1mm,13mm>*{\small 1}**@{},
 <28mm,0cm>*{\bullet};<8.5mm,13mm>*{s\hspace{-1mm}-\hspace{-1mm}1}**@{},
 <28mm,0cm>*{\bullet};<27mm,13mm>*{s\hspace{-1mm}
 +\hspace{-1mm}2}**@{},
 <28mm,0cm>*{\bullet};<37mm,13mm>*{p}**@{},
 <28mm,0cm>*{\bullet};<48mm,15mm>*{I_2}**@{},
 <28mm,0cm>*{\bullet};<49mm,12mm>*{\overbrace{\ \ \ \ \ \ \ \ \ \ \ }}**@{},
 <17mm,1cm>*{\bullet};<9mm,22mm>*{\bullet}**@{-},
 <17mm,1cm>*{\bullet};<15mm,22mm>*{\bullet}**@{-},
 <17mm,1cm>*{\bullet};<20mm,22mm>*{\ldots}**@{-},
 <17mm,1cm>*{\bullet};<25mm,22mm>*{\bullet}**@{-},
 <17mm,1cm>*{\bullet};<9mm,25mm>*{s}**@{},
 <28mm,0cm>*{\bullet};<20mm,27mm>*{I_1}**@{},
 <28mm,0cm>*{\bullet};<20mm,24mm>*{\overbrace{\ \ \ \ \ \ \ \ \ \ \ }}**@{},
 \end{xy}
 \\
 &&
 +
 \sum_{I_1\sqcup I_2=(p+1,\ldots,n) \atop \#I_1\geq 2, \#I_2\geq 2}
 \sum_{l+m=p+1\atop l,m\geq 2} \sum_{s=0}^{m-1}(-1)^{l+s(l+1)}
 \hspace{-12mm}
 \begin{xy}
 <28mm,0cm>*{{\bullet}};<0cm,1cm>*{\bullet}**@{~},
 <28mm,0cm>*{\bullet};<5mm,1cm>*{\cdots}**@{~},
 <28mm,0cm>*{\bullet};<10mm,1cm>*{\bullet}**@{~},
 <28mm,0cm>*{\bullet};<17mm,1cm>*{\bullet}**@{~},
 <28mm,0cm>*{\bullet};<27mm,1cm>*{\bullet}**@{~},
 <28mm,0cm>*{\bullet};<32mm,1cm>*{\ldots}**@{~},
 <28mm,0cm>*{\bullet};<38mm,1cm>*{\bullet}**@{~},
 <28mm,0cm>*{\bullet};<44mm,1cm>*{\bullet}**@{-},
 <28mm,0cm>*{\bullet};<49mm,1cm>*{\cdots}**@{-},
 <28mm,0cm>*{\bullet};<53mm,1cm>*{\bullet}**@{-},
 <28mm,0cm>*{\bullet};<0mm,13mm>*{\small 1}**@{},
 <28mm,0cm>*{\bullet};<10mm,13mm>*{s}**@{},
 <28mm,0cm>*{\bullet};<27mm,13mm>*{s\hspace{-1mm}+\hspace{-1mm}l\hspace{-1mm}
 +\hspace{-1mm}1}**@{},
 <28mm,0cm>*{\bullet};<37mm,13mm>*{p}**@{},
 <28mm,0cm>*{\bullet};<48mm,15mm>*{I_2}**@{},
 <28mm,0cm>*{\bullet};<49mm,12mm>*{\overbrace{\ \ \ \ \ \ \ \ \ }}**@{},
 <17mm,1cm>*{{\bullet}};<4mm,22mm>*{\bullet}**@{~},
 <17mm,1cm>*{\bullet};<9mm,22mm>*{\cdots}**@{~},
 <17mm,1cm>*{\bullet};<14mm,22mm>*{\bullet}**@{~},
 <17mm,1cm>*{\bullet};<20mm,22mm>*{\bullet}**@{-},
 <17mm,1cm>*{\bullet};<25mm,22mm>*{\cdots}**@{-},
 <17mm,1cm>*{\bullet};<30mm,22mm>*{\bullet}**@{-},
 <17mm,1cm>*{\bullet};<3mm,25mm>*{s\hspace{-1mm}+\hspace{-1mm}1}**@{},
 <17mm,1cm>*{\bullet};<14mm,25mm>*{s\hspace{-1mm}+\hspace{-1mm}l}**@{},
 <28mm,0cm>*{\bullet};<25mm,27mm>*{I_1}**@{},
 <28mm,0cm>*{\bullet};<25mm,24mm>*{\overbrace{\ \ \ \ \ \ \ \ \ \ \ }}**@{},
 \end{xy}
 \Eeqrn

 \bip

 Proof is straightforward though very tedious.

 \bip
 \bip

 {\bf Corollary 3.6.2.} {\em There is a canonical cofibration of
 operads}, $\caL_{\infty}\rar \cG erst_{\infty}$.

 \bip

 {\bf Corollary 3.6.3.} {\em A $\cG erst_{\infty}$-algebra is the data,
 $$
 \left(V, \ \{\nu_n\}_{n\geq 1}, \ \{\mu_{p,n}\}_{p\geq 2\atop n\geq 1 }
 \right),
 $$
 consisting of a graded vector space $V$ and two
 collections of homogeneous linear maps,
 \begin{itemize}
 \item $\nu_n: \odot^n V \rar V$ of degree $3-2n$,
 $n=1,2,3,\ldots$, and
 \item $\mu_{n;p}: (\otimes^n V)\otimes (\odot^p V)\rar V$ of
 degree $2-n-2p$, $k=2,3,4,\ldots$, $n=0,1,2,\ldots$,
 \end{itemize}
 which satisfy the equations, for any $a_1,\ldots,a_n, b_1,\dots,b_N\in V$,
 \bip
 $$
 \sum_{S_1\sqcup S_2=(1,\ldots N)}
 (-1)^{\sqcup}\nu_{|S_2|+1}(\nu_{|S_1|}(b_{S_1}), b_{S_2})=0, \ \ \ \  N\geq 1,
 $$
 and, for $ n\geq 2, N\geq 0$,
 $$
 \sum_{S_1\sqcup S_2=(1,\ldots N)}
 (-1)^{\sigma}\left\{
 \nu_{|S_2|+1}(\mu_{n;|S_1|}(a_1, \ldots, a_n; b_{S_1}), b_{S_2})\right. \hspace{65mm}
 $$
 $$
 - (-1)^{\tl{a}_1+\ldots \tl{a}_n-n}
 \mu_{n; |S_2|+1}(a_1, \ldots. a_n; \nu_{S_1}(b_{S_1}),b_{S_2})\hspace{75mm}
$$
 $$
 - \sum_{j=1}^n(-1)^{\tl{a}_1+\ldots \tl{a}_{j-1} + (\tl{a}_{j+1}+\ldots + \tl{a}_n)\tl{b}_{S_1}
 -n}\mu_{n; |S_2|}(a_1,\ldots,a_{j-1}, \nu_{|S_1|+1}(a_j,b_{S_1}),a_{j+1},\ldots, a_n; b_{S_2})
 { \left. \right\}}
 $$
 $$
 =\sum_{k+l=n+1\atop k,l\geq 2}\sum_{j=0}^k
 \sum_{S_1\sqcup S_2=(1,\ldots N)}
 (-1)^r\mu_{k;|S_2|}(a_1,\ldots, a_j, \mu_{l;|S_1|}(a_{j+1},\ldots,a_{j+l};b_{S_1}),
 a_{j+l+1},\ldots
 , a_n; b_{S_2}), \ \ \
 $$
 where $(-1)^{\sigma}$ is the standard Koszul sign of the shuffle permutation
 $b_1\ot\ldots\ot b_N\rar b_{S_1}\ot b_{S_2}$, and $r=j+l(n-j-l)+l(\tl{a}_1+\ldots+\tl{a}_j) +
 (\tl{a}_{j+l+1}+\ldots+\tl{a}_n)\tl{b}_{S_1}$.
 }

 \bip

 \sip

 Thus the operations $\nu_{\bullet}$ define on $V$ the structure
 of $L_{\infty}$ algebra, while the operations $\mu_{\bullet;0}$ define on $V$ the structure of
 $A_{\infty}$ algebra; the remaining operations $\mu_{\bullet;\bullet\geq 1}$ are homotopies,
 and homotopies of homotopies, and \ldots, which make these two basic structures Poisson-type
 consistent.
 In the special case when all operations but $\nu_2$ and $\mu_{2;0}$ vanish
 the $\cG erst_{\infty}$-algebra $V$ becomes nothing but a $\cG erst$-algebra,
 i.e.\ a (non-commutative, in general)
 Gerstenhaber algebra.

 \bip

 {\bf 3.7. Coalgebra interpretation}.
 The notion of $\caH$-algebra (which is a graded vector space $V$ equipped with two
 consistent associative multiplications, $\circ: V\odot V \rar V$ and $\bullet: V\ot V \rar V$,
 of degrees $0$ and $1$ respectively, see Sect.\ 1.10) can be naturally dualized
 leading to the notion of $\caH$-coalgebra.

 \sip

 Let $V$ be a graded vector space and $(\bar{B}V:=\ot^{\bullet\geq 1}V[1],\Delta_{\bullet})$
 the (reduced, as $\ot^0 V[1]$ is omitted) free tensor coalgebra generated by $V[1]$. Let
 $(\hat{\odot}^{\bullet}(\bar{B}V[1]), \Delta_{\circ})$ the completed (with respect to
 the natural filtration $F^{\geq r}:= \ot^{\geq r} V[1]$) free graded
 cocommutative coalgebra generated by $\bar{B}V[1]$ (the latter is $\bar{B}V$
 with shifted grading). Let
 $I\subset \hat{\odot}^{\bullet}(\ot^{\bullet}V[1]) $ be the coideal
 generated by $\ot^{\bullet\geq 2}V[1]$.
 The standard coassociative
 coproduct $\Delta_{\bullet}$
 in $\bar{B}V$ extends naturally to a degree -1 coproduct
  in the quotient $\hat{\odot}^{\bullet}(\bar{B}V[1])/I^2$
 making the latter into an $\caH$-coalgebra. (Here $I^2$ is the square of the coideal,
 i.e.\ the set of elements $x\in I$ such that $ \Delta_{\circ}x\in I\ot I$.)

 \bip

 {\bf 3.7.1. Proposition}. {\em
 There is  a one-to-one correspondence,}
 $$
 \left\{\Ba{c}
 Gerst_{\infty}\mbox{-}{\mathrm algebra}\\
 {\mathrm structures\  in }\ V
 \Ea
 \right\}
 \leftrightarrow
 \left\{\Ba{c}
 {\mathrm codifferentials\ in\ the\ cofree}\\
 {\caH}
 \mbox{-}
 {\mathrm  coalgebra }\
 \left(\hat{\odot}^{\bullet}(\bar{B}V[1])/I^2, \Delta_{\circ}, \Delta_{\bullet}\right)
 \Ea
 \right\}.
 $$

 \sip

 \Proof A coderivation $D$ of the coalgebra $\hat{\odot}^{\bullet}(\bar{B}V[1])/I^2$ is
 equivalent to a degree 1 linear map,
 $$
 \hat{\odot}^{\bullet}(\bar{B}V[1])/I^2 \lon V[2],
 $$
 vanishing on $k$, the direct summand of constants.
 As
 $$
 \hat{\odot}^{\bullet}(BV[1])/I^2 = \odot^{\bullet}(V[2]) \ \oplus \
 \left(\ot^{\bullet\geq 2}V[1]\right)[1]\ot \odot^{\bullet}(V[2]),
 $$
 this is the same as two collections of homogeneous linear maps,
 \{$\nu_n: \odot^n V \rar V$ of degree $3-2n$\}$_{n\geq 1}$
  and
 \{$\mu_{n;p}: (\otimes^n V)\otimes (\odot^p V)\rar V$ of
 degree $2-n-2p$\}$_{n\geq 2, p\geq 0}$.

 \sip

 It is straightforward to check that the condition $D^2=0$ translates precisely
 into the equations of Corollary 3.6.3.
 \hfill $\Box$

 \bip

 {\bf 3.8.  Geometric interpretation.} $\caL_{\infty}$-structure on a finite-dimensional
 vector space $V$ has a beautiful geometric interpretation
 \cite{Ko1} as
     a smooth degree 1 vector field $\vec{\nu}$ on the pointed affine graded formal
     manifold $(V[2],0)$ which satisfies the conditions
     $[\vec{\nu},\vec{\nu}]=0$ and  $\vec{\nu}|_0=0$.

 \sip

 Surprisingly enough, $\cG erst_{\infty}$-structure is also of purely geometric
 nature.

 \bip

 {\bf 3.8.1. Definitions.} (i) A {\em geometric $\cA_{\infty}$-structure}\, (respectively,
 {\em geometric $\cC_{\infty}$-structure}) on a graded manifold
 $\cM$ is the data $(\vec{\nu},\mu_{\bullet})$ consisting of
 \begin{itemize}
     \item[(a)]  a smooth homological vector field   $\vec{\nu}$ making $\cM$ into a dg manifold;
     \item[(b)]  a collection of maps,
     $$
     \mu_{\bullet}=\{\mu_n\}_{n\geq 1}: \ot^{\bullet}_{\f_{\cM}}\cT_{\cM}
     \rar \cT_{\cM},
      \ \ \ \
     {\mathrm with} \ \
     \mu_1=Lie_{\vec{\nu}},
     $$
 making the tangent sheaf $\cT_{\cM}$ into a sheaf of $\cA_{\infty}$-algebras
 (respectively, $\cC_{\infty}$-algebras).
 \end{itemize}

 \bip

 (ii) A geometric $\cA_{\infty}$-structure/$\cC_{\infty}$-structure on a
 pointed graded manifold $(\cM, *)$
 is called {\em minimal}\, if
 $\vec{\nu} I\subset I^2$ where $I$ is the ideal of the distinguished point $*\in \cM$.

 \bip

 \sip
 {\bf 3.8.2. Theorem.} {\em There is a one-to-one correspondence}
 $$
 \left\{\Ba{c}
 \cG erst_{\infty}\mbox{-}{\mathrm algebra}\\
 {\mathrm structures\  in }\ V
 \Ea
 \right\}
 \leftrightarrow
 \left\{\Ba{c}
 {\mathrm geometric}\ \cA_{\infty}\mbox{-}{\mathrm structures\ on \ the\ pointed}\\
 {\mathrm  affine\  formal\ manifold}\ \cM=(V[2],0)
 \Ea
 \right\}.
 $$

 \sip

 \Proof Let $\{e_{\al}, \al=1,\ldots, \dim V\}$ be basis of $V$ and
 $\{t^{\al}\}$ the associated dual basis of $V^*[-2]$
 which we identify with coordinate functions on $V[2]$. Set $t:= \sum_{\al} t^{\al}e_{\al}$,
 and let $\tau$ be the isomorphism of $\f_{\cM}$-modules
 $$
 \begin{array}{rccc}
   \tau: & \cT_{\cM} & \lon & \f_{\cM}\ot V \vspace{2mm}\\
    & X= \sum_{\al} X^{\al}(t)\p/\p t^{\al} & \lon  &
   \tau(X):= \sum_{\al} X^{\al}(t)e_{\al} .
    \\
 \end{array}
 $$

 \sip

 Let $\{\nu_n: \odot^n V \rar V\}_{n\geq 1}$  and
 $\{\mu_{n;p}: (\otimes^n V)\otimes (\odot^p V)\rar V\}_{n\geq 2, p\geq 0}$
 be a $\cG erst_{\infty}$-algebra structure on $V$.
 Define the vector field,
 $$
 \vec{\nu}:= \sum_{n=1}^{\infty} \frac{1}{n!} \tau^{-1}\circ \nu_n(t, \ldots,t),
 $$
 and the tensors, for $n\geq 2$,
 $$
 \begin{array}{rccc}
   \mu_n: & \ot^n_{\f_{\cM} }\cT_{\cM}  & \lon & \cT_{\cM} \vspace{2mm} \\
    & X_1\ot \ldots \ot X_n & \lon  &
 \sum_{p=1}^{\infty} \frac{1}{p!} \tau^{-1}\circ \mu_{n;p}
 \left(\tau(X_1),\ldots, \tau(X_n);t, \ldots,t\right)
     \\
 \end{array}
 $$
 on $\cM$. Then the equations of Corollary 3.6.3 translate precisely into the
 statement that $(\vec{\nu},\mu_{\bullet})$ is a geometric $\cA_{\infty}$-structure
 on $\cM$.

\sip

 This argument also works in the opposite direction through the Taylor decomposition
 of all the tensors
 at the distinguished point.
 \hfill $\Box$

 \sip
 \bip


 {\bf 3.9. $\cG erst_{\infty}$-manifolds.}
 Since the operad  $\cG erst_{\infty}$
 is minimal, its algebras are strong homotopy ones, i.e.\ can be transferred via
 quasi-isomorphisms. In this subsection we essentially give a purely
 geometric description of the derived(=homotopy) category of
 $\cG erst_{\infty}$-algebras.

 \bip

 {\bf 3.9.1. Definition.} A  $\cG erst_{\infty}$-{\em manifold}\,
 is a smooth manifold $\cM$ together with a homotopy class (in the sense of
 2.10), $(\eth,[\mu_{\bullet}])$, of minimal
 geometric $\cA_{\infty}$-algebra structures on the tangent
 sheaf.

 \sip

 Note that we do {\em not}\, assume in the above definition that $\cM$ is affine:
 the notion of  $\cG erst_{\infty}$-{ manifold} is built on a collection
 of {\em tensors}\, satisfying a system of {\em diffeomorphism covariant}\,
  differential equations. The following  tautologically formulated
   theorem is one of the main results
  of this section; it essentially describes a functor from the category of
 $\cG erst_{\infty}$-algebras to a subcategory of the category of formal dg manifolds.

 \bip

 {\bf 3.9.2. Theorem.} {\em If the operad $\cG erst_{\infty}$ acts on a dg
 vector space $(V,d)$,
 then the formal graded manifold associated with the cohomology vector space $H(V,d)$ is
 canonically a  $\cG erst_{\infty}$-manifold.}

 \bip

 \Proof Since we work over a field in this paper, there always exists a
 quasi-isomorphism
 of complexes $(V,d) \rar (H(V,d),0)$ and hence an induced structure of
 $\cG erst_{\infty}$-algebra on $H(V,d)$. This structure, however, is not canonical.
 What we have to show is that, first reinterpreting this induced structure as a geometric
 $\cA_{\infty}$-structure $(\eth,\mu_{\bullet})$ on the pointed
 flat formal manifold $(\cM,*)=(H(V,d)[2],0)$, then passing to the associated homotopy class,
 $(\eth,[\mu_{\bullet}])$, just in the sense of $\cA_{\infty}$-algebras,
 and finally forgetting the flat structure one gets at the end the structure on $\cM$
 which does not depend on any choices made. It is precisely this structure which was
 termed in 3.9.1 an $\cG erst_{\infty}$-{manifold}.

 \sip

 By Theorem 3.8.1, our input is a formal affine   dg manifold,
 $$
 \left(M\simeq V[2],\, \vec{\nu}=\{\nu_{\bullet}\} \
 \mbox{with}\ \nu_1=d ,\, *=0\right),
 $$
 together with a geometric $\cA_{\infty}$-structure $\mu_{\bullet}$.
 Let us choose a cohomological splitting of the complex $(V,d)$,
 i.e.\ a decomposition of the $\Z$-graded vector space $V$ into a direct sum,
 $$
 V= H(V,d) \oplus B \oplus B[-1],
 $$
 in such a way that the differential vanishes when restricted to the summands
 $ H(V,D)\oplus B[-1]$ while on the remaining summand it equals the shifted by $[1]$
 identity map $B \rar B[-1]$. According to Kontsevich \cite{Ko1}, such a splitting
 can be lifted to an isomorphism of formal affine dg manifolds,
 $$
 (M,  \vec{\nu}, *) \simeq (\cM, {\eth}, *)\times (\cB, \vec{d}_{DR}, *),
 $$
 where
 \begin{itemize}
     \item  $(\cM,{\eth}, *)$ is the formal affine {\em minimal}\,
     dg manifold
     whose tangent space at $*$ is $H(V,d)[2]$,
     \item $(\cB, \vec{d}_{DR}, *)$ is the formal affine dg manifold
     whose tangent space at $*$ is $B[2] \oplus B[1]$, homological vector field
     $\vec{d}_{DR}$ is linear and coincides precisely with the usual De Rham differential
     when one identifies the structure sheaf, $\f_{\cB, *}$, with the
     De Rham algebra
     of smooth formal differential forms on the vector space $B$. In particular,
     both the cohomology groups, $H(\f_{\cB, *},\vec{d}_{DR})$
     and $H(\cT_{\cB, *}, Lie_{\vec{d}_{DR}})$ are trivial.
 \end{itemize}

 \sip

 Let $\pi_1: M \rar \cM$ and   $\pi_2: M \rar \cB$ be the projections associated with the
 chosen above cohomological splitting.
 There is an associated decomposition of complexes of vector spaces
 (note that differentials are
 {\em not}  $\f_M$-linear),
 $$
 (\cT_{M, *},  Lie_{\vec{\nu}}) = (\pi_1^*\cT_{\cM, *},  Lie_{\vec{\nu}})\
 \oplus \ (\pi_2^*\cT_{\cB, *},  Lie_{\vec{\nu}}).
 $$
 The tangent vector space at $*\in \cB$ can be identified with $B[2] \oplus B[1]$.
 Let $H: B[2] \oplus B[1] \rar B[2] \oplus B[1]$ be a degree $-1$ linear
 map which is equal to zero on the summand $B[2]$ and is equal to
 the shifted by $[-1]$  identity map $B[1]\rar  B[2]$ on the remaining summand.
  Denote by the same letter $H$ its natural $\f_{\cB, *}$-linear extension
  to $\cT_{\cB, *}$. It is an easy calculation to check that the identity
  automorphism of the tangent sheaf to $M$ decomposes as follows,
 $$
 Id = pr\ \oplus\  Lie_{\vec{\nu}}\circ \pi_2^*(H)\ \oplus\ \pi_2^*(H)\circ Lie_{\vec{\nu}},
 $$
 where $\pi_2^*(H)$ is assumed to act as zero on the
 summand $\pi_1^*\cT_{\cM, *}$,
 and $pr$ stands for the canonical projection $\cT_{M, *}\rar \pi_1^*\cT_{\cM, *}$.
 Thus we have constructed an $\f_M$-{\em linear}\ homotopy,
 $\pi_2^*(H):\cT_{M, *} \rar \cT_{M, *}$, associated with
  quasi-isomorphic complexes of $k$-{\em linear}\, vector spaces
  $(\cT_{M, *},  Lie_{\vec{\nu}})$ and
 $(\pi_1^*\cT_{\cM, *},  Lie_{\vec{\nu}})$.

 \sip

 Next step is to employ, say, the explicit formulae of \cite{Me1,KS},
 to construct an $\cA_{\infty}$-algebra structure,
 $$
 \left(\hat{\mu}_{n\geq 2}:\ot^n_{\f_{M} }\pi_1^*\cT_{\cM}\rar \pi_1^*\cT_{\cM},\
 \mu_1=Lie_{\vec{\nu}}\right),
 $$
 on the lifted tangent sheaf $\pi_1^*\cT_{\cM}$. The key fact that the resulting
 $\hat{\mu}_{n\geq 2}$
 are tensors is ensured by $\f_M$-linearity of the constructed homotopy.

 \sip

 Finally one repeats the above procedure using
 the standard contraction homotopy (which is $\pi^{-1}\f_{\cM}$-linear)
 of the cohomologically trivial De Rham complex $(\f_{\cB, *},\vec{d}_{DR})$
 to induce a {\em geometric}\, $\cA_{\infty}$-algebra
 structure, $(\eth,\tilde{\mu}_{\bullet})$,  on $\cM$. It is easy to check that the
 associated homotopy class (in the sense of
  2.11 for $\eth$,  and 2.10 for $\tilde{\mu}_{\bullet}$) does not depend on
  the choice of a particular factorization
  of the dg manifold $(M, \vec{\nu}, 0)$ into a direct product
  of a minimal dg manifold and a linearly contractible one.
 \hfill $\Box$

 \bip

 {\bf 3.9.3. Remark}. The special case of the above Theorem when
 $(V,d, \cG erst_{\infty}\mbox{-action})$ is just a (non-commutative)
  Gerstenhaber algebra was proved in \cite{Me2} by explicit perturbative
  calculations.

 \bip

 {\bf 3.10. Homotopy commutative sibling of  $\cG erst_{\infty}$}.
 Theorem~3.9.2 motivates the following definition.

 \bip

 {\bf 3.10.1. Definition}. $\cG erst_{\infty}^c$ is the operad whose algebras are
 given by
 $$
 \left\{\Ba{c}
 {\mathrm A}\  \cG erst^c_{\infty}\mbox{-}{\mathrm algebra}\\
 {\mathrm structure\  in }\ V
 \Ea
 \right\}
 :=
 \left\{\Ba{c}
 {\mathrm A\ geometric}\ \cC_{\infty}\mbox{-}{\mathrm structure\ on \ the\ pointed}\\
 {\mathrm  affine\  formal\ manifold}\ (V[2],0)
 \Ea
 \right\}.
 $$

 \sip

 Almost repeating 3.9.1 and 3.9.2  one obtains  the notion of
 $\cG erst_{\infty}^c$-manifold (which was in fact introduced earlier \cite{Me2}
 under the name {\em Frobenius}$_{\infty}$ manifold) and the statement:

 \bip

 {\bf 3.10.2. Theorem.} {\em If the operad $\cG erst_{\infty}^c$ acts on a dg
 vector space $(V,d)$,
 then the formal graded manifold associated with the cohomology vector space $H(V,d)$ is
 canonically a  $\cG erst_{\infty}^c$-manifold.}

 \bip

 \bip

 \bip
 \sip

 \begin{center}
 {\bf \S 4  Deformation theory}
 \end{center}

 \bip

 {\bf 4.1.  Deformation functor.}
 The traditional approach to the deformation theory  of a mathematical
 structure  ${\mathsf A }$  is based on the idea of
 deformation functor, $\Def_{\fg}$, on Artinian rings.

 \sip

 Initially that idea was applied to
 Artinian rings concentrated in degree $0$ so that the tangent space,
  ${\Def}_{\fg}(k[\var]/\var^2)$, to the deformation functor
 (which is the same as the Zarisski tangent space, $\cT_{\mathsf A} \cM$,  to the moduli space
 at the distinguished point)
 equals some particular homogeneous bit,  $H^i(\fg,d)$, of the $\Z$-graded
 cohomology group, $H^{\bullet}(\fg,d)=\oplus_{i\in \Z}H^i(\fg,d) $,
 of the dg Lie algebra, $(\fg,d)$,
 controlling the deformations of $\mathsf A$. The next homogeneous bit, $H^{i+1}(\fg,d)$,
 absorbs the obstructions
 to exponentiating infinitesimal deformations from $H^i(\fg,d)$ to genuine ones.

 \sip

 Recent studies in mirror symmetry  led to an extension \cite{BK,Ma,Me2} of
 the deformation functor first to $\Z$-graded and then differential
 $\Z$-graded  Artinian rings.
 The dg extension of ${\Def}_{\fg}$ always
 produces {\em smooth}\,  formal
 dg moduli spaces $(\cM, \eth, *)$ with Zariski tangent space, $\cT_{\mathsf A} \cM$,
  isomorphic to the full
 cohomology group $H^{\bullet}(\fg,d)$
 and with obstructions encoded into a homological vector field
 on $\cM$ \cite{Ko1,Me2}. The table below compares the two deformation functors
 in three important examples,

 \bip

 \begin{center}
 \begin{tabular}{|c|c||c|}
 \hline
 ${\mathsf A}$  &  $\cT_{\mathsf A} \cM$ in classical ${\mathsf Def}$  &
 $\cT_{\mathsf A}\cM$ in extended ${\mathsf Def}$ \\
 \hline
 Complex manifold  $(M,J)$ & $H^1(M,T_M)$ & $H^{\bullet}(M, \wedge^{\bullet}
 T_M)$\\
 Symplectic manifold $(M,\omega)$ & $H^2(M,\R)$ & $H^{\bullet}(M,\R)$ \\
 Associative algebra\ \  $(A,\circ)$ & ${\mathrm Hoch}^2(A,A)$ &
 ${\mathrm Hoch}^{\bullet}(A,A)$\\
 \hline
 \end{tabular}
 \end{center}

 \bip

 These developments lead naturally to a question:
  What happens to ${\mathsf A}$ when it is deformed
 in the generic direction in $H^\bullet(\fg,d)$ (rather than in $H^i(\fg,d)$)?
 Or, equivalently, what is the universal structure $\cA$ over the extended moduli space $\cM$?
 Thanks to Stasheff \cite{St}, we know the answer to this question
 in the case ${\mathsf A}={\mathrm Associative\ algebra}$:
  deforming
 any given associative algebra $A$ along a  generic tangent vector
 in $\cT_A \cM={\mathrm Hoch}^{\bullet}(A,A)$ one obtains, if
 all obstructions vanish,
 an $\cA_{\infty}$-algebra. To  author's knowledge, infinity versions of such
 notions as complex and symplectic structure are still a mystery, and the
 deformation theory in the form of deformation functor gives no clue to its solution.
 Moreover, the dg extension of $\Def$ makes it evident that
 one does no´t really need Artinian rings to do the
 deformation theory
 --- the dg versal moduli space, $(\cM, \eth, *)$,
  representing the functor $\Def_{\fg}$ on dg Artinian rings
 is nothing but the image of $(\fg,d)$ under
 the canonical functor
 $$
 \left\{ \Ba{c} \mathrm the\ category\ of\\ \mathrm dg\ Lie\ algebras
          \Ea\right\}
          {\lon}
 \left\{ \Ba{c} \mathrm the\ derived\ category\ of\\ \mathrm dg\ Lie\ algebras
          \Ea\right\}.
 $$
 Thus the Artinian functor approach to  deformation theory is tantamount  to
  a perturbative computation of the
 minimal $\caL_{\infty}$-model for the controlling dg Lie algebra.

 \sip


 \bip

 {\bf 4.2. An operadic guide to deformation theory}.
 We loose too much information about the deformed mathematical structure
 ${\mathsf A}$ if we naively understand the extended deformation theory as
 outlined above in Sect.\ 4.1, i.e.,
 as the deformation functor
 $\Def_{\fg}$ extended to (differential)
 $\Z$-graded  Artinian rings .
 Here is a suggestion on what one should do instead.

 \Bi
 \item[Step 1:] Associate to the mathematical structure ${\mathsf A}$ we wish to deform a
 ``controlling" {\sf Deformation Algebra}, $(\fg, [\, ,\, ], d,
 {\mathsf ADD})$,
 consisting of a dg Lie algebra\footnote{More generally, one can replace the dg Lie algebra
 structure on $\fg$ with a $\caL_{\infty}$-algebra structure.}
 $(\fg=\bigoplus_{i\in \Z} \fg^i,  [\, ,\, ], d)$,
 and a collection of some additional algebraic operations, ${\mathsf ADD}$, on $\fg$.
 For example,
 \Bi\item if ${\mathsf A}$ is a symplectic or complex structure (see Examples 4.4 and
 and 4.5 below) , then {\sf Deformation Algebra} is
 a graded commutative Gerstenhaber algebra, and ${\mathsf ADD}$ is just
 a graded commutative product consistent with the dg Lie algebra
 structure via Poisson type identities.
 \item if ${\mathsf A}$ is an associative algebra structure (see Example 4.6 below),
 then {\sf Deformation Algebra} is what
 is called in \cite{GV} a  {homotopy Gerstenhaber algebra}, and ${\mathsf ADD}$
 is an infinite series of operations called braces.
 \Ei

 \item[Step 2:]  Find a cofibrant resolution, $\caD\cA_{\infty}$, of the
 operad $\caD\cA$ describing species {\sf Deformation Algebra}
 obtained in the previous step. In many
 important cases
 there exists {\em the}\, minimal cofibrant model $\caD\cA_{\infty}$ of
 $\caD\cA$ whose differential
 is decomposable.

 \item[Step 3:] As the operad $\caD\cA_{\infty}$ is cofibrant,
 its algebras are strong homotopy ones, that is, $\cP$-algebra
 structures can be transferred by quasi-isomorphisms of complexes (see
 Theorem~2.8.5). Then
 choosing a cohomological
 splitting of $(\fg, d)$,
 one induces on the cohomology space $H(\fg, d)$ a
 canonical
 homotopy class of $\caD\cA_{\infty}$-algebras (which is independent of the splitting used.)

 \item[Step 4:] Try to find a geometric interpretation (called {\em extended moduli
 space}) of the  homotopy class of
 minimal $\caD\cA_{\infty}$-algebra
 structures canonically induced on $H(\fg,d)$ in the previous step, that is, try to
 interpret the latter in terms of
  sections of some natural vector bundles (equipped with $\eth$-connections,
  tensorial algebraic
  structures, \ldots ) over a
 formal pointed dg manifold $(\cM,\eth, *)$ whose tangent space at the
 unique geometric point is
 precisely $H(\fg,d)$.
 \Ei

 \sip

 {\bf 4.3. Motivation and evidence.} All the steps above are functorial with respect
 to the choice of input, {\sf Deformation Algebra}. Making a particular choice means
 essentially a choice of precision with which we want to do the deformation theory.
 The most crude one is to apply first the obvious forgetful functor,
 $$
 F: \left\{ \mathsf Deformation \ Algebra\right\} \lon
 \left\{\mathsf dg\ Lie\ Algebra\right\}
 $$
 and then apply the operadic algorithm. All the four steps can be easily
 fulfilled (see Sect.\ 2.11)
 with the final
 outcome at Step 4 being  a smooth formal minimal dg manifold $(\cM,\eth, *)$.
 Note that this outcome
 is well defined only up to an action of the formal diffeomorphism group $Diff(\cM,*)$,
 $$
 (\cM,\eth, *)\sim (\cM,f_*\eth, *), \ \ \ \forall f\in Diff(\cM,*).
 $$

 \sip

 Thus, whatever the choice of {\sf Deformation Algebra} in Step 1, the outcome
 of Step 4 will  always be at least a smooth dg manifold with its
 automorphism group $Diff(\cM,*)$. As we expect $\mathsf ADD$ in the input of the
 deformation theory
 to be consistent, in some or other sense, with the underlying dg Lie algebra structure,
 it is natural to expect that what we get from $\mathsf ADD$ in Step 4 will
 also be consistent with the homological vector field $\eth$ and behave reasonably well under
 the action of $Diff(\cM,*)$. The main results of this paper (see the Introduction
 for the list) show that these expectations are met in such important cases
 as, e.g.,  extended deformations of complex and symplectic structures.
 This is why we believe that Step 4, the most questionable  one
  in the above programme, makes  sense.

 \bip

{\bf 4.4. Example (deformations of complex structures).}
The dg commutative Gerstenhaber algebra controlling extended deformations of a given complex
structure on a real $2n$-dimensional manifold $X$ is given by
$$
\fg=\left(\bigoplus_{i=0}^{2n} \fg^i,\   \fg^i= \bigoplus_{p+q=i}
\Gamma(X, \wedge^p \cT_X\ot \Omega^{0,q}_{X}),\ [\ \bullet\ ],
\bp,\ \circ \right)
$$
where $\cT_X$ stands for the sheaf of holomorphic vector fields,
$\Omega^{s,q}_X$ for the sheaf of smooth differential forms of
type $(s,q)$,  and $[\ \bullet\ ]= {\mathsf Schouten\ brackets}\ot
{\mathsf wedge\ product\ of\ forms}$,
and $\circ= {\mathsf wedge\ product\ of\ polyvector\ fields}\ot
{\mathsf wedge\ product\ of\ forms}$.


\bip

{\bf 4.5. Example (deformations of Poisson and symplectic
structures).} The dg commutative Gerstenhaber algebra controlling deformations of a
given Poisson structure, $\nu\in \Gamma(X, \wedge^2 T_{\R})$, on
a real smooth manifold $X$ is given by
$$
\left(\oplus_{i=0}^{\dim X} \Gamma(X, \wedge^i T_{X}),
[\ \bullet\ ]={\mathsf Schouten\ brackets}, d=[\nu\bullet\ldots ],
{\mathsf wedge\ product\ of\ polyvector\ fields}
\right),
$$
where $T_{X}$ stands for the sheaf of real tangent vectors.

\sip

If $\nu$ is non-degenerate, that is, $\nu=\omega^{-1}$
for some symplectic form $\omega$ on $X$, then the natural
``lowering of indices map'' $\omega^{\wedge i}:
\wedge^i T_{X} \rar \Omega^i_{X}$ sends $[\nu\bullet\ldots ]$
into the usual de Rham differential. The image of the Schouten brackets
under this isomorphism we denote
by $[\ \bullet\ ]_{\omega}$. In this way we make the de Rham complex of
$X$
into a dg commutative Gerstenhaber algebra,
$$
\fg= \left(\bigoplus_{i=0}^{\dim X} \Gamma(X, \Omega^i_{X}),\
[\ \bullet\ ]_{\omega}, \
d={\mathsf de\ Rham\ differential},\
 {\mathsf wedge\ product\ of\ forms}\right),
$$
which controls the extended deformations of the symplectic
structure $\omega$. More explicitly,
\[
 [\ka_1\bullet \ka_2]_{\omega}:= (-1)^{\tlk_1}[i_{\nu}, d]
(\ka_1\wedge \ka_2) -
(-1)^{\tlk_1}\left([i_{\nu}, d]\ka_1\right)\wedge \ka_2 -
\ka_1\wedge[i_{\nu}, d]\ka_2, \ \ \ \forall \ka_1,\ka_2\in
\wedge^*T_{\R},
\]
with  $i_{\nu}:\Omega^{\bullet}_X\rar \Omega^{\bullet-2}_X $ being
the natural contraction with the 2-vector $\nu$.

\sip


\bip

{\bf 4.6. Example (deformations of holomorphic vector bundles)}. Let
$E\rar X$ be a holomorphic vector bundle on a complex manifold $X$.
The standard Lie algebra structure in the endomorphism sheaf $End(E)$ extends naturally
to $\wedge^\bullet_{\f_X} End(E)\simeq  \odot^\bullet_{\f_X} (End(E)[1])$ in such a way
that $(\wedge^\bullet_{\f_X} End(E), \wedge,
[\, \bullet \, ]_E)$ becomes a sheaf of graded commutative Gerstenhaber algebras.

\sip

The dg commutative Gerstenhaber algebra controlling extended deformations of a given
holomorphic structure in $E\rar X$ is given by
$$
\fg=\left( \Gamma(X, \wedge^\bullet_{\f_X} End(E)\ot
\Omega^{0,\bullet}_{X}),\ [\ \bullet\ ], \bp,\ \circ \right)
$$
where  $[\, \bullet\, ]= [\, \bullet \,  ]_E\ot
{\mathsf wedge\ product\ of\ forms}$,
and $\circ= {\mathsf wedge\ product\ of\ endomorphisms}\ot
{\mathsf wedge\ product\ of\ forms}$.

\bip

{\bf 4.7. Homotopy Gerstenhaber algebras}.
Let $V$ be a graded vector space and $({B}V:=\ot^{\bullet\geq 0}V[1],\Delta)$
 the free tensor coalgebra cogenerated by $V[1]$.

\bip

{\bf 4.7.1. Definitions}. (i) A {\em $B_{\infty}$-algebra}\, structure on a graded
vector space $V$ is the structure of dg bialgebra,
$$
\left(BV, \Delta, \circ: BV\ot BV \rar BV , d:BV\rar BV\right),
$$
on the tensor coalgebra $(BV, \Delta)$ such that the element $1\in k=\ot^0V[1]$
is the identity element.

\sip

(ii) A {\em homotopy Gerstenhaber algebra}\, structure on a graded
vector space $V$ is a structure of $B_{\infty}$-algebra such that multiplication
$\circ$ preserves the filtration $F_r:= \ot^{\bullet\leq r}V[1]$.
We denote by $h\cG$ the operad whose algebras are homotopy Gerstenhaber
algebras.

\sip

The $h\cG$-algebra structure on $V$ can be described by a collection of
homogeneous maps,
\Beqrn
M_k: && \ot^k V \lon V, \\
M_{1;k}: && \ot^{k+1} V \lon V,
\Eeqrn
satisfying a system of quadratic equations
 written explicitly in \cite{V}. In particular, the operations $M_2$ and $M_{1;1}$
induce on $H(V,M_1)$ the structure of graded commutative Gerstenhaber algebra.
Thus there is a canonical map of operads, $p: h\cG \rar \cG$.

\bip

{\bf 4.7.2. Example (higher order Steenrod operations)}.
Let  $S_\bullet X$ be the singular chain complex,
of a topological space $X$. Elements of $S_n X$ are  formal linear combinations
of continuous maps, $\sigma: \Delta[n]\rar X$, from the standard $n$-simplex
$\Delta[n]$ to $X$. For such a map $\sigma\in S^n X$ and a
$k$-face, $f:\Delta[k] \rar \Delta[n]$, of the standard simplex spanned by vertices
$n_0=f(0), \ldots, n_k=f(k)$ ($f$ being injective and monotone),
denote by $\sigma[n_0,\ldots,n_k]\in S_kX$ the associated composition,
$$
\sigma[n_0,\ldots,n_k]: \Delta[k] \stackrel{f}{\lon} \Delta[n]
\stackrel{\sigma}{\lon} X.
$$

\sip

Let $S^\bullet X:=\Hom_k(S_\bullet, k)$ be the associated cochain complex of $X$.
It was noted by Gerstenhaber and Voronov in \cite{GV} using earlier results
of Baues that the data,
\Beqrn
M_1(\phi)(\sigma)&:=& \sum_{k=0}^{n+1}(-1)^k\phi(\sigma[0,1, \ldots, \hat{k}, \ldots,n]),\\
M_2(\phi,\psi)(\sigma) &:=& \sum_{k=0}^n \phi(\sigma[0,\ldots,k])
\psi(\sigma[k,\ldots,n]), \ \ \forall\ \phi,\psi\in S^\bullet X,\\
M_k &:=& 0\ \mbox{for}\ k\geq 3,\\
M_{1;k}(\phi_0; \phi_1,\ldots,\phi_n)(\sigma) &:=&
\phi_0(\sigma[0,n_1,n_1+n_2, \ldots, n_1+\ldots + n_k])\\
&&\phi_1(\sigma[0,1,\ldots, n_1)\phi_2(\sigma[n_1,\ldots, n_1+n_2)\dots\\
&&\phi_k(\sigma[n_1\ldots +n_{k-1}, \ldots, n_1+\ldots + n_k])\\
&&
\forall\ \phi_0\in S^k X, \phi_1\in S^{n_1}X, \ldots, \phi_k\in S^{n_k}X,
\Eeqrn
make $S^\bullet X$ into an $h\cG$-algebra. The operation $M_{1;1}$ is nothing but
the Steenrod operation $\cup_1$.

\bip

{\bf 4.7.3.  Example (deformations of associative algebras)}.
Let $A$ be an associative algebra, and $C^{\bullet}(A,A):=\Hom_k(A^\bullet, A)$
its {\em Hochschild complex}\, with the differential,
\Beqrn
(d\phi)(a_1, \ldots, a_{n+1})&:=& a_1\phi(a_2, \ldots,a_{n+1})\\
&& \sum_{k=1}^n (-1)^k \phi(a_1, \ldots, a_{k-1},a_k,a_{k+1}, \ldots,a_{n+1})\\
&& + (-1)^{k+1}\phi(a_1,a_2, \ldots,a_n)a_{n+1}.
\Eeqrn
It was shown in \cite{GV} that the data,
\Beqrn
M_1&:=& d,\\
M_2(\phi,\psi)(a_1,\ldots,a_{k+l}) &:=&\phi(a_1,\ldots,a_k)
\psi(a_{k+1},\ldots, a_{k+l}) \ \ \forall\ \phi\in C^k(A,A), \psi\in
C^l(A,A)\\
M_k &:=& 0\ \mbox{for}\ k\geq 3,
\Eeqrn
and, for any  $\phi_0,\phi_1, \ldots, \phi_k\in C^{\bullet}(A,A)$,
$$
M_{1;k}(\phi_0; \phi_1,\ldots,\phi_n)(a_1,\ldots, a_m) :=
\hspace{9cm}
$$
$$
\hspace{3cm}\sum (-1)^{\sum_{p=1}^n(|\phi_{p}|-1)i_p} \phi_0(a_1, \ldots, a_{i_1},
\phi_1(a_{i_1}, \ldots),\dots, a_{i_n},
\phi_n(a_{i_n},\ldots ),\ldots, a_m)
$$
where the summation runs over all possible ordered substitutions of $\phi_1,\ldots,\phi_n$
into $\phi_0$, makes the Hochschild complex into an $h\cG$-algebra.

\sip

This $h\cG$-algebra controls deformations of the associative algebra structure in $A$.

\bip

{\bf 4.8. Deligne's conjecture.} The operations $M_2$ and $M_{1;1}$ induce on the Hochschild
cohomology, $H^\bullet (A,A)$, the structure of $\cG$-algebra. Deligne conjectured
that this action of the operad $\cG$ on $H^\bullet (A,A)$ can be lifted
to the action of $\cG_{\infty}$ on $C^\bullet (A,A)$. This conjecture, which was recently
proved in \cite{Ko3,KS0,MS,Ta,V}, is essentially the same as the following statement.

\bip

{\bf 4.8.1. Theorem} \cite{Ta,TT,V}. {\em There is a natural morphism of operads},
$f:\cG_{\infty}\rar h\cG$, {\em such that the diagram},
$$
 \xymatrix{
 \cG_{\infty} \ar[r]^f  &  h\cG \ar[d]^p \\
   & \cG \ar@{<-}[ul]^{q.-iso.}
 }
$$
{\em commutes.}

\bip

{\bf 4.9. Approximations to deformation theory.} The cofibrant
 resolution, $\caD\cA_{\infty}$, of the
 operad $\caD\cA$ describing species {\sf Deformation Algebra}
could be so complicated  that it would be unrealistic to ask for an
immediate geometric
interpretation of the homotopy class of minimal $\caD\cA_{\infty}$-algebra
structures as in Step 4.

\sip

In Steps 3 and 4 one can therefore replace $\caD\cA_{\infty}$ by its
``approximation", a cofibrant operad
$\widetilde{\caD\cA}_{\infty}$ fitting the commutative diagram
$$
 \xymatrix{
 \widetilde{\caD\cA}_{\infty} \ar[r]^i  &  \caD\cA_{\infty} \ar[d]^{q.-iso.} \\
   & \cG \ar@{<-}[ul]^j
 }
$$
for some natural morphisms of operads $i$ and $j$, $i$ being preferably a cofibration.

\bip


{\bf 4.10. Approximations to $\cG_{\infty}$ deformation theory}. As discussed in
Sections~4.4-4.6,
extended deformations of basic geometric structures are described by the operad
$\cG_{\infty}$. In view of theorem 4.8.1, the same operad can be applied to the
deformation theory of associative algebras.

\sip

At present we have no complete picture of the geometric object behind a homotopy
class of minimal $\cG_{\infty}$-algebras. We appeal instead to the infinite
tower of cofibrant approximations to $\cG_{\infty}$ introduced in Sect.\ 2.7,
 $$
 \caL_{\infty}=\cG_{\infty}^{(1)} \lon  \cG_{\infty}^{(2)} \lon \cG_{\infty}^{(3)}
  \lon \ldots \lon \cG_{\infty}^{(n)}\lon  \ldots \lon \cG_{\infty},
 $$
and, in the rest of this paper, attempt to give such a picture for the first three
floors of this tower using all the previous results.
The ground floor, $\cG_{\infty}^{(1)}$, corresponds, in view of 2.7.2,
to the forgetful functor 4.3 so that $\cG_{\infty}^{(1)}$-approximation to the
deformation theory of examples 4.4-4.6 and  4.7.3 simply says that the extended
moduli space is a formal minimal dg manifold $(\cM, \eth, *)$. If its non-linear
cohomology (see Introduction),
$$
 M\simeq \frac{{\mathsf Zeros}(\eth)}{\Img \eth},
 $$
 makes sense in some geometric category,
it is precisely the  versal moduli space associated with the
$\Z$-graded extension of the classical deformation functor $\Def_\fg$.


\bip

{\bf 4.11. Theorem}. {\em There is a canonical isomorphism of
operads,
$$
\cG_{\infty}^{(2)}=\cG erst_{\infty}^{c},
$$
 where $\cG erst_{\infty}^{c}$ is defined in Sect.\ 3.10.1.}

 \sip

 \Proof The statement follows immediately from the definition 3.10.1 and the proof
 of Theorem 3.4.2 in \cite{Me2} (see also Sect.\ 4.12 for a reconstruction
 of that argument).
 \hfill $\Box$

\bip

{\bf 4.11.1.  Corollary}.
{\em If the operad $\cG_{\infty}$ acts on a dg vector space $(V,d)$,
 then the formal graded manifold associated with the cohomology vector space $H(V,d)$
 is canonically a  $\cG erst_{\infty}^{c}$-manifold.}

 \sip

 \Proof The statement follows from 4.11  and 3.9.2.
 \hfill $\Box$

\bip

{\bf 4.11.2. Corollary} \cite{Me2}. {\em Let $\mathsf A$ be one of the structures 4.4-4.6 or 4.7.3.
The extended  moduli space of deformations of $\mathsf A$ is naturally
 a $\cG erst_{\infty}^{c}$-manifold, i.e.\ a dg manifold equipped with a homotopy class
 of geometric $C_{\infty}$-structures on the tangent bundle}.

\bip

{\bf 4.11.3. Notation}. For a graded module $V$ over a graded ring $\f$ we set
$$
\circledast^\bullet_\f V := \left(\ot^\bullet_\f V[1]\right)[-1]
$$
and
$$
\frac{\circledast^{\bullet}_\f V }{\mathrm shuffle\ products}
:=  \left(\frac{\ot^\bullet_\f V[1]}{\mathrm shuffle\ products}\right)[-1].
$$

\bip
{\bf 4.12. Proof of Theorems A, B and E}.
The $\cG_{\infty}$-algebra structure on $V$ induces canonically
a homotopy class of minimal $\cG_{\infty}$-algebra structures on its cohomology,
$H(V,d)$. Let ${\mathbf \Theta}$ be any representative of this homotopy class.

\sip

Let $\{\p_{\al}, \al=1,\ldots, \dim H(V,d)\}$ be a homogeneous basis of $H(V,d)[2]$ and
 $\{t^{\al}\}$ the associated dual basis of $H(V,d)^*[-2]$.
 We identify the latter  with coordinate functions on the formal manifold $(\cM,*)$
 associated with $H(V,d)[2]$, and the former with basis vector fields $\p/ \p t^{\al}$.
  We can assume without loss of generality\footnote{For example, we could opt to work
  in  the category
of graded manifolds over  graded base spaces, ``sources of $\Z$-graded constants";
then all the signs lost under our assumption in ``natural over the base" calculations
can be easily restored
through the condition that the expression under study
 is functorial with respect to the base space change.} that
  degrees of all $\p_{\al}$ vanish $\bmod 2\Z$.

\sip

Let $t^{B_i}$, $i\geq 2$, stand for the Lie polynomial,
$$
t^{B_i}:= \left[\ldots \left[[t^{\beta_1}\bullet t^{\beta_2}]\bullet t^{\beta_3}\right]
\ldots \bullet t^{\beta_i} \right],
$$
and
$$
\p_{B_i}\equiv (\p_{\beta_1}| \p_{\beta_2})
\ldots | \p_{\beta_i})
$$
for the image of the tensor product,
$$
\p_{\beta_1}\ot \p_{\beta_2}\ot
\ldots \ot \p_{\beta_i},
$$
under the degree $n-1$ composition
$$
\ot^{n} H(V,d)[2]\lon \circledast^{n} H(V,d)[2].
$$
Note that under our assumption the parities of $t^{B_i}$ and $\p_{B_i}$ are
both equal to $(i-1)\bmod 2\Z$. Here and below we use ${B_i}$ to denote
the multi-index $\be_1\be_2\ldots\be_i$.

\sip

By Proposition~2.6.1, the $\cG_{\infty}$-algebra  structure $\Theta$
is the same as the differential, $\delta$, of the Gerstenhaber algebra
$\hat{\odot}^{\bullet} {\mathsf Lie}({H(V,d)}^*[-2])$. The latter is uniquely
determined
on the generators,
$$
\delta t^{\al}=\sum_{k\geq 0, n\geq 0
\atop {B_{i_1}, \ldots, B_{i_k}\atop \be_1\ldots \be_n}} \frac{1}{n!}
\Theta_{B_{i_1}\ldots B_{i_k};\be_1\ldots\be_n}^\al
t^{B_{i_1}}\ldots t^{B_{i_k}} t^{\be_1} \ldots t^{\be_n},
$$
for some homogeneous constants $\Theta_{B_{i_1}\ldots
B_{i_k};\be_1\ldots\be_n}^\al\in k$ (see footnote 9). Here and
below juxtaposition of $t^\bullet$s means their symmetric product
$\odot$.

We re-arrange the above data into a smooth ``vector field",
$$
\eth := \sum_{n\geq 0
\atop \be_1\ldots \be_n} \frac{1}{n!}
\Theta_{0;\be_1\ldots\be_n}^\al
 t^{\be_1} \ldots t^{\be_n}\frac{\p}{\p t^\al}
$$
and a collection of maps, for $k\geq 2$,
 $$
 \begin{array}{rccc}
   \Theta^{[k]}: & \odot^{k}_{\Bbb C}
   \left({\circledast^{\bullet}_{\Bbb C}
    \cT_{\cM} }\right)
   & \lon &  \cT_{\cM} \vspace{4mm} \\
    & \p_{B_1}\odot \ldots \odot \p_{B_k} & \lon  &
    \sum_{\atop { n\geq 0
\atop \be_1\ldots \be_n}} \frac{1}{n!}
\Theta_{B_{i_1}\ldots B_{i_k};\be_1\ldots\be_n}^\al
 t^{\be_1} \ldots t^{\be_n}\frac{\p}{\p t^\al}.
     \\
 \end{array}
 $$

\sip

\noindent
Thus the defining equation for the differential $\delta$ takes the form,
$$
\sum_{\al}(\delta t^\al) \p_\al = \eth + \sum_{k\geq 1
\atop B_{i_1}, \ldots, B_{i_k}}
\Theta^{[k]}( \p_{B_1},\ldots , \p_{B_k})
t^{B_{i_1}}\ldots t^{B_{i_k}}.
$$
The ideal $I:= <{\mathsf Lie}^{\geq 2}{V}^*[-2]>$ is generated by $t^{B_{i}}$.

\sip

Using the identities,
$$
\delta f(t) = \sum_{\al}(\delta t^\al) \frac{\p f}{\p t^\al},
$$
and
$$
[t^\al\bullet f(t)] = \sum_{\be}[t^\al\bullet t^\be]\frac{\p f}{\p t^\be},
$$
where  $f(t)$ is an arbitrary smooth function  on $(\cM,*)$,
it is not hard to study the transformation properties of the defined above fields
$\eth$ and $\Theta^{[k]}$ under an arbirary formal (non-linear) change of coordinates,
$$
t^\al \lon \tilde{t}^\al=f^\al(t),
$$
and conclude that
\begin{itemize}
    \item[(i)] $\eth$ is indeed a vector field on $\cM$;
    \item[(ii)] the map $\Theta^{[1]}$ factors through the composition
$$
\Theta^{[1]}:
   {\circledast^{\bullet}_{\Bbb C}
    \cT_{\cM} } 
   \lon
 {\circledast^{\bullet}_{\f_{\cM} }
    \cT_{\cM} } 
    \lon
    \cT_{\cM},
$$
i.e.\ represents a family of {\em tensors}, $\mu_\bullet: \ot^{\bullet}_{\f_{\cM}}
    \cT_{\cM}\rar \cT_{\cM}$, vanishing on shuffle products;
    \item[(iii)] the maps $\Theta^{[k]}$ are sections of certain jet bundles (of order $k-1$)
    on $\cM$.
\end{itemize}

\sip

\noindent A similar calculation shows
that
\begin{itemize}
    \item[(iv)] the equation $\sum_{\al}(\delta^2 t^\al) \p_\al=0 \bmod I$
    implies $[\eth, \eth]=0$;
    this is essentially Lemma~2.7.2;

    \item[(v)] the equation $\sum_{\al}(\delta^2 t^\al) \p_\al=0 \bmod I^2$ implies
    that the data $\Theta^{[1]}\simeq \mu_\bullet$ is nothing but a geometric
    $\cC_\infty$-structure on the dg manifold $(\cM,*, \eth)$; this is essentially Theorem~4.11.
\end{itemize}
To prove the statements A,B and E we have to move one level up and  study
the equation
$$
\sum_{\al}(\delta^2 t^\al) \p_\al=0 \bmod I^3,
$$
which now involves the non-tensorial object $\Theta^{[2]}$, a section of the bundle
$J^1(\ot^{\bullet}_{\f_{\cM}} \cT_{\cM}^*)\ot_{\f_{\cM}} (\ot^{\bullet}_{\f_{\cM}} \cT_{\cM}^*)
\ot_{\f_{\cM}} \cT_{\cM}$. It is convenient for our purposes
to choose a splitting of this bundle,
say, an affine torsion-free affine connection $\nabla$ on $\cM$. Then we can
replace non-tensorial objects
$$
t^{B_n}= \left[\ldots \left[[t^{\beta_1}\bullet t^{\beta_2}]\bullet t^{\beta_3}\right]
\ldots \bullet t^{\beta_n} \right],
$$
by tensorial ones,
$$
\breve{t}^{B_n}= \left[\ldots \left[[t^{\beta_1}\bullet t^{\beta_2}]'\bullet t^{\beta_3}\right]'
\ldots \bullet t^{\beta_n} \right]',
$$
which are given recursively by
\Beqrn
[t^{\beta_1}\bullet t^{\beta_2}]'&=& [t^{\beta_1}\bullet t^{\beta_2}] \\
\left[[t^{\beta_1}\bullet t^{\beta_2}]'\bullet t^{\beta_3}\right]'
&=& \left[[t^{\beta_1}\bullet t^{\beta_2}]'\bullet
t^{\beta_3}\right] + \Gamma^{\be_1}_{\mu\nu}\odot [t^{\mu}\bullet
t^{\beta_2}]'\odot [t^{\nu}\bullet t^{\beta_3}]' \\
&& + \Gamma^{\be_2}_{\mu\nu}\odot [t^{\mu}\bullet
t^{\beta_1}]'\odot [t^{\nu}\bullet t^{\beta_3}]'
\\
\ldots &=& \ldots \\
\breve{t}^{B_n} &=& \left[\breve{t}^{B_{n-1}} \bullet
t^{\beta_n}\right] +
\sum_{i=1}^{n-1}\Gamma^{\be_i}_{\mu\nu}\odot\breve{t}^{B_{n-1}(\beta_i,\mu)}\odot
[t^{\nu}\bullet t^{\beta_n}]', \Eeqrn where $B_{n-1}(\beta_i,\mu)$
stands for the multi-index $\be_1\ldots \be_{i-1}\mu
\be_{i+1}\ldots \be_{n-1}$ and $\Gamma^{\be}_{\mu\nu}$ are the
Christoffel symbols of the connection $\nabla$ in the coordinate
system $\{t^\al\}$. Thus the defining equation for the
differential $\delta$ can be written in the form,
$$
\sum_{\al}(\delta t^\al) \p_\al = \eth + \sum_{k\geq 1 \atop
B_{i_1}, \ldots, B_{i_k}}\breve{\Theta}^{[k]}( \p_{B_1},\ldots ,
\p_{B_k})\odot \breve{t}^{B_{i_1}}\odot\ldots\odot
\breve{t}^{B_{i_k}}.
$$
where all the coefficients, $\breve{\Theta^{[k]}}$, are tensors rather than
sections of the jet bundles.

\sip

Now, to prove
Theorems A and B it is enough to understand the first non-trivial component,
$$
[t^{\beta_1}\bullet t^{\beta_2}]\odot [t^{\ga_1}\bullet t^{\ga_2}],
$$
 of the formal power series
$\sum_{\al}(\delta^2 t^\al) \p_\al \bmod I^3$. The only terms of
 $\sum_{\al}(\delta t^\al) \p_\al \bmod I^3$ which contribute to that component
 are written explicitly below,
\Beqrn \sum_{\al}(\delta t^\al) \p_\al &=& \eth + \sum
\breve{\Theta^{[1]}}(\p_{\be_1}| \p_{\be_2})\odot
[t^{\beta_1}\bullet t^{\beta_2}] + \sum
\breve{\Theta}^{[1]}\left(\p_{\be_1}|
\p_{\be_2}|\p_{\be_3}\right)\odot
\left[[t^{\beta_1}\bullet t^{\beta_2}]\bullet t^{\beta_3}\right]' \\
&& +\, \sum \breve{\Theta}^{[2]}\left((\p_{\be_1}|
\p_{\be_2}),(\p_{\ga_1}| \p_{\ga_2})\right) \odot
[t^{\beta_1}\bullet t^{\beta_2}]\odot [t^{\ga_1}\bullet t^{\ga_2}]
+ \ldots \Eeqrn Applying to the shown terms the differential
$\delta$ and ignoring all components of the equation
$\sum_{\al}(\delta^2 t^\al) \p_\al =0 \bmod I^3$ except the chosen
one, one gets an equation,
$$
\hspace{57mm} [\mu_2,\mu_2]^{\nabla} = Lie_\eth \mu_{2,2}, \hspace{57mm} (*)
$$
where
\Beqrn
[\mu_2,\mu_2]^{\nabla}(X,Y,Z,W)&=&
[\mu_2,\mu_2]_{HM}(X,Y,Z,W) \\
&&
- (Lie_\eth \mu_{3})(W,\nabla_Z X + \nabla_X Z, Y)
- (Lie_\eth \mu_{3})(Z,\nabla_W X + \nabla_X W, Y)\\
&&
- (Lie_\eth \mu_{3})(W,\nabla_Z Y + \nabla_Y Z, X)
- (Lie_\eth \mu_{3})(Z,\nabla_W Y + \nabla_Y W, X)\\
&&
+ \mu_{3}(W, (Lie_\eth \nabla)_Z X + (Lie_\eth \nabla)_X Z, Y)\\
&&
+ \mu_{3}(W, (Lie_\eth \nabla)_Z X + (Lie_\eth \nabla)_X Z, Y)\\
&&
+ \mu_{3}(W, (Lie_\eth \nabla)_Z Y + (Lie_\eth \nabla)_Y Z, X)\\
&&
+ \mu_{3}(W, (Lie_\eth \nabla)_Z Y + (Lie_\eth \nabla)_Y Z, X),
\Eeqrn
$\mu_2$, $\mu_3$ and $\mu_{2,2}$ are, respectively, the compositions,
$$
\mu_2: \odot^2_{\f_M}\cT_\cM \stackrel{\simeq}{\lon}
 \frac{\circledast^{2}_{\f_{\cM} }  \cT_{\cM} }{\mathrm shuffle\ products}
 \hook  \frac{\circledast^{\bullet}_{\f_{\cM} }  \cT_{\cM} }{\mathrm shuffle\ products}
 \stackrel{\Theta^{[1]}}{\lon} \cT_\cM,
$$
$$
\mu_3:  \frac{\ot^{3}_{\f_{\cM} }  \cT_{\cM} }{\mathrm shuffle\ products}
\stackrel{\simeq}{\lon} \frac{\circledast^{3}_{\f_{\cM} }  \cT_{\cM} }
{\mathrm shuffle\ products}
 \hook  \frac{\circledast^{\bullet}_{\f_{\cM} }  \cT_{\cM} }{\mathrm shuffle\ products}
 \stackrel{\Theta^{[1]}}{\lon} \cT_\cM
$$
$$
\mu_{2,2}: \wedge^2_{\f_\cM}(\odot^2_{\f_M}\cT_\cM)
\stackrel{\simeq}{\lon}
\odot^2_{\f_M} \frac{\circledast^{2}_{\f_{\cM} }  \cT_{\cM} }{\mathrm shuffle\ products}
 \hook  \frac{\circledast^{\bullet}_{\f_{\cM} }  \cT_{\cM} }{\mathrm shuffle\ products}
 \stackrel{\breve{\Theta}^{[2]}}{\lon} \cT_\cM,
$$
$[\mu_2,\mu_2]_{HM}$ is the Hertling-Manin bracket (see Sect.\ 0.1),
and $X,Y,Z,W$ are arbitrary smooth even (for simplicity of presentation)
vector fields on $\cM$.
As the r.h.s. of the  equation $(*)$ is obviously
a tensor, the l.h.s.\ must represent a tensor as well.
(The original Hertling-Manin bracket can not possibly  be a tensor without correction
terms as the product $\mu_2$ is
associative only up to homotopy.)
Changing the affine connection,
$$
\nabla  \rar \nabla' = \nabla + S,
$$
for some symmetric  tensor $S: \odot^2_{\f_\cM}\cT_\cM\rar \cT_\cM$,
we get
$$
[\mu_2,\mu_2]^{\nabla'}= [\mu_2,\mu_2]^{\nabla} -
 2 Lie_\eth \nu,
$$
where
\Beqrn \nu(X,Y,Z,W)&=& \mu_{3}(W,S(Z,X),Y) +
\mu_{3}(Z,S(W,X),Y)\\
&&+ \mu_{3}(W,S(Z,Y),X) +
\mu_{3}(Z,S(W,Y),X).
\Eeqrn
Thus the cohomology class,
$[[\mu_2,\mu_2]]\in  H(\wedge^2_{\f_\cM}(\odot^2_{\f_M}\cT_\cM))$,
associated with $[\mu_2,\mu_2]^{\nabla}$ does not depend on the choice of $\nabla$.
 The equation $(*)$ implies Theorems A and B almost
immediately.

\sip

To prove Theorem E one needs to understand  a general structure of  the components,
$$
\breve{t}^{B_{k}}\odot\breve{t}^{B_{l}},
$$
 of the formal power series
$\sum_{\al}(\delta^2 t^\al) \p_\al \bmod I^3$.
Here are  the terms of
 $\sum_{\al}(\delta t^\al) \p_\al \bmod I^3$ contributing to that component,
$$
\sum_{\al}(\delta t^\al) \p_\al = \eth + \sum_{k\geq 2,
B_{k}}
\breve{\Theta}^{[1]}( \p_{B_k})\odot
\breve{t}^{B_{k}} +
 \sum_{k,l\geq 2\atop B_{k}, B_{l}}
\breve{\Theta}^{[2]}( \p_{B_k}, \p_{B_l})\odot
\check{t}^{B_{k}} \odot\breve{t}^{B_{l}} + \ldots\ .
$$
It is not hard  to see that the projection of
the equation $\sum_{\al}(\delta^2 t^\al) \p_\al =0
\bmod I^3$ to that component gives the equation of the form,
$$
[\mu_k(X_1,\ldots,X_k), \mu_l(Y_1,\ldots,Y_k)] + {\mathrm correction\ terms}
= (Lie_\eth \mu_{k,l})(X_1,\ldots,X_k,Y_1,\ldots,Y_k),
$$
where the brackets stand for the standard commutator of vector fields, and
 $\mu_k,\mu_l$ ($\mu_{k,l}$) are homogeneous components of
$\Theta^{[1]}$ (respectively, $\breve{\Theta}^{[2]}$). Again,
as the r.h.s. of the above equation is obviously
a tensor, the l.h.s.\ must be a tensor as well. This tensor,
$[\mu_k,\mu_l]^\nabla$, is $Lie_\eth$-closed and, moreover,
defines a vanishing cohomology class
$[[\mu_k,\mu_l]]$ of the complex of
 $\f_M$-modules, $(\ot^{\bullet}_{\f_{\cM}}\cT^*_{\cM}\ot \cT_{\cM},
 Lie_{\eth})$.
 \hfill $\Box$

\bip

{\small
{\em Acknowledgement.} A part of this work was done during author's visit
to the Max Planck Institute for Mathematics in Bonn
which was supported by the ``Sofia Kovalevskaya'' award of Matilde Marcolli.  Excellent working
conditions in the MPIM are gratefully acknowledged.
I would like
to thank E.\ Getzler, Yu.I.\ Manin, M.\ Markl and A.A.\ Voronov for  valuable
correspondence and discussions, and M.\ Marcolli for supporting the project and the visit.
Thanks are also due to the referee for valuable remarks.}

  {\small

  \end{document}